# Stable transport of information near essentially unstable localized structures

THIERRY GALLAY[1], GUIDO SCHNEIDER[2], HANNES UECKER[2]

[1] Institut Fourier, Université de Grenoble I, F - 38402 Saint-Martin d'Hères, France
[2] Mathematisches Institut I, Universität Karlsruhe, D - 76128 Karlsruhe, Germany

June 23, 2003

## Abstract

When the steady states at infinity become unstable through a pattern forming bifurcation, a travelling wave may bifurcate into a modulated front which is time-periodic in a moving frame. This scenario has been studied by B. Sandstede and A. Scheel for a class of reaction-diffusion systems on the real line. Under general assumptions, they showed that the modulated fronts exist and are spectrally stable near the bifurcation point. Here we consider a model problem for which we can prove the nonlinear stability of these solutions with respect to small localized perturbations. This result does not follow from the spectral stability, because the linearized operator around the modulated front has essential spectrum up to the imaginary axis. The analysis is illustrated by numerical simulations.

# Contents

# 1 Introduction

Localized structures such as pulses and fronts play an important role in the mathematical theory of information transport. Typical situations where such nonlinear phenomena arise are the propagation of electromagnetic waves in wires or fibers [AA97, NM92], and the motion of electric pulses along nerve axons [Hux52]. An important issue, both from a theoretical and a practical point of view, is the robustness of these solutions with respect to small inhomogeneities of the propagation medium.

In a remarkable paper [SS99], B. Sandstede and A. Scheel studied a new bifurcation scenario for traveling pulses in reaction-diffusion systems on the real line. They investigated the situation where the homogeneous steady state at infinity becomes unstable and bifurcates to a spatially periodic Turing pattern. The originally stable pulse thus undergoes an "essential instability", in the sense that the essential spectrum of the linearized operator crosses the imaginary axis at the bifurcation point. Under general assumptions, the authors showed that the original pulse bifurcates to a "modulated pulse" which is time-periodic in a uniformly translating frame. They also proved that this bifurcating solution is spectrally stable [SS00]. However, since the spectrum of the linearization extends all the way to the imaginary axis (without gap), this last result does not immediately imply the nonlinear stability of the modulated pulse. The analysis of [SS99, SS00] can be generalized to front solutions connecting two different stable equilibria [SS01a, SS01b]. In this case, modulated fronts may bifurcate from an existing traveling wave if one or both of the rest states at infinity become unstable.

In this paper, we go beyond the linear stability analysis of [SS00, SS01b] and we show,

at least on a specific example, that modulated fronts are *nonlinearly* stable with respect to spatially localized perturbations. In simple terms, this result implies that information can be transported in a stable manner even if the propagation medium becomes unstable through a Turing bifurcation. To keep the analysis as simple as possible, we do not consider abstract reaction-diffusion systems as in [SS99], but we prefer to concentrate on a model problem that exhibits all features of the general case. Although this has not been proved so far, we certainly expect that all results below hold true for general reaction-diffusion systems under the same assumptions as in [SS00] (for pulses) or [SS01b] (for fronts).

Our model problem is a Chaffee-Infante equation for the first variable $u$ coupled to a Swift-Hohenberg equation for the second variable $v$, namely:

$$\begin{aligned} \partial_t u &= \partial_x^2 u + \tfrac{1}{2}(u - c_0)(1 - u^2) + v \, , \\ \partial_t v &= -(1 + \partial_x^2)^2 v + \alpha v - v^3 - \gamma v F(u) \, , \end{aligned} \quad (1)$$

where $u(x,t), v(x,t) \in \mathbb{R}$, $x \in \mathbb{R}$, and $t \geq 0$. This system is especially convenient to analyze, because it couples two scalar equations which are rather well understood. In what follows, the speed parameter $c_0$ and the coupling parameter $\gamma$ will be fixed, with $0 < c_0 < 1$ and $\gamma > 0$ not too big (see Theorem 2.3 below). Our bifurcation parameter $\alpha$ will then vary in a neighborhood of the bifurcation point $\alpha = 0$. To cover all interesting cases, we shall consider three different functions $F$, namely

$$\text{I) } F(u) = 1 - u^2, \qquad \text{II) } F(u) = 1 - u, \qquad \text{III) } F(u) = 1 + u.$$

For all choices of $F$, system (1) possesses two spatially homogeneous equilibria $(u, v) = (\pm 1, 0)$ and a one-parameter family of front solutions

$$(u, v) = (\tanh((x - c_0 t - x_0)/2), 0), \qquad x_0 \in \mathbb{R}, \quad (2)$$

connecting these equilibria. For $\alpha < 0$, the equilibria and the family of front solutions are asymptotically stable with some exponential rate. When $\alpha$ crosses the origin from left to right, some of the equilibria become unstable, depending on the particular choice of $F$. In case I, the steady states ahead of and behind the front undergo a Turing bifurcation and spatially periodic equilibria are created. In case II, this happens only for the steady state $(u, v) = (1, 0)$ ahead of the front, and in case III only for the steady state $(u, v) = (-1, 0)$ behind the front. In this respect, case I is close to the case of a pulse.

At the bifurcation point $\alpha = 0$, the front solutions (2) become essentially unstable and, in cases I and II, a family of modulated fronts is created. These solutions are time-periodic in a moving frame with speed $c \approx c_0$, and they connect a spatially periodic Turing pattern at $x = +\infty$ to another Turing pattern (case I) or to the uniform steady state $(-1, 0)$ (case II) at $x = -\infty$, see Fig. 1. We shall not consider case III any longer, since the analysis in [SS01a] shows that at least generically no modulated fronts exist in that case; typically the pattern is outrun by the front, see fig.6 on page 44 for an illustration. As for the stability, it turns out that in case II the family of modulated fronts is asymptotically stable with exponential rate. This can be proved rather easily using weighted spaces, see section 6. Thus the challenging case in proving stability is case I. In this situation, the linearization around the modulating fronts has continuous spectrum up to the imaginary axis. It is the purpose of this paper to explain how nevertheless the nonlinear stability of these solutions can be shown.

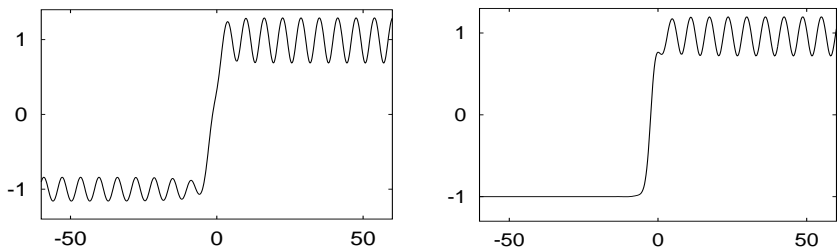

Figure 1: Modulated fronts for (1) in cases I (left) and II (right). The snapshots show the $u$-component obtained from generic initial data at some large time; see also section 7.

**Remark 1.1** *In case I, the model problem (1) seems non-generic, since both homogeneous equilibria $(\pm 1, 0)$ undergo a Turing bifurcation at the same value of the parameter $\alpha$. In fact, we implicitly restrict our analysis to systems for which the destabilization of both equilibria has the same origin (in our example, this is the coupling of the bistable $u$-equation to a single Swift-Hohenberg equation). This also explains why the wavelengths of the bifurcating patterns ahead of and behind the front coincide. As was observed by one of the referees of this paper, it would then be more natural to consider the case of a pulse instead of a front. But then we would have to replace the scalar $u$-equation in (1) by a $2 \times 2$ system, which makes the analysis even more intricate. Also, we found it interesting to encompass all possible cases (I, II, and III) in a single, relatively simple model.*

**Acknowledgments:** This work was supported by the French-German cooperation project PROCOPE 00307TK entitled "Attractors for extended systems". The authors also thank B. Sandstede and A. Scheel for stimulating discussions, and both referees for useful comments and suggestions.

## 2 Main Results

In this section, we give our results in the most interesting case, i.e. when $F(u)=1-u^2$ in (1). To simplify the notation, we rewrite (1) in the form

$$\partial_t U = L(\partial_x)U + N(U) , \qquad (3)$$

where $U = (u, v)$ and

$$L(\partial_x)\begin{pmatrix} u \\ v \end{pmatrix} = \begin{pmatrix} \partial_x^2 u \\ -(1+\partial_x^2)^2 v \end{pmatrix} , \quad N\begin{pmatrix} u \\ v \end{pmatrix} = \begin{pmatrix} \frac{1}{2}(u-c_0)(1-u^2) + v \\ \alpha v - v^3 - \gamma v(1-u^2) \end{pmatrix} .$$

In the invariant subspace $\{(u, v) \mid v = 0\}$, system (1) has exactly three homogeneous equilibria, namely $U_0 = (c_0, 0)$ and $U_\pm = (\pm 1, 0)$. In addition, there exists a family of traveling waves $U_h(x, t) = (h(x - c_0 t, 0))$ connecting $U_-$ to $U_+$. The profile $h$ satisfies the ordinary differential equation

$$h'' + c_0 h' + (h - c_0)(1 - h^2)/2 = 0 , \qquad (4)$$

together with the boundary conditions $h(\pm\infty) = \pm 1$. Up to translations, the unique solution is $h(y) = \tanh(y/2)$. The function $U_h(x, t)$ will be referred to as the "original front solution", as opposed to the "modulated front" which will be considered below.

To study the stability of the front solution $U_h$, it is advantageous to go to a comoving frame. The new space variable will be denoted by $y$, i.e. $y = x - c_0 t$. System (1) then reads

$$\partial_t U = L(\partial_y) U + c_0 \partial_y U + N(U) . \qquad (5)$$

We first investigate the stability of the homogeneous equilibria $U_\pm = (\pm 1, 0)$. Linearizing (5) at $U = U_\pm$, we obtain $\partial_t U = L(\partial_y)U + c_0 \partial_y U + DN(U_\pm)U$, or explicitly

$$\begin{aligned}
\partial_t u &= \partial_y^2 u + c_0 \partial_y u - (1 \mp c_0) u + v , \\
\partial_t v &= -(1 + \partial_y^2)^2 v + c_0 \partial_y v + \alpha v .
\end{aligned} \qquad (6)$$

The spectrum of the linear operator in the right-hand side is easily computed if we decompose $u, v$ in Fourier modes $e^{iky}$. We find $\Sigma^\pm = \{\hat\lambda_1^\pm(k) \mid k \in \mathbb{R}\} \cup \{\hat\lambda_2(k) \mid k \in \mathbb{R}\}$, where

$$\hat\lambda_1^\pm(k) = -k^2 - (1 \mp c_0) + c_0 \mathrm{i} k , \quad \hat\lambda_2(k) = -(1 - k^2)^2 + \alpha + c_0 \mathrm{i} k . \qquad (7)$$

Since $0 < c_0 < 1$ we immediately conclude that the trivial equilibria $U_\pm = (\pm 1, 0)$ are stable for $\alpha < 0$ and unstable for $\alpha > 0$.

We next consider the stability of the original front $U_h(y) = (h(y), 0)$, which is a steady state of (5).

**Notation:** For $n \in \mathbb{N}_0 = \mathbb{N} \cup \{0\}$ we denote by $C_b^n(\mathbb{R})$ the space of all functions $u : \mathbb{R} \to \mathbb{R}$ which are bounded and uniformly continuous together with their first $n$ derivatives. We equip $C_b^n(\mathbb{R})$ with the norm $\|u\|_{C_b^n} = \sum_{j=0}^n \sup_{x \in \mathbb{R}} |\partial_x^j u(x)|$.

**Theorem 2.1** *For $\alpha < 0$ the family $\{U_h(\cdot - y_0) \mid y_0 \in \mathbb{R}\}$ of front solutions is asymptotically stable with exponential rate $\mu > 0$ (depending on $\alpha$). More precisely, given any $C > 0$, there exists $\delta > 0$ such that, for all $U_0 \in [C_b^0(\mathbb{R})]^2$ with*

$$\inf_{y_0 \in \mathbb{R}} \sup_{y \in \mathbb{R}} \|U_0(y) - U_h(y - y_0)\|_{\mathbb{R}^2} \le \delta ,$$

*system (5) has a unique global solution $U \in C^0([0, +\infty), [C_b^0(\mathbb{R})]^2)$ with initial data $U_0$, and there exists $y_1 \in \mathbb{R}$ such that*

$$\sup_{y \in \mathbb{R}} \|U(y, t) - U_h(y - y_1)\|_{\mathbb{R}^2} \le C e^{-\mu t} \quad \text{for all } t \ge 0 .$$

**Proof.** The strategy is standard, see [Sat77] or section 5.4 in [Hen81]. Linearizing (5) at $U_h$, we obtain $\partial_t U = \Lambda U$, where $\Lambda = L(\partial_y) + c_0 \partial_y + DN(U_h)$. Explicitly,

$$\begin{aligned}
\partial_t u &= \partial_y^2 u + c_0 \partial_y u + \tfrac{1}{2}(1 + 2 c_0 h - 3 h^2) u + v , \\
\partial_t v &= -(1 + \partial_y^2)^2 v + c_0 \partial_y v + \alpha v - \gamma(1 - h^2) v .
\end{aligned}$$

We have to study the spectrum of the linear operator $\Lambda$ in $C_b^0(\mathbb{R})$, or equivalently in $L^2(\mathbb{R})$. Due to translation invariance of the original system, the operator $\Lambda$ has a zero eigenvalue with eigenfunction $U_h' = (h', 0)$. If $\alpha < 0$, we claim that the rest of the spectrum is strictly contained in the left half-plane of the complex plane. Indeed, the essential spectrum is determined by the linearization around the steady states $U_\pm$, hence it follows from (7) that $\Sigma_{\text{ess}}(\Lambda) \subset \{z \in \mathbb{C} \,|\, \Re(z) \le \lambda_0\}$, where $\lambda_0 = \max(\alpha, -1 + c_0) < 0$. Assume now that $\lambda$ is an isolated eigenvalue of $\Lambda$ with $\Re(\lambda) > \lambda_0$, and let $\hat{U} = (\hat{u}, \hat{v})$ be a (nonzero) eigenfunction. Then $\hat{U}(y)$ decays exponentially as $|y| \to \infty$, and $\hat{v}$ satisfies the decoupled equation

$$-(1 + \partial_y^2)^2 \hat{v} + c_0 \partial_y \hat{v} + \alpha \hat{v} - \gamma(1 - h^2)\hat{v} = \lambda \hat{v} \,.$$

Taking the scalar product of both sides with $\hat{v}$ and using the fact that $\gamma(1 - h^2) > 0$, we obtain the inequality

$$\Re(\lambda) \|\hat{v}\|_{L^2}^2 \le -\|(1 + \partial_y^2)\hat{v}\|_{L^2}^2 + \alpha\|\hat{v}\|_{L^2}^2 \le \alpha\|\hat{v}\|_{L^2}^2 \,,$$

which implies that $\hat{v} = 0$ due to $\alpha < 0$. It follows that $A\hat{u} = \lambda\hat{u}$, where $A$ is the second order differential operator

$$A = \partial_y^2 + c_0 \partial_y + \frac{1}{2}(1 + 2c_0 h - 3h^2) \,.$$

We know that 0 is a simple eigenvalue of $A$, and that the corresponding eigenfunction $h'$ is positive. By Sturm-Liouville theory, the other isolated eigenvalues of $A$ are all strictly negative. We conclude that either $\lambda = 0$ (in which case $\hat{u} = Ch'$ for some $C > 0$) or $\lambda < 0$. Thus, there exists $\mu > 0$ such that $\Sigma(\Lambda) \subset \{0\} \cup \{z \in \mathbb{C} \,|\, \Re(z) \le -\mu\}$. Now, applying for instance the center manifold theorem [Hen81], we obtain the desired result. $\square$

According to Theorem 2.1, for $\alpha < 0$ information can be transported in the system using the stable fronts $U_h$. We now consider the bifurcation that occurs when $\alpha$ crosses zero from left to right. In what follows, we set

$$\alpha = \varepsilon^2 > 0 \,,$$

where $\varepsilon > 0$ is a small parameter. It is clear from (7) that the homogeneous steady states $U_\pm$ are now unstable, and so is the front solution $U_h$. Remark that, when $\alpha$ crosses zero, the *essential spectrum* of the linearized operator $\Lambda$ crosses the imaginary axis, so that classical bifurcation theory is not applicable.

For later use we remark that, when $\varepsilon > 0$ is not too big, the spectrum of $\Lambda$ can be "stabilized" if we introduce exponentially weighted spaces, see [Sat77]. Indeed, if we set $U(y, t) = e^{-\beta y}\tilde{U}(y, t)$ for some $\beta > 0$, the linear equation $\partial_t U = \Lambda U$ becomes $\partial_t \tilde{U} = \Lambda_\beta \tilde{U}$, where

$$\Lambda_\beta = L(\partial_y - \beta) + c_0(\partial_y - \beta) + DN(U_h) \,.$$

**Proposition 2.2** *Fix $0 < c_1 < c_0$. There exists $\beta_0 > 0$ such that, if $0 < \beta < \beta_0$ and $\varepsilon^2 \le c_1 \beta$, there exists $\nu > 0$ (depending on $\beta$) such that the spectrum of $\Lambda_\beta$ satisfies*

$$\Sigma(\Lambda_\beta) \subset \{0\} \cup \{z \in \mathbb{C} \,|\, \Re(z) \le -2\nu\} \,.$$

**Proof.** The essential spectrum of $\Lambda_\beta$ is

$$\Sigma_{\text{ess}}(\Lambda_\beta) = \{\hat{\lambda}_1^+(k+\mathrm{i}\beta) \,|\, k \in \mathbb{R}\} \cup \{\hat{\lambda}_1^-(k+\mathrm{i}\beta) \,|\, k \in \mathbb{R}\} \cup \{\hat{\lambda}_2(k+\mathrm{i}\beta) \,|\, k \in \mathbb{R}\}\,,$$

see (7). Clearly, $\Re(\hat{\lambda}_1^\pm(k+\mathrm{i}\beta)) \leq -1+c_0(1-\beta)+\beta^2 < 0$ if $\beta$ is sufficiently small. Similarly,

$$\Re(\hat{\lambda}_2(k+\mathrm{i}\beta)) = \varepsilon^2 - c_0\beta + 4k^2\beta^2 - (1-k^2+\beta^2)^2 \leq \varepsilon^2 - c_0\beta + 4\beta^2(1+2\beta^2)\,.$$

Thus, if $\varepsilon^2 \leq c_1\beta$ and $\beta > 0$ is sufficiently small, then $\Sigma_{\text{ess}}(\Lambda_\beta) \subset \{z \in \mathbb{C} \,|\, \Re(z) \leq \lambda_0\}$ for some $\lambda_0 < 0$ depending on $\beta$.

Assume now that $\lambda$ is an isolated eigenvalue of $\Lambda_\beta$ with $\Re(\lambda) > \lambda_0$, and let $\tilde{U} = (\tilde{u}, \tilde{v})$ be a nonzero eigenfunction. Proceeding exactly as in Theorem 2.1, we show that $\tilde{v} = 0$ and that $u(y) = e^{-\beta y}\tilde{u}(y)$ is an eigenfunction of the Sturm-Liouville operator $A$, so that either $\lambda = 0$ (in which case $u = Ch'$ for some $C > 0$) or $\lambda < 0$. This concludes the proof. □

Of course, Proposition 2.2 does not imply stability of the front $U_h$ when $\alpha > 0$, because the nonlinear terms cannot be controlled in the weighted space. Nevertheless, the spectral stabilization property will be one of the key ingredients in the stability proof of the modulated fronts, see section 5.1 and section 6. The proof of Proposition 2.2 is also the only place where a particular structure of our model system is really used, see Remark 2.10.

In order to find new stable structures near the original front $U_h$, we first consider the bifurcation scenario for the homogeneous steady states $U_\pm$. If we restrict ourselves to the space of periodic functions with period $2\pi$, we can apply classical bifurcation theory. Indeed, the spectrum of the linearization $L(\partial_x) + DN(U_\pm)$ consists of the eigenvalues

$$\lambda_1^\pm(k) = -k^2 - (1\mp c_0)\,, \quad \lambda_2(k) = -(1-k^2)^2 + \alpha\,, \tag{8}$$

where $k \in \mathbb{Z}$ (or $k \in \mathbb{N}_0$ if we further restrict the space to *even* functions). In the latter case, as $\alpha$ crosses zero, a single eigenvalue $\lambda_2(1) = \alpha$ crosses the imaginary axis, while all the other ones stay negative and bounded away from the origin. As is easy to verify, this is a supercritical pitchfork bifurcation. Thus, for $\alpha = \varepsilon^2 > 0$ small enough, there exist stable periodic equilibria $U_{\text{per}}^\pm(x)$ satisfying $|U_{\text{per}}^\pm - U_\pm| = \mathcal{O}(\varepsilon)$. This pattern forming bifurcation is often referred to as a "Turing bifurcation", see [Tur52].

**Theorem 2.3** *Fix $c_0 \in (0,1)$. There exist $\varepsilon_0 > 0$ and $\gamma_0 > 0$ such that, for all $\varepsilon \in (0, \varepsilon_0)$ and all $\gamma \in (0, \gamma_0)$, there exist two families $\{U_{\text{per}}^\pm(x-x_0) \,|\, x_0 \in \mathbb{R}\}$ of smooth periodic equilibria of (3), satisfying $U_{\text{per}}^\pm(x) = U_{\text{per}}^\pm(x+2\pi)$ for all $x \in \mathbb{R}$, and*

$$U_{\text{per}}^\pm(x) = \left(\pm 1 + \frac{\varepsilon a^\pm}{2\mp c_0}\cos x, \varepsilon a^\pm \cos x\right) + \mathcal{O}(\varepsilon^2)\,,$$

*with $a^\pm(c_0, \gamma) = \frac{2}{\sqrt{3}} + \mathcal{O}(\gamma)$.*

**Proof.** See section 3.1. □

For simplicity, we restrict ourselves to the case of periodic equilibria with period $2\pi$, but the bifurcation argument above also applies to periodic functions with a nearby period, see [CE90] for a complete discussion in the case of the Swift-Hohenberg equation. As is well-known [Eck65], the bifurcating periodic equilibria are linearly stable if and only if their period is close enough to $2\pi$.

In what follows, we fix $c_0 \in (0,1)$, $\gamma \in (0, \gamma_0)$, and we always assume that $\varepsilon > 0$ is sufficiently small (in particular, $0 < \varepsilon < \varepsilon_0$). Although the linearization around the bifurcating equilibria $U_{\text{per}}^\pm$ has continuous spectrum all the way to the imaginary axis, the nonlinear stability of these solutions with respect to spatially localized perturbations can be shown using the techniques developed in [Sch96, Sch98a, Sch98b].

**Notation.** For $n \in \mathbb{N}$, let $H^n(\mathbb{R})$ be the (Sobolev) space of all functions $u \in L^2(\mathbb{R})$ whose first $n$ derivatives are also in $L^2(\mathbb{R})$, equipped with the norm $\|u\|_{H^n} = (\sum_{j=0}^{n} \|\partial_x^j u\|_{L^2}^2)^{1/2}$. For $s \geq 0$, we set $\mathrm{H}_s^n = \{u \in H^n(\mathbb{R}) \mid \rho^s u \in H^n(\mathbb{R})\}$ where $\rho(x) = (1+x^2)^{1/2}$. The space $\mathrm{H}_s^n$ is equipped with the norm $\|u\|_{\mathrm{H}_s^n} = \|\rho^s u\|_{H^n}$.

**Theorem 2.4** *Let $\varepsilon > 0$ be sufficiently small, and let $U_{\text{per}} = U_{\text{per}}^+$ or $U_{\text{per}} = U_{\text{per}}^-$, where $U_{\text{per}}^\pm$ are the periodic equilibria constructed in Theorem 2.3. There exist $C, \delta > 0$ such that, for all $V_0 \in (\mathrm{H}_2^2)^2$ with $\|V_0\|_{(\mathrm{H}_2^2)^2} \leq \delta$, equation (3) has a unique global solution $U(x,t) = U_{\text{per}}(x) + V(x,t)$ with initial data $U_{\text{per}} + V_0$. Moreover, $\|V(t)\|_{(L^\infty(\mathbb{R}))^2} \leq C(1+t)^{-1/2}$ for all $t \geq 0$.*

**Proof.** See section 3.4. □

**Remark 2.5** *Much more is known about the asymptotic behavior of the perturbation $V(x,t)$ as $t \to \infty$. Under the assumptions of Theorem 2.4 there exists $V_* \in \mathbb{R}$ and $d > 0$ such that*

$$\sup_{x \in \mathbb{R}} \left| V(x,t) - \frac{1}{\sqrt{t}} V_* \exp\left(-\frac{x^2}{4dt}\right) \partial_x U_{\text{per}}(x) \right| = \mathcal{O}(t^{-1+\eta}) \quad \text{as } t \to \infty,$$

*for some arbitray but fixed $\eta > 0$, see Theorem 2.8 below. Thus, spatially localized perturbations vanish asymptotically as a solution of a linear diffusion equation.*

Finally, we study the bifurcation that the front $U_h$ undergoes when $\alpha$ crosses zero. For $\alpha = \varepsilon^2$ sufficiently small, in addition to the (unstable) original front $U_h$, equation (3) has a family of modulated fronts connecting the stable equilibria $U_{\text{per}}^-$ and $U_{\text{per}}^+$. These bifurcating solutions are time-periodic in a frame moving with speed $c = c_0 + \mathcal{O}(\varepsilon^2)$, and their profile is $\mathcal{O}(\varepsilon)$-close to the original front $U_h$. This bifurcation scenario has been thoroughly studied by B. Sandstede and A. Scheel for general reaction-diffusion systems in [SS99, SS01a]. Unfortunately, our model problem (1) does not exactly fit into this abstract framework, because the fourth order Swift-Hohenberg equation is not a reaction-diffusion system. For this reason, the proof of the following result will be outlined in section 4.

**Theorem 2.6** *For $\varepsilon > 0$ sufficiently small there exists a modulated front solution of (3) of the form*

$$U(x,t) = U_{\text{mf}}(x-ct, x), \quad x \in \mathbb{R}, \quad t \in \mathbb{R},$$

where $U_{\mathrm{mf}}(\xi, x)$ is $2\pi$-periodic in its second argument and $c = c_0 + \mathcal{O}(\varepsilon^2)$. Moreover, there exist positive constants $C, \beta_1, \beta_2$ (independent of $\varepsilon$) such that

$$\sup_{\xi, x \in \mathbb{R}} |U_{\mathrm{mf}}(\xi, x) - U_h(\xi)| \leq C\varepsilon \,,$$

and

$$\|U_{\mathrm{mf}}(\xi, \cdot) - U_{\mathrm{per}}^+(\cdot + x_+)\|_{(H^2(0,2\pi))^2} \leq Ce^{-\beta_1 \xi} \,, \quad \xi \geq 0 \,,$$
$$\|U_{\mathrm{mf}}(\xi, \cdot) - U_{\mathrm{per}}^-(\cdot + x_-)\|_{(H^2(0,2\pi))^2} \leq Ce^{\varepsilon \beta_2 \xi} \,, \quad \xi \leq 0 \,,$$

for some $x_\pm \in [0, 2\pi)$.

**Proof.** See section 4. □

**Remark 2.7** *Due to translation invariance of the original problem, $U_{\mathrm{mf}}(x - ct - x_0, x - x_1)$ is also a solution of (3) for all $x_0, x_1 \in \mathbb{R}$. Thus, without loss of generality, we can assume that $x_- = 0$ in Theorem 2.6. We may also break the translation invariance in the variable $\xi$ by imposing $U_{\mathrm{mf}}(0, 0) = 0$.*

We are now able to state our main result, which shows that the family of modulated fronts is asymptotically stable with respect to small, localized perturbations. We recall that $\beta_0$ is the positive constant defined in Proposition 2.2.

**Theorem 2.8** *For $\beta \in (0, \beta_0)$ and $\varepsilon > 0$ sufficiently small, there exist positive constants $C$, $\nu$, $\delta$, $d$ such that the following holds. For all $V_0 : \mathbb{R} \to \mathbb{R}^2$ with $\|V_0(x)(x^2 + e^{\beta x})\|_{(H^2)^2} \leq \delta$ there exists a unique global solution $U(x, t)$ of (3) with initial data $U(x, 0) = U_{\mathrm{mf}}(x, x) + V_0(x)$. Moreover, there exists a shift function $q : \mathbb{R}_+ \to \mathbb{R}$ and two real constants $q_*, V_*$ such that $U(x, t)$ can be represented as*

$$U(x, t) = U_{\mathrm{mf}}(x - ct - q(t), x) + V(x, t) \,, \quad x \in \mathbb{R}, \, t \geq 0 \,,$$

*where*

$$\sup_{x \in \mathbb{R}} \left| V(x, t) - \frac{1}{\sqrt{t}} V_* \exp\left(-\frac{x^2}{4dt}\right) \partial_x U_{\mathrm{per}}^-(x) \right| \leq \frac{C}{t^{3/4}} \,, \quad t \geq 1 \,, \quad (9)$$

*and*

$$|q(t) - q_*| + \sup_{\xi \in \mathbb{R}} |V(\xi + ct, t)e^{\beta \xi}| \leq Ce^{-\nu t} \,, \quad t \geq 0 \,. \quad (10)$$

**Proof.** See section 5. □

**Remark 2.9** *From the proof it will be clear that the decay in (9) can be improved to $t^{-1+\eta}$ with arbitrary small $\eta > 0$. For simplicity we stick to $t^{-3/4}$.*

**Remark 2.10** *As we explain in section 5.1, the essential properties of system (3) that we use are the stability of the Turing pattern $U_{\mathrm{per}}^-$ behind the front, and the fact that the spectrum of the linearized operator can be stabilized using appropriate weighted spaces. Thus, Theorem 2.8 will hold for any of the reaction-diffusion systems considered in [SS01a] provided one can prove the analogue of Proposition 2.2 and Theorem 2.4.*

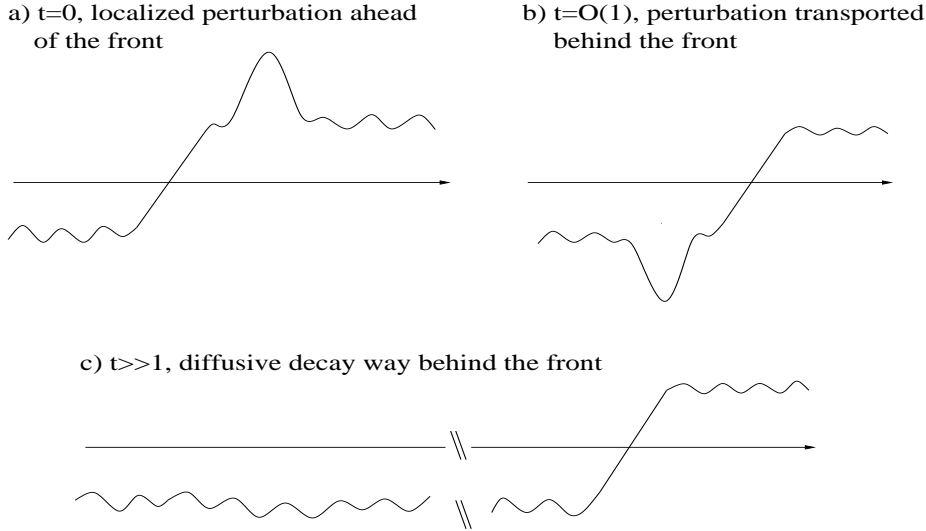

Figure 2: Spatially localized perturbations (a) are transported behind the modulated front (b), where they are damped diffusively (c).

Theorem 2.8 shows that small, spatially localized perturbations of the modulated front do not destroy the form of $U_{\mathrm{mf}}$, but lead to a finite shift $q_*$. In fact, the perturbation is transported behind the front, where it vanishes diffusively, see fig.2. Thus, information can be transported in a stable manner even in essentially unstable media exhibiting a Turing bifurcation.

Exponential weights have been widely used to prove the stability of fronts propagating into *unstable* states, see [Sat77, CE90, Gal94, ES02]. In our case, the invaded steady state $U_{\mathrm{per}}^+$ is in fact stable, but this fact is *not* used in the proof of Theorem 2.8. If we use in addition the stability of the equilibrium $U_{\mathrm{per}}^+$, it should be possible to replace the exponential weight $e^{\beta x}$ with a polynomial one, in which case the convergence of $q(t)$ and the decay of $V(x+ct)$ in (10) will be algebraic (i.e. like some inverse power of $t$). In any case, the decay of $V(x,t)$ in the laboratory frame will always be algebraic, because this is what we have for localized perturbations of the periodic steady state $U_{\mathrm{per}}^-(x)$, see Theorem 2.4.

This paper stands in line with [BK92, ES02, ES00], where the diffusive stability of a ground state, here the spatially periodic equilibria $U_{\mathrm{per}}^\pm$, has been used to prove diffusive stability of more complicated structures. In contrast to these papers, we have to deal here with an additional zero eigenvalue which leads to the shift $q(t)$ in Theorem 2.8.

We proceed as follows. In section 3 we prove the existence and the stability of the Turing patterns $U_{\mathrm{per}}^\pm(x)$. Section 4 contains a short proof of Theorem 2.6, i.e. the construction of the modulated front solutions. In section 5 we prove Theorem 2.8. The proof is based on a strategy similar to the one used in [ES02]. The method is improved in the sense that in contrast to [ES02], where $\beta = \mathcal{O}(\varepsilon)$, here we can choose $\beta = \mathcal{O}(1)$ which simplifies the proof. Section 6 is devoted to the stability of modulated fronts in case II. As already said, these solutions are asymptotically stable with an exponential rate. In section 7 we illustrate our analysis with numerical simulations, also showing computer experiments for a model with modulated pulses. These results indicate that the assertions of Theorem 2.8 remain true for relatively large values of $\varepsilon$ and $\delta$.

# 3 Existence and stability of Turing patterns

The aim of this section is to give a short proof of Theorems 2.3 and 2.4. This will be done by adapting to system (1) some known results about existence and stability of periodic solutions for the Swift-Hohenberg equation. The existence part is a rather standard bifurcation argument, see e.g. [CE90] for more details. The stability proof is based on previous results by the second author [Sch96, Sch98a, Sch98b].

## 3.1 Existence

We are interested in $2\pi$-periodic stationary solutions of (3) that bifurcate from the homogeneous steady states $U_\pm = (\pm 1, 0)$ as the parameter $\alpha$ becomes positive. Since equation (3) remains invariant if $(U, c_0)$ is replaced with $(-U, -c_0)$, it is sufficient to investigate the solutions that bifurcate from $U_- = (-1, 0)$. Setting $U = U_- + W$, we obtain the system

$$\partial_t W = \mathcal{L}W + \mathcal{N}(W), \tag{11}$$

where

$$\mathcal{L}\begin{pmatrix} w_1 \\ w_2 \end{pmatrix} = \begin{pmatrix} \partial_x^2 w_1 - (1+c_0)w_1 + w_2 \\ -(1+\partial_x^2)^2 w_2 + \alpha w_2 \end{pmatrix}, \quad \mathcal{N}\begin{pmatrix} w_1 \\ w_2 \end{pmatrix} = \begin{pmatrix} \tfrac{1}{2}(3+c_0)w_1^2 - \tfrac{1}{2}w_1^3 \\ -w_2^3 - \gamma w_2(2w_1 - w_1^2) \end{pmatrix}. \tag{12}$$

We study system (11) in the Hilbert space $X = H^1_{\text{per},+}(\mathbb{R})^2$, where

$$H^1_{\text{per},+}(\mathbb{R}) = \{u \in H^1_{\text{loc}}(\mathbb{R}) \mid u(x) = u(x+2\pi),\ u(x) = u(-x),\ \forall x \in \mathbb{R}\}.$$

In other words, in the space of $2\pi$-periodic functions, we freeze the translation invariance by assuming the function $W$ to be even, a symmetry that is preserved under evolution. In the space $X$, the linear operator $\mathcal{L}$ in (11) has compact resolvent, hence purely discrete spectrum. Its eigenvalues are $\{\lambda_1^-(k)\}_{k\in\mathbb{N}_0}$ and $\{\lambda_2(k)\}_{k\in\mathbb{N}_0}$, where

$$\lambda_1^-(k) = -k^2 - (1+c_0), \quad \lambda_2(k) = -(1-k^2)^2 + \alpha,$$

see (8). As we shall show, when the largest eigenvalue $\lambda_2(1) = \alpha$ crosses the origin (from left to right), a supercritical pitchfork bifurcation occurs: the origin $W = (0,0)$ looses its stability, and a pair of stable equilibria is created at a distance $\mathcal{O}(\sqrt{\alpha})$ of the origin.
If $\alpha > -(1+c_0)$, the largest eigenvalue $\alpha$ of the operator $\mathcal{L}$ in $L^2_{\text{per},+}$ is simple, with eigenvector

$$\Phi(x) = (A\cos(x), \cos(x)), \quad \text{where } A = \frac{1}{2+c_0+\alpha}.$$

Let $E_c = \{r\Phi \mid r \in \mathbb{R}\}$ and $E_s = (1-P)X$, where $P : X \to X$ is the spectral projection onto the one-dimensional eigenspace $E_c$ of $\mathcal{L}$. By construction, $PW = (\Pi W)\Phi$ for all $W \in X$, where $\Pi : X \to \mathbb{R}$ is the bounded linear form defined by

$$\Pi\begin{pmatrix} w_1 \\ w_2 \end{pmatrix} = \frac{1}{\pi}\int_0^{2\pi} w_2(x)\cos(x)\,\mathrm{d}x.$$

From now on, we assume that $|\gamma| \le \gamma_0$ for some $\gamma_0 > 0$ which will be fixed later. By the center manifold theorem, if $|\alpha|$ is sufficiently small, system (11) has a one-dimensional local center manifold of the form

$$\mathcal{V}_c = \{r\Phi + f(r) \,|\, |r| < r_0\},$$

where $r_0 > 0$ and $f : (-r_0, r_0) \to E_s$ is a $C^3$ function satisfying $f(0) = 0$, $f'(0) = 0$. In addition, $f$ maps $(-r_0, r_0)$ into the domain of $\mathcal{L}$, and the following identity holds:

$$f'(r)\Big(\alpha r + \Pi \mathcal{N}(r\Phi + f(r))\Big) = \mathcal{L}f(r) + (1 - P)\mathcal{N}(r\Phi + f(r)), \quad |r| < r_0. \tag{13}$$

The evolution defined by (11) on $\mathcal{V}_c$ is given by the reduced system

$$\dot r = \alpha r + \Pi \mathcal{N}(r\Phi + f(r)). \tag{14}$$

Since $f(r) = \mathcal{O}(r^2)$ as $r \to 0$, it follows from (12) that $\mathcal{N}(r\Phi + f(r)) = r^2 \Psi + \mathcal{O}(r^3)$, where

$$\Psi(x) = \Big(\frac{3 + c_0}{2} A^2 \cos^2(x), -2\gamma A \cos^2(x)\Big).$$

Remark that $\Pi\Psi = 0$, so that $\Pi \mathcal{N}(r\Phi + f(r)) = \mathcal{O}(r^3)$. On the other hand, since $f \in C^3$, there exists $\Xi \in E_s$ such that $f(r) = \Xi r^2 + \mathcal{O}(r^3)$. Inserting this expansion into (13) and keeping only the lowest order terms in $r$, we obtain the relation $(\mathcal{L} - 2\alpha)\Xi + \Psi = 0$. It follows that

$$\Xi(x) = (B\cos(2x) + D, b\cos(2x) + d),$$

where

$$b = \frac{-\gamma A}{9 + \alpha}, \quad d = \frac{-\gamma A}{1 + \alpha}, \quad B = \frac{b + \frac{3+c_0}{4}A^2}{5 + c_0 + 2\alpha}, \quad D = \frac{d + \frac{3+c_0}{4}A^2}{1 + c_0 + 2\alpha}.$$

Using this information, we conclude that $\Pi \mathcal{N}(r\Phi + f(r)) = -ar^3 + \mathcal{O}(r^4)$, where

$$a = \frac{3}{4} + \gamma\Big(B + 2D + A(b + 2d) - \frac{3}{4}A^2\Big).$$

It is clear that $a > 0$ if $\gamma \in (0, \gamma_0)$ for $\gamma_0 > 0$ sufficiently small. In this case, if $\alpha = \varepsilon^2$ is small enough, equation (14) has exactly three equilibria in a neighborhood of size $r_0$ of the origin: $r = 0$ and $r = r_\pm$, where

$$r_\pm = \pm \frac{\varepsilon}{\sqrt{a}} + \mathcal{O}(\varepsilon^2), \quad \text{as } \varepsilon \to 0.$$

By the center manifold theorem, equation (11) has also three equilibria in a neighborhood of zero in $X$, namely $W = 0$ and $W = W_\pm = r_\pm \Phi + f(r_\pm)$. Since (11) is translation invariant in the space variable $x$, it is straightforward to verify that $W_-(x) = W_+(x + \pi)$ for all $x \in \mathbb{R}$. Thus, equation (11) has in fact a unique family $\{W_{\text{per}}(x - x_0) \,|\, x_0 \in \mathbb{R}\}$ of non-constant $2\pi$-periodic solutions in a neighborhood of the origin, where

$$W_{\text{per}} = \frac{\varepsilon}{\sqrt{a}}\Phi + \mathcal{O}(\varepsilon^2).$$

Setting $U_{\text{per}}^- = (-1, 0) + W_{\text{per}}$, we obtain spatially periodic equilibria of (3) with the desired properties. As was already mentioned, $U_{\text{per}}^+$ is obtained by replacing $c_0$ with $-c_0$ in the expression of $-U_{\text{per}}^-$. This concludes the proof of Theorem 2.3.

## 3.2 Bloch waves

The first step in the stability analysis of the Turing patterns is studying the linearization of (3) around $U^{\pm}_{\mathrm{per}}$. Since $U^{\pm}_{\mathrm{per}}(x)$ are spatially periodic functions the associated eigenvalue problem is naturally formulated in terms of Bloch waves. In this section we briefly recall the definition of the Bloch wave transform and list a few identities that will be useful in the sequel. For a rigorous introduction to Bloch waves techniques, we refer to [RS72].

The starting point of Bloch wave analysis in the case of a $2\pi$-periodic underlying pattern is the following (formal) relation

$$\begin{aligned} u(x) &= \int_{-\infty}^{+\infty} e^{ikx} \tilde{u}(k) \, \mathrm{d}k = \sum_{n \in \mathbb{Z}} \int_{-1/2}^{1/2} e^{i(n+\ell)x} \tilde{u}(n+\ell) \, \mathrm{d}\ell \\ &= \int_{-1/2}^{1/2} \sum_{n \in \mathbb{Z}} e^{i(n+\ell)x} \tilde{u}(n+\ell) \, \mathrm{d}\ell = \int_{-1/2}^{1/2} e^{i\ell x} \hat{u}(\ell, x) \, \mathrm{d}\ell \,, \end{aligned} \qquad (15)$$

where $\tilde{u}$ is the Fourier transform of $u$ defined by

$$(\mathcal{F}u)(k) \equiv \tilde{u}(k) = \frac{1}{2\pi} \int_{-\infty}^{+\infty} u(x) e^{-ikx} \, \mathrm{d}x \,,$$

and $\hat{u}$ is the Bloch wave transform of $u$ defined by

$$(\mathcal{T}u)(\ell, x) \equiv \hat{u}(\ell, x) = \sum_{n \in \mathbb{Z}} e^{inx} \tilde{u}(n+\ell) \,. \qquad (16)$$

From Plancherel's theorem and Parseval's identity we easily deduce the relation

$$\int_{\mathbb{R}} |u(x)|^2 \, \mathrm{d}x = \int_{-1/2}^{1/2} \left( \int_0^{2\pi} |\hat{u}(\ell, x)|^2 \, \mathrm{d}x \right) \mathrm{d}\ell \,,$$

which shows that the Bloch wave transform $\mathcal{T}$ defined in (16) is an isomorphism between $L^2(\mathbb{R})$ and $L^2((-1/2, 1/2) \times (0, 2\pi))$. The inverse transform is given by (15), namely

$$u(x) = (\mathcal{T}^{-1}\hat{u})(x) = \int_{-1/2}^{1/2} e^{i\ell x} \hat{u}(\ell, x) \, \mathrm{d}\ell \,.$$

We note the useful elementary properties

$$\begin{aligned} \hat{u}(\ell, x) &= e^{ix} \hat{u}(\ell+1, x) \,, \\ \hat{u}(\ell, x) &= \hat{u}(\ell, x+2\pi) \,, \\ \hat{u}(\ell, x) &= \overline{\hat{u}(-\ell, x)} \quad \text{for real-valued } u \,. \end{aligned}$$

The Bloch wave transform of the product $uv$ is a convolution

$$(\widehat{uv})(\ell, x) = \int_{-1/2}^{1/2} \hat{u}(\ell - \ell', x) \hat{v}(\ell', x) \, \mathrm{d}\ell' \equiv (\hat{u} \star \hat{v})(\ell, x) \,.$$

On the other hand, if $v$ is $2\pi$–periodic we simply have

$$(\widehat{uv})(\ell, x) = \hat{u}(\ell, x) v(x) \,.$$

Let $n, s \in \mathbb{N}$, and let $\mathrm{H}^n_s = \mathrm{H}^n_s(\mathbb{R})$ be the weighted Sobolev space defined by the norm

$$\|u\|_{\mathrm{H}^n_s} \;=\; \left( \sum_{m=0}^n \int_{-\infty}^\infty |\partial_x^m u(x)|^2 (1+x^2)^s \,\mathrm{d}x \right)^{1/2} .$$

The image of $\mathrm{H}^n_s$ under the Bloch wave transformation is the space $\hat{\mathrm{H}}^n_s$ defined by the norm

$$\|\hat{u}\|_{\hat{\mathrm{H}}^n_s} \;=\; \left( \sum_{j=0}^s \sum_{m=0}^n \int_{-1/2}^{1/2} \int_0^{2\pi} |\partial_\ell^j \partial_x^m \hat{u}(\ell, x)|^2 \,\mathrm{d}x \,\mathrm{d}\ell \right)^{1/2} .$$

It is not difficult to show that there exists $C \geq 1$ such that

$$C^{-1} \|u\|_{\mathrm{H}^n_s} \;\leq\; \|\hat{u}\|_{\hat{\mathrm{H}}^n_s} \;\leq\; C \|u\|_{\mathrm{H}^n_s} ,$$

for all $u \in \mathrm{H}^n_s$.

In the sequel, we mainly work with the spaces corresponding to $s = n = 2$ or $s = 0, n = 2$. To bound the nonlinear terms we need the following estimates. If $u \in \mathrm{H}^2_2$ and $v \in \mathrm{H}^2_0$, then $uv \in \mathrm{H}^2_2$ and

$$\|\widehat{uv}\|_{\hat{\mathrm{H}}^2_2} \;=\; \|\hat{u} \star \hat{v}\|_{\hat{\mathrm{H}}^2_2} \;\leq\; C \|\hat{u}\|_{\hat{\mathrm{H}}^2_2} \,\|\hat{v}\|_{\hat{\mathrm{H}}^2_0} .$$

If $u \in \mathrm{H}^2_0$ and $f \in \mathcal{C}^2_\mathrm{b}(\mathbb{R})$, then $fu \in \mathrm{H}^2_0$ and

$$\|\widehat{fu}\|_{\hat{\mathrm{H}}^2_0} \;\leq\; C \|f\|_{\mathcal{C}^2_\mathrm{b}} \,\|\hat{u}\|_{\hat{\mathrm{H}}^2_0} , \quad \text{where} \quad \|f\|_{\mathcal{C}^2_\mathrm{b}} = \sum_{j=0}^2 \sup_{x \in \mathbb{R}} |\partial_x^j f(x)| .$$

**Notation.** If $\mathcal{A}$ is a linear operator, we define its Bloch wave transform by $\hat{\mathcal{A}} = \mathcal{T}\mathcal{A}\mathcal{T}^{-1}$. For instance, if $T_\zeta$ is the translation operator defined by $(T_\zeta f)(x) = f(x - \zeta)$, then

$$(\widehat{T_\zeta}\hat{u})(\ell, x) \;=\; e^{-i\ell\zeta} \hat{u}(\ell, x - \zeta) .$$

## 3.3 Linear stability

In this section, we study the linearization of (3) around the spatially periodic equilibria $U^\pm_{\mathrm{per}}$. Since only the stability of $U^-_{\mathrm{per}}$ will be used in the proof of Theorem 2.8, we shall concentrate here on this case. Setting $U = U^-_{\mathrm{per}} + V$ in (3), we obtain for $V$ the linearized equation

$$\partial_t V \;=\; \mathcal{M} V , \qquad \text{or equivalently} \qquad \partial_t \hat{V} \;=\; \hat{\mathcal{M}} \hat{V} , \tag{17}$$

with $\mathcal{M} = L(\partial_x) + DN(U^-_{\mathrm{per}})$. Our goal is to localize the spectrum of the linear operator $\mathcal{M}$ in the space $L^2(\mathbb{R})$ (actually, the result would be the same in $\mathrm{H}^2_2$.) Since this question is well-documented in the literature, we just summarize here the results. On this occasion, we also introduce some notations which we will use in section 5.

As the linearized problem has periodic coefficients, the operator $\hat{\mathcal{M}} = \mathcal{T}\mathcal{M}\mathcal{T}^{-1}$ can be represented as direct integral $\int^\oplus \mathcal{M}(\ell) \,\mathrm{d}\ell$, where, for each $\ell \in [-1/2, 1/2]$, $\mathcal{M}(\ell)$ is the linear operator on $H^2_{\mathrm{per}}([0, 2\pi])$ defined by $\mathcal{M}(\ell) w = e^{-i\ell x} \mathcal{M}(e^{i\ell x} w)$. The spectrum of $\mathcal{M}(\ell)$ is a sequence of eigenvalues $\{\mu_n(\ell)\}_{n\in\mathbb{N}}$, where $\Re(\mu_n(\ell)) \to -\infty$ as $n \to \infty$. The corresponding eigenfunctions, which are $2\pi$-periodic, will be denoted by $w_{\ell,n}$. Then the spectrum of $\mathcal{M}$ is given by

$$\sigma(\mathcal{M}) = \left\{\mu_n(\ell) \,\big|\, \ell \in [-1/2, 1/2],\, n \in \mathbb{N}\right\}.$$

By construction, the Bloch waves $e^{i\ell x} w_{\ell,n}$ satisfy $\mathcal{M}(e^{i\ell x} w_{\ell,n}) = \mu_n(\ell) e^{i\ell x} w_{\ell,n}$. As in [Sch96], one has the following result:

**Lemma 3.1** *For $\varepsilon > 0$ sufficiently small, there exist $\ell_0 > 0$ and $\nu_0 > 0$ such that*
**a)** *If $|\ell| \geq \ell_0$, then $\Re(\mu_n(\ell)) \leq -\nu_0$ for all $n \in \mathbb{N}$.*
**b)** *If $|\ell| < \ell_0$, the principal eigenvalue $\mu_1(\ell)$ is isolated and has the expansion*

$$\mu_1(\ell) = -d\ell^2 + \mathcal{O}(\ell^4), \quad as\ \ell \to 0, \tag{18}$$

*where $d = 4 + \mathcal{O}(\varepsilon)$. The corresponding eigenfunction $\varphi(\ell) = w_{\ell,1}$ depends smoothly on $\ell$ and satisfies $\varphi(0) = c_N \partial_x U_{\mathrm{per}}^-$, where $c_N > 0$ is a normalization such that $\|\varphi(0)\|_{L^2(0,2\pi)} = 1$. Finally, $\Re(\mu_n(\ell)) \leq -\nu_0$ for all $n \geq 2$.*

**Proof.** For all $\ell \in [-1/2, 1/2]$, the linear operator $\mathcal{M}(\ell)$ is explicitly given by

$$\mathcal{M}(\ell) = \begin{pmatrix} (\partial_x + i\ell)^2 + f_1 & 1 \\ f_2 & -(1 + (\partial_x + i\ell)^2)^2 + f_3 \end{pmatrix},$$

where $U_{\mathrm{per}}^- = (u_{\mathrm{per}}, v_{\mathrm{per}})$ and

$$f_1 = \tfrac{1}{2}(1 + 2c_0 u_{\mathrm{per}} - 3u_{\mathrm{per}}^2) = -(1 + c_0) + \mathcal{O}(\varepsilon),$$
$$f_2 = 2\gamma u_{\mathrm{per}} v_{\mathrm{per}} = \mathcal{O}(\varepsilon), \quad f_3 = \varepsilon^2 - \gamma(1 - u_{\mathrm{per}}^2) - 3v_{\mathrm{per}}^2 = \mathcal{O}(\varepsilon).$$

Thus, $\mathcal{M}(\ell)$ is a small, bounded perturbation of the constant coefficients operator obtained by setting $\varepsilon = 0$, namely

$$\mathcal{M}^0(\ell) = \begin{pmatrix} (\partial_x + i\ell)^2 - (1 + c_0) & 1 \\ 0 & -(1 + (\partial_x + i\ell)^2)^2 \end{pmatrix}.$$

The eigenvalues of $\mathcal{M}^0(\ell)$ are given by

$$\lambda_1^-(k + \ell) = -(k + \ell)^2 - (1 + c_0), \quad \lambda_2(k + \ell) = -(1 - (k + \ell)^2)^2, \quad k \in \mathbb{Z}.$$

Observe that these eigenvalues are all bounded away from zero, except for $\lambda_2(\pm 1 + \ell)$ which touch the origin when $\ell = 0$. Therefore, if $0 < \ell_1 < 1/2$ and $\varepsilon > 0$ is sufficiently small, the following holds for the eigenvalues $\mu_n(\ell)$ of the perturbed operator $\mathcal{M}(\ell)$:
i) If $|\ell| \geq \ell_1$, then $\Re(\mu_n(\ell)) \leq -\ell_1^2$ for all $n \in \mathbb{N}$.
ii) If $|\ell| < \ell_1$, then $\Re(\mu_n(\ell)) \leq -1/2$ for all $n \geq 3$.
As is clear from ii), we denote by $\mu_1(\ell)$ and $\mu_2(\ell)$ the two eigenvalues of $\mathcal{M}(\ell)$ that bifurcate from $\lambda_2(\pm 1 + \ell)$ when $\varepsilon$ is nonzero.

On the other hand, when $\varepsilon > 0$ and $\ell = 0$, we know that $\mathcal{M}(0)$ has a zero eigenvalue with eigenfunction $\partial_x U_{\text{per}}^-$, so that $\mu_1(0) = 0$. Moreover, it follows from the analysis in section 3.1 that $\mu_2(0) = -2\varepsilon^2 + \mathcal{O}(\varepsilon^3)$. Indeed, by construction, $\mu_2(0)$ is the convergence rate in time of $2\pi$-periodic solutions of (3) towards the circle of equilibria $\{U_{\text{per}}^-(x-x_0) \,|\, x_0 \in [0, 2\pi]\}$. This rate can be computed on the one-dimensional center manifold $\mathcal{V}_c$. The motion on $\mathcal{V}_c$ is given by $\dot r = g(r) = \varepsilon^2 r - ar^3 + \mathcal{O}(r^4)$, and the steady state $U_{\text{per}}^-$ corresponds to $r = r_+ = \varepsilon/\sqrt{a} + \mathcal{O}(\varepsilon^2)$. Then $g(r_+) = 0$ and $\mu_2(0) = g'(r_+) = -2\varepsilon^2 + \mathcal{O}(\varepsilon^3)$. Thus, in contrast to $\mathcal{M}^0(\ell)$, the operator $\mathcal{M}(\ell)$ has a *simple* zero eigenvalue when $\ell = 0$.

Now, a straightforward expansion in the parameter $\ell$ shows that the first eigenvalue $\mu_1(\ell)$ satisfies (18), and that the corresponding eigenfunction $\varphi(\ell)$ is proportional to

$$\begin{pmatrix} \frac{1}{2+c_0} \\ 1 \end{pmatrix} \sin(x) + \mathrm{i}\ell \begin{pmatrix} \frac{2}{(2+c_0)^2} \\ 0 \end{pmatrix} \cos(x) + \mathcal{O}(\varepsilon + \ell^2) \,.$$

In particular, $\varphi(0) = c_N \partial_x U_{\text{per}}^-$. Finally, the second eigenvalue $\mu_2(\ell)$ has the expansion

$$\mu_2(\ell) = -2\varepsilon^2 - d\ell^2 + \mathcal{O}(\varepsilon^3 + \ell^4) \,, \tag{19}$$

where $d = 4 + \mathcal{O}(\varepsilon)$.

If $\ell_1 > 0$ and $\varepsilon > 0$ are small enough, it follows from (18), (19) that

iii) If $|\ell| \leq \ell_1$, then $\mu_1(\ell) \leq -2\ell^2$ and $\mu_2(\ell) \leq -\varepsilon^2 - 2\ell^2$.

iv) If $|\ell| \leq \varepsilon$, then $0 \geq \mu_1(\ell) > \mu_2(\ell) > \Re(\mu_n(\ell))$ for all $n \geq 3$.

Combining i)–iv), we see that Lemma 3.1 holds with, for instance, $\ell_0 = \varepsilon$ and $\nu_0 = \varepsilon^2$. $\square$

From now on, we fix $\varepsilon > 0$ small enough so that the conclusion of Lemma 3.1 holds. We define the central projections $\hat{P}_{\text{c}}(\ell) : \mathrm{H}_{\text{per}}^2 \to \mathrm{H}_{\text{per}}^2$ by

$$\hat{P}_{\text{c}}(\ell)f = \langle \psi(\ell), f \rangle \varphi(\ell) \,, \quad |\ell| < \ell_0 \,, \tag{20}$$

where $\langle \cdot, \cdot \rangle$ is the usual scalar product in $\mathrm{L}^2([0, 2\pi])^2$ and $\psi(\ell)$ is the solution of the adjoint eigenvalue problem $\mathcal{M}^*(\ell)\psi(\ell) = \mu_1(\ell)\psi(\ell)$ normalized so that $\langle \psi(\ell), \varphi(\ell) \rangle = 1$.

Since we work in $\hat{\mathrm{H}}_2^2$, we will also need a version of the projection that depends smoothly on the variable $\ell$. To do that, we fix once and for all a nonnegative smooth cut-off function $\chi$ with support in $[-\ell_0/2, \ell_0/2]$ which is equal to 1 on $[-\ell_0/4, \ell_0/4]$. Then we define the operators $\hat{E}_{\text{c}}, \hat{E}_{\text{s}} : \mathrm{H}_{\text{per}}^2 \to \mathrm{H}_{\text{per}}^2$ by

$$\hat{E}_{\text{c}}(\ell) = \chi(\ell)\hat{P}_{\text{c}}(\ell) \,, \quad \hat{E}_{\text{s}}(\ell) = \mathbf{1}(\ell) - \hat{E}_{\text{c}}(\ell) \,. \tag{21}$$

It will be useful to define auxiliary mode filters $\hat{E}_{\text{c}}^{\text{h}}$ and $\hat{E}_{\text{s}}^{\text{h}}$ by

$$\hat{E}_{\text{c}}^{\text{h}}(\ell) = \chi(\ell/2)\hat{P}_{\text{c}}(\ell) \,, \quad \hat{E}_{\text{s}}^{\text{h}}(\ell) = \mathbf{1}(\ell) - \chi(2\ell)\hat{P}_{\text{c}}(\ell) \,.$$

These definitions are made in such a way that $\hat{E}_{\text{c}}^{\text{h}} \hat{E}_{\text{c}} = \hat{E}_{\text{c}}$ and $\hat{E}_{\text{s}}^{\text{h}} \hat{E}_{\text{s}} = \hat{E}_{\text{s}}$. As a consequence of Lemma 3.1, there exists $C > 0$ and $\nu_1 > 0$ (depending on $\varepsilon$) such that, for all $\hat{u} \in \hat{\mathrm{H}}_2^2$,

$$\|e^{t\hat{\mathcal{M}}} \hat{E}_{\text{s}}^{\text{h}} \hat{u}\|_{\hat{\mathrm{H}}_2^2} \leq C e^{-2\nu_1 t} \|\hat{u}\|_{\hat{\mathrm{H}}_2^2} \,, \quad t \geq 0 \,. \tag{22}$$

Let $V_0 \in H_2^2$, and let $V(t) = e^{t\mathcal{M}} V_0$ be the solution of (17) with initial data $V_0$. Then the Bloch wave transform of $V(t)$ can be decomposed as $\hat{V}(t) = \hat{E}_{\mathrm{c}}\hat{V}(t) + \hat{E}_{\mathrm{s}}\hat{V}(t)$. In view of (22), the stable part $\hat{E}_{\mathrm{s}}\hat{V}(t)$ converges exponentially to zero as $t \to +\infty$. On the other hand, by Lemma 3.1, the central part $\hat{v}(t) = \hat{E}_{\mathrm{c}}\hat{V}(t)$ satisfies

$$\hat{v}(\ell t^{-1/2}, x, t) = e^{-d\ell^2} \hat{v}_0(0, x) + \mathcal{O}(t^{-1/2}), \quad t \to +\infty. \tag{23}$$

To formulate this result more precisely, we introduce a few notations which will be useful to handle the nonlinear problem also.

For $\sigma \in (0, 1]$, we introduce the rescaling operator $\widehat{\mathcal{L}}_\sigma$ defined by

$$(\widehat{\mathcal{L}}_\sigma \hat{u})(\ell, x) = \hat{u}(\sigma \ell, x).$$

Note that the scaling does not act on the $x$ variable, only on the Bloch variable $\ell$. Since the domain for $\ell$ is finite, it will change with the scaling. Therefore, we introduce the function space

$$\mathcal{K}_{\sigma,\rho} = \left\{ \hat{u} : (-1/(2\sigma), 1/(2\sigma)) \times (0, 2\pi) \to \mathbb{C} \,\Big|\, \|\hat{u}\|_{\mathcal{K}_{\sigma,\rho}} < \infty \right\}, \tag{24}$$

where

$$\|\hat{u}\|^2_{\mathcal{K}_{\sigma,\rho}} = \sum_{j=0}^{2} \sum_{m=0}^{2} \int_{-1/(2\sigma)}^{1/(2\sigma)} \int_0^{2\pi} |\partial_\ell^j \partial_x^m \hat{u}(\ell, x)|^2 (1+\ell^2)^{2\rho} \,\mathrm{d}x \,\mathrm{d}\ell. \tag{25}$$

The polynomial weight in the Bloch variable $\ell$ will be used in section 5 to control the nonlinear terms. Indeed, if $\hat{u}, \hat{v} \in K_{\sigma,\rho}$ for some $\rho > 1/2$ and if

$$\hat{w}(\ell, x) = (\hat{u} \star \hat{v})(\ell, x) = \int_{-1/2\sigma}^{1/2\sigma} \hat{u}(\ell - \ell', x) \hat{v}(\ell', x) \,\mathrm{d}\ell',$$

then there exists $C > 0$ *independent of* $\sigma$ such that $\|\hat{w}\|_{K_{\sigma,\rho}} \leq C \|\hat{u}\|_{K_{\sigma,\rho}} \|\hat{v}\|_{K_{\sigma,\rho}}$.

It follows from the definitions that $\mathcal{K}_{1,\rho} = \hat{H}_2^2$ with equivalent norms (note, however, that the constants in the equivalence relation depend on $\rho$). Moreover, $\widehat{\mathcal{L}}_\sigma$ is an isomorphism from $\hat{H}_2^2$ to $\mathcal{K}_{\sigma,\rho}$, or more generally from $\mathcal{K}_{\sigma^{n-1},\rho}$ to $\mathcal{K}_{\sigma^n,\rho}$ for any $n \in \mathbb{N}$. In particular,

$$\|\widehat{\mathcal{L}}_\sigma \hat{f}\|_{\mathcal{K}_{\sigma^n,\rho}} \leq \sigma^{-1/2-2\rho} \|\hat{f}\|_{\mathcal{K}_{\sigma^{n-1},\rho}} \quad \text{and} \quad \|\widehat{\mathcal{L}}_\sigma^{-1} \hat{f}\|_{\mathcal{K}_{\sigma^{n-1},\rho}} \leq \sigma^{-3/2} \|\hat{f}\|_{\mathcal{K}_{\sigma^n,\rho}}, \tag{26}$$

where in the first estimate the additional factor $\sigma^{-2\rho}$ is due to the weight in the $\ell$-variable.

Using the definitions above together with Lemma 3.1, it is not difficult to verify that (23) can be written in the more precise form

$$\|\widehat{\mathcal{L}}_{1/\sqrt{t}} \hat{E}_{\mathrm{c}} \hat{V}(t) - e^{-d\ell^2} \hat{P}_{\mathrm{c}}(0) \hat{V}_0(0, \cdot)\|_{\mathcal{K}_{1/\sqrt{t},\rho}} \leq \frac{C}{\sqrt{t}} \|\hat{V}_0\|_{\hat{H}_2^2}, \quad t \geq 1.$$

On the other hand, using (22) and (26), we find

$$\|\widehat{\mathcal{L}}_{1/\sqrt{t}} \hat{E}_{\mathrm{s}} \hat{V}(t)\|_{\mathcal{K}_{1/\sqrt{t},\rho}} \leq t^{\rho+1/4} \|\hat{E}_{\mathrm{s}} \hat{V}(t)\|_{\hat{H}_2^2} \leq C e^{-\nu_1 t} \|\hat{V}_0\|_{\hat{H}_2^2}, \quad t \geq 0.$$

Thus we have

**Proposition 3.2** *Fix $\rho \geq 0$. If $\varepsilon > 0$ is small enough, the solution $\hat{V}(t)$ of (17) with initial data $\hat{V}_0$ satisfies*

$$\|\hat{V}(\ell t^{-1/2}, x, t) - A e^{-d\ell^2} \varphi(0, x)\|_{\mathcal{K}_{1/\sqrt{t},\rho}} \leq \frac{C}{t^{1/2}} \|\hat{V}_0\|_{\hat{H}_2^2},$$

*for all $t \geq 1$, where $A = \langle \psi_0, \hat{V}_0(0, \cdot) \rangle$. Moreover, there exists $\nu_1 > 0$ such that*

$$\|(\hat{E}_s \hat{V})(\ell t^{-1/2}, x, t)\|_{\mathcal{K}_{1/\sqrt{t},\rho}} \leq C e^{-\nu_1 t} \|\hat{V}_0\|_{\hat{H}_2^2},$$

*for all $t \geq 1$.*

To translate this result into the original variables, we observe that

$$\begin{aligned}
V(x,t) &= \frac{1}{\sqrt{t}} \int_{-\sqrt{t}/2}^{\sqrt{t}/2} e^{i\ell t^{-1/2} x} \hat{V}(\ell t^{-1/2}, x, t)\, d\ell \\
&= \frac{1}{\sqrt{t}} \Big( \int_{-\infty}^{\infty} A e^{-d\ell^2} \varphi(0, x) e^{i\ell t^{-1/2} x}\, d\ell + \mathcal{O}(1/\sqrt{t}) \Big) \\
&= \frac{A\, c_N}{\sqrt{4\pi d t}} e^{-x^2/(4dt)} \partial_x U_{\mathrm{per}}^-(x) + \mathcal{O}(1/t),
\end{aligned}$$

as $t \to +\infty$. This proves the analogue of Theorem 2.4 for the linearized system (17), see also Remark 2.5. Since $V(x,t)$ behaves for large times like a solution of the linear heat equation $\partial_t V = d\partial_x^2 V$ multiplied by the derivative of the steady state $U_{\mathrm{per}}^-$, we say that $V(x,t)$ converges "diffusively" to zero.

### 3.4 Nonlinear stability

In [Sch96, EWW97, Sch98a, Sch98b] it has been observed that spatially periodic equilibria with the above spectral properties are also nonlinearly stable with respect to spatially localized perturbations. The proof relies on the fact that the nonlinear terms are "asymptotically irrelevant" when compared to the linear diffusion. This is not obvious a priori, because the nonlinearity contains quadratic terms that are potentially dangerous, but this happen to be true due to nontrivial cancellations. Then a standard renormalization procedure [BK92] can be used to prove that the perturbations $V(x,t)$ converge diffusively to zero in the nonlinear case also. This is the statement of Theorem 2.4.

In section 5, we shall apply this renormalization procedure to a more difficult problem, namely the stability of the modulated fronts $U_{\mathrm{mf}}$. The proof of Theorem 2.4 is a particular case of this more complicated argument, see Remark 5.1 below. So, for the sake of brevity, we shall not discuss Theorem 2.4 any longer, and refer to section 5 for more details.

## 4 Construction of modulated fronts

In this section, we follow closely [SS99] where modulated pulse solutions are constructed for general reaction-diffusion systems. However, since the assumptions of [SS99] are not exactly satisfied in our case, we give here a short proof of Theorem 2.6. Throughout this section we refer to [PSS97], [SS01b] for the functional analytic background and the relation between temporal and spatial dynamics.

## 4.1 The idea

As already said, the modulated fronts are time-periodic in a frame moving with a speed $c$ close to the speed $c_0$ of the original front. In this frame, we shall denote the spatial variable by $\xi = x - ct$, to distinguish it from $y = x - c_0 t$. Equation (3) then becomes

$$\partial_t U = L(\partial_\xi) U + c \partial_\xi U + N(U) \,, \tag{27}$$

and we look for solutions $U(\xi, t)$ that are periodic in their second argument.

The key idea in the construction of the modulated fronts is *spatial dynamics*, i.e., system (27) will be considered as an evolution system for $U$ with respect to the variable $\xi = x - ct$, in a space of functions $U(\xi, t)$ which are periodic in $t$. The use of spatial dynamics goes back at least to [Kir82] and is nowadays a well established method for the construction of fronts and pulses.

Writing (27) as a first order system with respect to $\xi$ yields

$$\partial_\xi V = \partial_t L_1 V + G(V), \tag{28}$$

with

$$V = (v_1, v_2, v_3, v_4, v_5, v_6) = (u, \partial_\xi u, v, \partial_\xi v, \partial_\xi^2 v, \partial_\xi^3 v)$$

and

$$L_1 V = \begin{pmatrix} 0 \\ +v_1 \\ 0 \\ 0 \\ 0 \\ -v_3 \end{pmatrix}, \quad G(V) = \begin{pmatrix} v_2 \\ -\frac{1}{2}(v_1 - c_0)(1 - v_1^2) - v_3 - c v_2 \\ v_4 \\ v_5 \\ v_6 \\ (-1 + \alpha) v_3 + c v_4 - 2 v_5 - v_3^3 - \gamma(1 - v_1^2) v_3 \end{pmatrix}.$$

For a fixed $s \geq 2$, this system will be considered in the infinite-dimensional phase space

$$\begin{aligned} \mathcal{X} = {} & H^s_{\text{per}}(0, T) \times H^{s-1/2}_{\text{per}}(0, T) \times H^s_{\text{per}}(0, T) \\ & \times H^{s-1/4}_{\text{per}}(0, T) \times H^{s-1/2}_{\text{per}}(0, T) \times H^{s-3/4}_{\text{per}}(0, T) \end{aligned}$$

where $T = (2\pi)/c$. In the spatial dynamics formulation (28) we will easily find equilibria $V_\pm$ corresponding to the equilibria $U_\pm$ of (3), periodic solutions $V^\pm_{\text{per}}$ which correspond to the spatially periodic equilibria $U^\pm_{\text{per}}$ of (3), and we also find a trivial heteroclinic connection $V_h$ between $V_-$ and $V_+$ which corresponds to the trivial front $U_h$ of (3).

The linearization of the spatial dynamics formulation (28) around the trivial equilibria $V_\pm$ possesses two eigenvalues close to the imaginary axis and infinitely many eigenvalues with positive real part and infinitely many eigenvalues with negative real part. Hence there will be a two-dimensional center manifold $W_c(V_\pm)$, an infinite-dimensional stable manifold $W_s(V_\pm)$ and an infinite-dimensional unstable manifold $W_u(V_\pm)$ of $V_\pm$. For $\alpha > 0$ sufficiently small, the periodic solution $V_{\text{per}}^\pm$ lies in the two-dimensional center manifold $W_c(V_\pm)$.

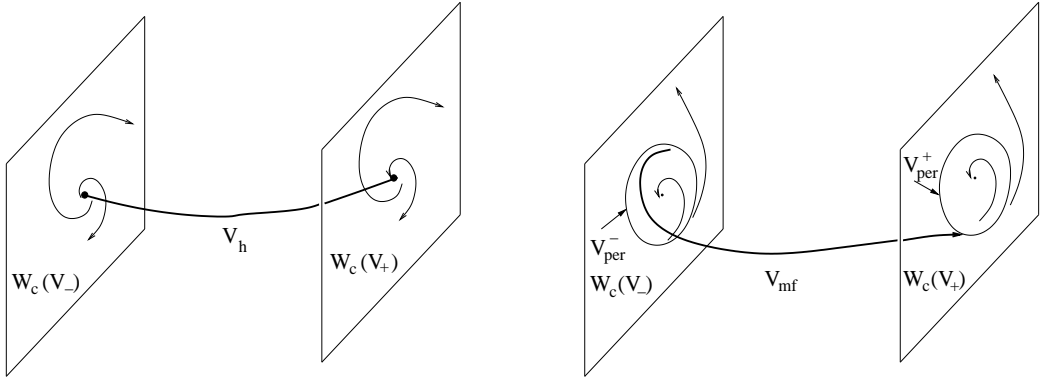

Figure 3: The original front and the modulated front. The vertical planes symbolize the two-dimensional center manifolds at $V_-$ and $V_+$ for the spatial dynamics formulation (28). Left picture ($\alpha < 0$): $V_\pm$ are unstable equilibria, and the solid line represents the trivial front $V_h$. Right picture ($\alpha > 0$): $V_\pm$ are stable, $V_{\text{per}}^\pm$ are unstable, and the solid line is the modulated front we want to construct.

The modulated front solutions we are interested in will be found in the intersection of the *center-unstable* manifold $W_{cu}(V_-)$ of $V_-$ and the *stable* manifold $W_s(V_{\text{per}}^+)$ of $V_{\text{per}}^+$. The reason is as follows. Since $V_{\text{per}}^+$ is unstable in the two-dimensional center manifold $W_c(V_+)$ of $V_+$, solutions converging towards $V_{\text{per}}^+$ for $\xi \to \infty$ have to be in the stable manifold $W_s(V_{\text{per}}^+)$ of $V_{\text{per}}^+$. On the other hand for $\xi \to -\infty$ the periodic solution $V_{\text{per}}^-$ attracts all solutions in the two-dimensional center manifold $W_c(V_-)$ of $V_-$, except for $V_-$ itself. Therefore, to converge towards $V_{\text{per}}^-$ for $\xi \to -\infty$ it is sufficient to be in the center-unstable manifold $W_{cu}(V_-)$ of $V_-$. Thus, we will almost be done if we show that the center-unstable manifold $W_{cu}(V_-)$ of $V_-$ and the stable manifold $W_s(V_{\text{per}}^+)$ of $V_{\text{per}}^+$ intersect.

As is explained in [SS99], this is the case only if the parameter $c$ is chosen appropriately. To cope with this problem, we consider (27) as a dynamical system for the pair $(c, V) \in \mathbb{R} \times \mathcal{X}$, where $c$ obeys the trivial equation $\dot{c} = 0$, and we look for an intersection of the corresponding manifolds in the extended space $\mathbb{R} \times \mathcal{X}$. We proceed in three steps:

i) In Lemma 4.1 we prove that, for $\alpha = 0$, the center-unstable manifold $W_{cu}(I \times \{V_-\})$ of the family of fixed points $I \times \{V_-\}$ and the stable manifold $W_s(I \times \{V_+\})$ of the family of fixed points $I \times \{V_+\}$ intersect transversally in the extended space $\mathbb{R} \times \mathcal{X}$ (here $I \subset \mathbb{R}$ is a closed neighborhood of $c_0$). Moreover, there exists exactly one $c \in I$ such that $W_{cu}(V_-)$ intersects $W_s(V_+)$, namely the velocity $c_0$ of the trivial front.

ii) In Lemma 4.2 we reformulate the bifurcation of the spatially periodic equilibria $U_{\text{per}}^\pm$ for (3) at $\alpha = 0$ as the bifurcation of time-periodic solutions $V_{\text{per}}^\pm$ for the spatial dynamics formulation (28).

iii) By a perturbation argument, we show that for $\alpha > 0$ sufficiently small the center-unstable manifold $W_{cu}(I \times \{V_-\})$ and the stable manifold $W_s(I \times \{V_{\text{per}}^+\})$ intersect transversally in the extended space $\mathbb{R} \times \mathcal{X}$. Since $V_{\text{per}}^-$ is attractive for $\xi \to -\infty$, this implies the existence of a modulated front solution with a velocity $c$ close to the original velocity $c_0$.

## 4.2 Transversality in extended space

For all $\alpha, c \in \mathbb{R}$, system (28) possesses the fixed points $V_\pm = (\pm 1, 0, 0, 0, 0, 0)$. Moreover, for all $\alpha$ and for $c = c_0$ we have the heteroclinic connection

$$V(\xi, t) = V_h(\xi) = (h(\xi), h'(\xi), 0, 0, 0, 0)$$

between $V_-$ and $V_+$, i.e. $\lim_{\xi \to \pm\infty} V_h(\xi) = V_\pm$. These solutions lie in the invariant subspaces $\mathcal{X}_{00} \subset \mathcal{X}_0 \subset \mathcal{X}$, where $\mathcal{X}_0$ is the six-dimensional subspace that consists of $t$-independent solutions satisfying $\partial_\xi V = G(V)$, and $\mathcal{X}_{00}$ is the two-dimensional subspace $\mathcal{X}_{00} = \{V \in \mathcal{X}_0 \mid v_3 = v_4 = v_5 = v_6 = 0\}$. We now prove the following result.

**Lemma 4.1** *Fix $c_0 \in (0,1)$. There exist $\alpha_0 > 0$ and $\gamma_0 > 0$ such that for all $\alpha \in \mathbb{R}$ with $|\alpha| \leq \alpha_0$ and all $\gamma \in [0, \gamma_0)$ the following holds. There exists a closed neighborhood $I$ of $c_0$ such that the center-unstable manifold $W_{cu}(I \times \{V_-\})$ of the family of fixed points $I \times \{V_-\}$ and the stable manifold $W_s(I \times \{V_+\})$ of the family of fixed points $I \times \{V_+\}$ intersect transversally in the extended space $\mathbb{R} \times \mathcal{X}$.*

**Proof. a) Existence of invariant manifolds.** The linearization of (28) at $V_\pm$ is given by

$$\partial_\xi V = \Lambda^\pm V \qquad \text{with} \qquad \Lambda^\pm = \partial_t L_1 + DG(V_\pm) \,. \tag{29}$$

Since $V(\xi, t)$ is $T$-periodic in its second argument with $T = 2\pi/c$, it is natural to decompose it using Fourier series:

$$V(\xi, t) = \sum_{m \in \mathbb{Z}} V_m(\xi) e^{im\omega t} \,, \quad \text{where } \omega = c \,. \tag{30}$$

The space $\mathcal{X}$ then splits into a direct sum $\oplus_{m \in \mathbb{Z}} \mathcal{X}_m$, where for each $m \in \mathbb{Z}$ the six-dimensional subspace $\mathcal{X}_m$ is invariant under (29). Restricting (29) to $\mathcal{X}_m$ yields

$$\partial_\xi V_m = \Lambda_m^\pm V_m \qquad \text{with} \qquad \Lambda_m^\pm = im\omega L_1 + DG(V_\pm) \,,$$

or explicitly

$$\Lambda_m^\pm = \begin{pmatrix} 0 & 1 & 0 & 0 & 0 & 0 \\ i\omega m + (1 \mp c_0) & -c & -1 & 0 & 0 & 0 \\ 0 & 0 & 0 & 1 & 0 & 0 \\ 0 & 0 & 0 & 0 & 1 & 0 \\ 0 & 0 & 0 & 0 & 0 & 1 \\ 0 & 0 & -i\omega m - 1 + \alpha & c & -2 & 0 \end{pmatrix} .$$

In order to solve the associated eigenvalue problem $\Lambda_m^\pm V_m = \lambda V_m$, we analyze the condition $\det(\Lambda_m^\pm - \lambda \mathrm{Id}) = 0$. This equation splits into two parts and reduces to the system

$$\lambda(\lambda + c) - (1 \mp c_0) - \mathrm{i}\omega m = 0 \quad \text{or} \quad -(1+\lambda^2)^2 + c\lambda - \mathrm{i}\omega m + \alpha = 0, \qquad (31)$$

which can also be obtained directly from the temporal formulation (27), see [HS99].

For $\alpha = 0$, we have exactly two eigenvalues on the imaginary axis, i.e. $\lambda = \mathrm{i}$ for $m = 1$ and $\lambda = -\mathrm{i}$ for $m = -1$. Moreover there are infinitely many stable and infinitely many unstable directions, and the associated eigenvalues are contained in two sectors. see fig.2 in [SS99]. We define projections $P_s^\pm$, $P_u^\pm$ and $P_c^\pm$ on the stable, unstable and central part of the linear operators $\Lambda^\pm$. This can be done for each $\Lambda_m^\pm$ with $m \in \mathbb{Z}$, and so these operators are well-defined. Then the space $\mathcal{X}$ splits into

$$\mathcal{X} = \mathcal{X}_c^\pm \oplus \mathcal{X}_s^\pm \oplus \mathcal{X}_u^\pm \qquad \text{with} \qquad \mathcal{X}_j^\pm = P_j^\pm \mathcal{X}.$$

The restrictions $\Lambda_j^\pm = \Lambda^\pm|_{\mathcal{X}_j^\pm}$ generate analytic semigroups in $\mathcal{X}_j^\pm$ satisfying

$$\begin{aligned}
\|e^{\Lambda_c^\pm \xi}\|_{\mathcal{X}_c^\pm \to \mathcal{X}_c^\pm} &\leq C & &\text{for all } \xi \in \mathbb{R}, \\
\|e^{\Lambda_s^\pm \xi}\|_{\mathcal{X}_s^\pm \to \mathcal{X}_s^\pm} &\leq Ce^{-\beta\xi} & &\text{for all } \xi \geq 0, \\
\|e^{\Lambda_u^\pm \xi}\|_{\mathcal{X}_u^\pm \to \mathcal{X}_u^\pm} &\leq Ce^{-\beta|\xi|} & &\text{for all } \xi \leq 0,
\end{aligned}$$

for some constants $C, \beta > 0$. This can be proved by a rescaling argument, see Lemma 3.1 in [SS99]. In the present case, the appropriate scaling is:

$$\begin{aligned}
&V_{1,m} = \tilde{V}_{1,m},\ V_{2,m} = (1+m^2)^{1/4}\tilde{V}_{2,m},\ V_{3,m} = \tilde{V}_{3,m}, \\
&V_{4,m} = (1+m^2)^{1/8}\tilde{V}_{4,m},\ V_{5,m} = (1+m^2)^{1/4}\tilde{V}_{5,m}, \quad \text{and} \quad V_{6,m} = (1+m^2)^{3/8}\tilde{V}_{6,m}.
\end{aligned}$$

By construction, $V \in \mathcal{X}$ if and only if $\tilde{V} \in H^s_{\mathrm{per}}(0,T)$.

Since the nonlinear term $\mathcal{N}(V) = \partial_t L_1 V + G(V) - \Lambda^\pm V$ is smooth from $\mathcal{X}$ into $\mathcal{X}$, [Hen81] guarantees the existence of the following local smooth invariant manifolds in $\mathcal{X}$: the center manifolds $W_c(V_\pm)$ which are tangent to $\mathcal{X}_c^\pm$ at $V_\pm$, the stable manifolds $W_s(V_\pm)$ tangent to $\mathcal{X}_s^\pm$, the unstable manifolds $W_u(V_\pm)$ tangent to $\mathcal{X}_u^\pm$, the center-stable manifolds $W_{cs}(V_\pm)$ tangent to $\mathcal{X}_{cs}^\pm = \mathcal{X}_c^\pm \oplus \mathcal{X}_s^\pm$, and finally the center-unstable manifolds $W_{cu}(V_\pm)$ tangent to $\mathcal{X}_{cu}^\pm = \mathcal{X}_c^\pm \oplus \mathcal{X}_u^\pm$. In addition, as in [SS99] the center-unstable manifold $W_{cu}(V_-)$ and the stable manifold $W_s(V_+)$ can be extended to a whole neighborhood of the trivial heteroclinic connection $V_h$. These global manifolds will still be denoted by $W_{cu}(V_-)$ and $W_s(V_+)$.

**b) Transversality.** The manifolds $W_{cu}(V_-)$ and $W_s(V_+)$ do not intersect transversally in the space $\mathcal{X}$, see e.g. [SS99]. To obtain transversality the space has to be extended to $\mathbb{R} \times \mathcal{X}$ by adding one direction corresponding to the velocity $c$.

The transversality of $W_{cu}(I \times \{V_-\})$ and $W_s(I \times \{V_+\})$ in the extended space $\mathbb{R} \times \mathcal{X}$ for $(c, V)$ means that the tangent spaces to these manifolds at a point $(c_0, V_h(0))$ span $\mathbb{R} \times \mathcal{X}$, namely

$$\mathbb{R} \times \mathcal{X} = \{u + v \mid u \in T_{(c_0, V_h(0))}(W_{cu}(I \times \{V_-\})),\ v \in T_{(c_0, V_h(0))}(W_s(I \times \{V_+\}))\}.$$

Such a situation is robust under small perturbations, and hence slightly perturbed manifolds will also possess a non empty intersection.

To prove transversality, we follow again [SS99] and consider the linearization of (28) around the heteroclinic connection $V_h$:

$$\partial_\xi \tilde{V} = \partial_t L_1 \tilde{V} + DG(V_h)\tilde{V} \ . \tag{32}$$

We first observe that, for $c = c_0$ and $\alpha = 0$, (32) has a unique nontrivial solution $\tilde{V} \in C_b^0(\mathbb{R}, \mathcal{X})$ that converges exponentially to zero as $\xi \to +\infty$, namely $\tilde{V} = \partial_\xi V_h$. Indeed, if $\tilde{V}(\xi, t)$ is any such solution and if $\tilde{U}(\xi, t) = (\tilde{v}_1(\xi, t), \tilde{v}_3(\xi, t))$, then by construction $\tilde{U}(\xi, t)$ is a time-periodic solution of the linearization of (27) around $U_h$:

$$\partial_t \tilde{U} = L(\partial_\xi)\tilde{U} + c_0 \partial_\xi \tilde{U} + DN(U_h)\tilde{U} \ .$$

Moreover, for $\beta > 0$ sufficiently small, $\tilde{U} \in C^0_{\text{per}}([0,T], X_\beta)$ where $X_\beta = L^2(\mathbb{R}, e^{\beta\xi} \, \mathrm{d}\xi)$. By Proposition 2.2, the spectrum of the linear operator $L(\partial_\xi) + c_0 \partial_\xi + DN(U_h)$ in $X_\beta$ is strictly contained in the left half-plane, except for a simple eigenvalue 0 with eigenvector $\partial_\xi U_h$. Thus we necessarily have $\tilde{U} = \partial_\xi U_h$, hence also $\tilde{V} = \partial_\xi V_h$.

Next we observe that, for any $m \in \mathbb{Z}$, the six-dimensional subspace $\mathcal{X}_m$ is left invariant by the time-independent equation (32). This property allows to study the transversality of the manifolds $W_{cu}(V_-)$ and $W_s(V_+)$ in each subspace $\mathcal{X}_m$ separately, since the tangent spaces $T_{V_h(\xi)}(W_{cu}(V_-))$ (resp. $T_{V_h(\xi)}(W_s(V_+))$) for different values of $\xi$ are mapped onto each other by solutions of (32). Remark that (32) defines a well-posed evolution in each subspace $\mathcal{X}_m$, but not in the whole space $\mathcal{X}$.

For any $m \in \mathbb{Z}$, let

$$E_{m,-}^{cu} = T_{V_h(0)}(W_{cu}(V_-)) \cap \mathcal{X}_m \ , \quad E_{m,+}^s = T_{V_h(0)}(W_s(V_+)) \cap \mathcal{X}_m \ .$$

From (31), it is not difficult to verify that $\dim(E_{m,-}^{cu}) = \dim(E_{m,+}^s) = 3$ for all $m \in \mathbb{Z}$. Moreover, if $m \neq 0$, the argument above shows that $E_{m,-}^{cu} \cap E_{m,+}^s = \{0\}$, so that $\mathcal{X}_m = E_{m,-}^{cu} \oplus E_{m,+}^s$. If $m = 0$, then $\dim(E_{0,-}^{cu} \cap E_{0,+}^s) = 1$, so that $E_{0,-}^{cu} + E_{0,+}^s$ has codimension 1 in $\mathcal{X}_0$. Summarizing, we have shown that there exists $\Psi \in \mathcal{X}_0$ such that

$$\mathcal{X} = \{u + v \mid u \in T_{V_h(0)}(W_{cu}(V_-)), \ v \in T_{V_h(0)}(W_s(V_+))\} + \text{span}\{\Psi\}. \tag{33}$$

If we further restrict system (28) to the two-dimensional invariant subspace $\mathcal{X}_{00}$, we obtain the simple differential equation

$$\partial_\xi v_1 = v_2 \ , \quad \partial_\xi v_2 = -cv_2 - \frac{1}{2}(v_1 - c_0)(1 - v_1^2)$$

which governs the travelling waves of the Chaffee-Infante equation. For $c = c_0$ the one-dimensional manifolds $W_{cu}(V_-)|_{\mathcal{X}_{00}}$ and $W_s(V_+)|_{\mathcal{X}_{00}}$ intersect (non transversally) along the heteroclinic orbit $V_h|_{\mathcal{X}_{00}}$. Thus the vector $\Psi$ can be chosen as the unit normal to $V_h'(0)$ in $\mathcal{X}_{00}$. Moreover, it is easy to see that a variation with respect to the parameter $c$ shifts the stable and unstable manifolds through each other, namely $W_{cu}(I \times \{V_-\})|_{\mathcal{X}_{00}}$ and $W_s(I \times \{V_+\})|_{\mathcal{X}_{00}}$ intersect transversally in the three-dimensional space $\mathbb{R} \times \mathcal{X}_{00}$. Together with (33), this implies that $W_{cu}(I \times \{V_-\})$ and $W_s(I \times \{V_+\})$ intersect transversally in the extended space $\mathbb{R} \times \mathcal{X}$. □

## 4.3 Bifurcation of Turing patterns

The following lemma is the analogue of Theorem 2.3 for the spatial formulation (28).

**Lemma 4.2** *Fix $c_0 \in (0,1)$. There exist $\varepsilon_0 > 0$ and $\gamma_0 > 0$ such that for all $\varepsilon \in (0, \varepsilon_0)$, all $\gamma \in (0, \gamma_0)$ and all $c > 0$ equation (28) with $\alpha = \varepsilon^2$ has two families $\{V_{\text{per}}^{\pm}(\xi + ct + \xi_0) \mid \xi_0 \in \mathbb{R}\}$ of periodic solutions satisfying $V_{\text{per}}^{\pm}(x) = V_{\text{per}}^{\pm}(x + 2\pi)$ and $\|V_{\text{per}}^{\pm}(x) - V_{\pm}\|_{\mathbb{R}^6} = \mathcal{O}(\varepsilon)$.*

**Remark 4.3** *We may prove Lemma 4.2 directly as follows. For the eigenvalue problem (31) in the proof of Lemma 4.1 all eigenvalues $\lambda$ except for $m = \pm 1$ are uniformly bounded away from the imaginary axis and contained in two sectors. For $\alpha = 0$ and $\omega = c$ we have two eigenvalues on the imaginary axis, namely $\lambda = \mathrm{i}$ for $m = 1$ and $\lambda = -\mathrm{i}$ for $m = -1$. With $\alpha$ changing sign from $-$ to $+$ the real parts of these eigenvalues change sign from $+$ to $-$. To see this we insert $\lambda = \mathrm{i} + \mu$ into the second equation of (31) and obtain*

$$-(\mu^2 + 4\mathrm{i}\mu - 4)\mu^2 + c\mu + \alpha = 0$$

*and therefore*

$$\mu = -\frac{\alpha}{c} + \mathcal{O}(\alpha^2).$$

*Now the $\mathcal{S}^1$–equivariant Hopf-bifurcation theorem applies [GSS88], which shows the required result. On the two-dimensional center manifolds $W_c(V_{\pm})$ we have an unstable origin for $\alpha < 0$, and for $\alpha > 0$ a stable origin and an unstable periodic solution $V_{\text{per}}^{\pm}$ with, e.g., $V_{\text{per}}^{-}(\xi, t) = V_{-} + \mathcal{O}(\sqrt{\alpha}) \left( e^{\mathrm{i}ct} + c.c. \right) + \mathcal{O}(\alpha)$.*

## 4.4 Bifurcation of the modulated front

In order to complete the proof of Theorem 2.6 it remains to establish point iii) of the program of section 4.1.

As a consequence of the continuity with respect to $\alpha = \varepsilon^2$, the stable manifold $W_s(I \times \{V_{\text{per}}^+\})$ of the bifurcating periodic solution $V_{\text{per}}^+$ stays $\mathcal{O}(\varepsilon)$-close to $W_s(I \times \{V_+\})$ for $\xi \in \mathbb{R}^+$, see for instance [GH83, SS99]. Therefore, the intersection of $W_{cu}(I \times \{V_-\})$ and $W_s(I \times \{V_{\text{per}}^+\})$ is still transversal for $\varepsilon > 0$ sufficiently small. By [SS99, Lemma 3.9], there exists $c = c_0 + \mathcal{O}(\varepsilon^2)$ such that $W_{cu}(V_-)$ and $W_s(V_{\text{per}}^+)$ intersect, which follows from

$$\text{dist}(W_{cu}(V_-), W_s(V_{\text{per}}^+)) = \mathcal{O}(|c - c_0|) + \mathcal{O}(\varepsilon^2). \tag{34}$$

This distance can be measured in $\mathcal{X}_{00}$ by the distance of the intersection points of these manifolds with the $v_2$-axis in $\mathcal{X}_{00}$. The dynamics in $\mathcal{X}_{00}$ is mainly described by the first equation in (1). The linear term $+v$ in the first equation of (1) is $\mathcal{O}(\varepsilon)$. However, the $\mathcal{O}(\varepsilon)$-terms in $v$ belong to $\mathcal{X}_1$, and by quadratic interaction they couple $\mathcal{O}(\varepsilon^2)$ back to $\mathcal{X}_0$. In all other $\mathcal{X}_j$ with $j \neq 0$ due to transversality for $\varepsilon = 0$ the tangent spaces of $W_{cu}(V_-) \cap \mathcal{X}_j$ and $W_s(V_{\text{per}}^+) \cap \mathcal{X}_j$ also span $\mathcal{X}_j$ for $\varepsilon > 0$ sufficiently small and (34) follows.

A solution of (28) in the intersection converges for $\xi \to \infty$ to $V_{\text{per}}^+$ with some rate $e^{-\beta \xi}$ where $\beta > 0$ can be chosen independent of $\varepsilon$. For $\xi \to -\infty$ all solutions in the center-unstable manifold of $V_-$ converge with some exponential rate to the two-dimensional center manifold of $V_-$. On this manifold except for the origin $V_-$ all solutions converge for $\xi \to -\infty$ with some rate $\mathcal{O}(e^{-\mathcal{O}(\varepsilon)|\xi|})$ towards the unstable periodic solution $V_{\mathrm{per}}^-$, cf. figure 3.

It remains to prove that $W_s(V_{\mathrm{per}}^+)$ is not connected with $W_u(V_-)$. This can be shown as in [SS99] by remarking that the distance between $W_s(V_+)$ and $W_u(V_-)$ near $V_h(0)$ is of order $\varepsilon^2$, since both manifolds are smooth in $\alpha = \varepsilon^2$. One the other hand, since $\mathrm{dist}(V_{\mathrm{per}}^+, V_+) \geq C\varepsilon$, one can verify by looking at $\mathcal{X}_1$ that the distance between $W_s(V_{\mathrm{per}}^+)$ and $W_s(V_+)$ near $V_h(0)$ is also bounded from below by $C\varepsilon$. Thus, if $\varepsilon > 0$ is small enough, there cannot be an intersection between $W_s(V_{\mathrm{per}}^+)$ and $W_u(V_-)$. In detail

$$\mathrm{dist}(W_u(V_-), W_s(V_{\mathrm{per}}^+)) \geq \mathrm{dist}(W_s(V_+), W_s(V_{\mathrm{per}}^+)) - \mathrm{dist}(W_u(V_-), W_s(V_+))$$
$$\geq C_1\varepsilon - C_2\varepsilon^2 \geq C_1\varepsilon/2$$

for $\varepsilon > 0$ sufficiently small.

Summarizing, we have found modulated front solutions $V_{\mathrm{mf}}(\xi, t)$ of (28) connecting $V_{\mathrm{per}}^-$ to $V_{\mathrm{per}}^+$ with a velocity $c = c_0 + \mathcal{O}(\varepsilon^2)$. Moreover, $\sup_{\xi, t \in \mathbb{R}} |V_{\mathrm{mf}}(\xi, t) - V_h(\xi)| \leq C\varepsilon$. Setting $U_{\mathrm{mf}}(\xi, \xi + ct) = (v_{\mathrm{mf},1}(\xi, t), v_{\mathrm{mf},3}(\xi, t))$ we obtain a modulated front $U_{\mathrm{mf}}$ satisfying the conclusions of Theorem 2.6. □

## 4.5 Remarks on case II and case III

The following two remarks are an adaption of the theory given in [SS01a] to our situation.

**Remark 4.4** *In case II there is no Hopf-bifurcation at $V_-$, but the same argument as above shows the transversality of $W_s(I \times \{V_{\mathrm{per}}^+\})$ and $W_u(I \times \{V_-\})$ in $\mathbb{R} \times \mathcal{X}$. Therefore, for $c$ close to $c_0$, there are modulated front solutions $U(x, t) = U_{\mathrm{mf}}(x - ct, x)$ of (3) which are $2\pi$-periodic in the second argument and satisfy*

$$\lim_{\xi \to \infty} U_{\mathrm{mf}}(\xi, x) = U_{\mathrm{per}}^+(x) \quad \text{and} \quad \lim_{\xi \to -\infty} U_{\mathrm{mf}}(\xi, x) = U_- \ .$$

**Remark 4.5** *In case III there is no Hopf-bifurcation at $V_+$. There is still a transversal intersection between $W_{cu}(I \times \{V_-\})$ and $W_s(I \times \{V_+\})$. But for $c = c_0$ we know that $W_{cu}(V_-)$ and $W_s(V_+)$ intersect along the heteroclinic connection $V_h$. By uniqueness there is no other connection between $W_{cu}(V_-)$ and $W_s(V_+)$ and so no modulated front in this case.*

# 5 The nonlinear stability proof

This section contains the proof of Theorem 2.8. As already said, the linearization around the modulated fronts has continuous spectrum up to the imaginary axis, because the periodic pattern $U_{\mathrm{per}}^-$ behind the front is only diffusively stable. In addition, we have an embedded zero eigenvalue which is responsible for the shift along the family of modulated fronts. This coexistence of discrete and continuous spectrum on the imaginary axis is the main technical difficulty in proving the stability of the modulated fronts. To handle this problem, we follow the approach developed in [ES02]. The basic idea is to assume that initially perturbations decay with an exponential rate as $x \to +\infty$. Even if the steady state ahead of the front is weakly unstable (in our case, it is not), such perturbations will be "overtaken" by the propagating front before they have enough time to grow. Once they are far behind the front, they vanish diffusively because the periodic pattern $U_{\mathrm{per}}^-$ is stable. In case II, the steady state $U_-$ behind the front is exponentially stable, so that the perturbations decay to zero exponentially as $t \to \infty$, see section 6.

Throughout this section, we fix $\varepsilon > 0$ and $\beta > 0$ small enough so that the conclusions of Proposition 2.2, Theorem 2.3, Theorem 2.6 and Lemma 3.1 hold. Unless explicitly stated, the constants below will depend on $\varepsilon$ and $\beta$. The small parameters we shall use are the size of the initial data, which we denote by $\delta$, and the rescaling parameter $\sigma$ which will be introduced in section 5.3.

## 5.1 The idea

The zero eigenvalue of the linearized operator will lead to a spatial shift of the front. Therefore, we write solutions $U(x,t)$ of (3) as a shift of the modulated front plus some perturbation $V(x,t)$, namely

$$U(x,t) = U_{\mathrm{mf}}(x - ct - q(t), x) + V(x,t), \quad x \in \mathbb{R}, \quad t \geq 0, \tag{35}$$

where $q(t) \in \mathbb{R}$ is the spatial shift. This representation is clearly not unique, but we shall use this freedom below to impose a condition on $V(x,t)$ that will determine $q(t)$ uniquely. The perturbation $V(x,t)$ satisfies

$$\partial_t V(x,t) = L(\partial_x) V(x,t) + DN(U_{\mathrm{mf}}(x - ct - q(t), x)) V(x,t) \tag{36}$$
$$+ N_1(U_{\mathrm{mf}}(x - ct - q(t), x), V(x,t)) + \dot{q}(t) \partial_1 U_{\mathrm{mf}}(x - ct - q(t), x),$$

with $N_1(U_{\mathrm{mf}}, V) = N(U_{\mathrm{mf}} + V) - N(U_{\mathrm{mf}}) - DN(U_{\mathrm{mf}})V = \mathcal{O}(V^2)$. To analyze the behavior of $V(x,t)$ near $x = +\infty$, it is convenient to go to the comoving frame and to use exponential weights [Sat77]. Let $\beta \in (0, \beta_0)$, where $\beta_0$ is defined in Proposition 2.2. We introduce the weighted variable $W$ defined by

$$W(\xi,t) = V(\xi + ct, t) e^{\beta \xi}, \quad \xi \in \mathbb{R}, \quad t \geq 0. \tag{37}$$

The equation for $W$ is

$$\partial_t W(\xi,t) = L_{\beta,c} W(\xi,t) + DN(U_{\mathrm{mf}}(\xi - q(t), \xi + ct)) W(\xi,t)$$
$$+ N_2(U_{\mathrm{mf}}(\xi - q(t), \xi + ct), V(\xi + ct, t), W(\xi,t)) \tag{38}$$
$$+ \dot{q}(t) \partial_1 U_{\mathrm{mf}}(\xi - q(t), \xi + ct) e^{\beta \xi},$$

where
$$L_{\beta,c} = L(\partial_\xi - \beta) + c(\partial_\xi - \beta) \quad \text{and} \quad N_2(U_{\mathrm{mf}}, V, W) = \mathcal{O}(VW).$$

Eq. (38) is coupled to (36) through the nonlinear terms $N_2(U_{\mathrm{mf}}, V, W)$, which are in fact *linear* with respect to $W$, see section 5.2. Although $V$ and $W$ are simply related via (37), we find it convenient to keep both variables in the sequel. We shall assume that $W \in H^2(\mathbb{R})$ and $V \in \mathrm{H}_2^2 = \{f \in H^2(\mathbb{R}) \mid x^2 f \in H^2(\mathbb{R})\}$, which is equivalent to requiring that $(x^2 + e^{\beta x})V$ lies in $H^2(\mathbb{R})$. This defines a weighted perturbation space which is an algebra for the product of functions.

To understand the spectrum of the linear operator in (38) we write

$$L_{\beta,c} W + DN(U_{\mathrm{mf}})W = L_{\beta,c_0} W + DN(U_h)W + (c - c_0)\partial_\xi W + L_\Delta(U_{\mathrm{mf}})W \;,$$

where $L_\Delta(U_{\mathrm{mf}})W = DN(U_{\mathrm{mf}})W - DN(U_h)W$. Since $|c - c_0| = \mathcal{O}(\varepsilon^2)$ and $|DN(U_{\mathrm{mf}}) - DN(U_h)| = \mathcal{O}(\varepsilon + |q(t)|)$ by Theorem 2.6, we see that, provided $\varepsilon$ and $q$ are small enough, the time-dependent operator $L_{\beta,c} + DN(U_{\mathrm{mf}})$ can be considered as a perturbation of the simpler operator $\Lambda_\beta = L_{\beta,c_0} + DN(U_h)$. By Proposition 2.2, the spectrum of $\Lambda_\beta$ is strictly contained in the left half-plane, except for a zero eigenvalue which is due to translation invariance and is not affected by the exponential weight. Let $\Pi_c : H^2 \to H^2$ be the spectral projection onto the one-dimensional eigenspace of $\Lambda_\beta$ corresponding to the eigenvalue zero. We also set $\Pi_s = \mathbf{1} - \Pi_c$. If $W$ is defined by (37), we impose the condition $\Pi_c W(t) = 0$ for all $t$, which amounts to fixing the shift $q(t)$ in (35). Under this assumption on $W$, a standard argument shows that any solution $U(x, t)$ that stays sufficiently close to the modulated front $U_{\mathrm{mf}}(x - ct, x)$ can be decomposed in a unique way according to (35).

The evolution system for $W$ and $q$ reads

$$\begin{aligned}
\partial_t W &= \Pi_s \left( L_{\beta,c} W + DN(U_{\mathrm{mf}})W + N_2(U_{\mathrm{mf}}, V, W) + \dot{q}\partial_1 U_{\mathrm{mf}} e^{\beta \xi} \right) \;, \\
\dot{q} &= -(\Pi_c(\partial_1 U_{\mathrm{mf}} e^{\beta \xi}))^{-1} \Pi_c((c - c_0)\partial_\xi W + L_\Delta(U_{\mathrm{mf}})W + N_2(U_{\mathrm{mf}}, V, W)) \;,
\end{aligned} \quad (39)$$

where $\Pi_c(\partial_1 U_{\mathrm{mf}} e^{\beta \xi}) = 1 + \mathcal{O}(\varepsilon + |q|)$. By Proposition 2.2, there exists $\nu = \nu(\beta) > 0$ such that $\|e^{\Lambda_\beta t} \Pi_s\|_{H^2} = \mathcal{O}(e^{-2\nu t})$ as $t \to \infty$. If $V(x,t)$ and $q(t)$ are bounded and remain sufficiently small, due to $\dot{q} = \mathcal{O}(W)$ we expect that the solutions of (39) will satisfy $|\dot{q}(t)| + \|W(t)\|_{H^2} = \mathcal{O}(e^{-2\nu t})$ as $t \to +\infty$.

In order to use this exponential decay of the auxiliary variable $W$, we rewrite (36) in the form

$$\partial_t V = L(\partial_x) V + DN(U_{\mathrm{per}}^-) V + N_1(U_{\mathrm{per}}^-, V) + \dot{q}\partial_1 U_{\mathrm{mf}} + G(U_{\mathrm{mf}} - U_{\mathrm{per}}^-, V) \;, \quad (40)$$

where

$$G(U_{\mathrm{mf}} - U_{\mathrm{per}}^-, V) = N_1(U_{\mathrm{mf}}, V) - N_1(U_{\mathrm{per}}^-, V) + (DN(U_{\mathrm{mf}}) - DN(U_{\mathrm{per}}^-))V \;.$$

To bound the last term in (40), we have to control the quantity

$$\begin{aligned}
\Omega(x,t) &= |U_{\mathrm{mf}}(x - ct - q(t), x) - U_{\mathrm{per}}^-(x)|\,|V(x,t)| \\
&= |U_{\mathrm{mf}}(x - ct - q(t), x) - U_{\mathrm{per}}^-(x)| e^{-\beta(x-ct)}|W(x - ct, t)| \;.
\end{aligned} \quad (41)$$

Take $\eta > 0$ small enough so that $\eta \beta \leq \nu$. By Theorem 2.6, we have

$$\Omega(x,t) \leq \begin{cases} C e^{-\varepsilon \beta_2 \eta t} |V(x,t)| & \text{if } x - ct - q(t) \leq -\eta t \;, \\ C e^{\beta(\eta t - q(t))} |W(x,t)| & \text{if } x - ct - q(t) > -\eta t \;. \end{cases}$$

It follows that $G(U_{\mathrm{mf}} - U_{\mathrm{per}}^-, V) = \mathcal{O}(e^{-\varepsilon\beta_3 t}V) + \mathcal{O}(e^{\nu t}W)$, where $\beta_3 = \beta_2\eta$. Thus, since $W = \mathcal{O}(e^{-2\nu t})$, we conclude that $G(U_{\mathrm{mf}} - U_{\mathrm{per}}^-, V)$ decays exponentially as $t \to \infty$. As a consequence of these observations, (36) has the form

$$\partial_t V = LV + DN(U_{\mathrm{per}}^-)V + N_1(U_{\mathrm{per}}^-, V) + \mathcal{O}(\dot{q}) + \mathcal{O}(e^{-\varepsilon\beta_3 t}V) + \mathcal{O}(e^{\nu t}W) \,. \qquad (42)$$

This is a small perturbation of the equation

$$\partial_t V = LV + DN(U_{\mathrm{per}}^-)V + N_1(U_{\mathrm{per}}^-, V) \,, \qquad (43)$$

which governs the evolution of the perturbations $V(x,t) = U(x,t) - U_{\mathrm{per}}^-(x)$ of the periodic equilibrium $U_{\mathrm{per}}^-$. As was already mentioned in section 3.4, sufficiently small solutions of (43) in $H_2^2$ converge diffusively to zero as $t \to \infty$, and due to the analysis in [ES02] the exponentially decreasing terms $\mathcal{O}(\dot{q})$, $\mathcal{O}(e^{-\varepsilon\beta_3 t}V)$, $\mathcal{O}(e^{\nu t}W)$ in (42) do not change the result. This is the idea of the proof of Theorem 2.8, which we now develop in more detail.

## 5.2 The unscaled equations

We first give explicit formulas for the various quantities appearing in equations (36), (38), (39), and (40). In what follows, $U$ will either denote $U_{\mathrm{mf}}$ or $U_{\mathrm{per}}^-$. If $V = (V_1, V_2)$ and $U = (U_1, U_2)$, equation (36) holds with

$$L(\partial_x)V = \begin{pmatrix} \partial_x^2 V_1 \\ -(1+\partial_x^2)^2 V_2 \end{pmatrix} \,,$$

$$DN(U)V = \begin{pmatrix} \frac{1}{2}(1 + 2c_0 U_1 - 3U_1^2) & 1 \\ 2\gamma U_1 U_2 & \alpha - \gamma + \gamma U_1^2 - 3U_2^2 \end{pmatrix} \begin{pmatrix} V_1 \\ V_2 \end{pmatrix} \,,$$

$$N_1(U,V) = \begin{pmatrix} \frac{1}{2}(c_0 - 3U_1)V_1^2 - \frac{1}{2}V_1^3 \\ U_2(\gamma V_1^2 - 3V_2^2) + 2\gamma U_1 V_1 V_2 - V_2^3 + \gamma V_1^2 V_2 \end{pmatrix} \,.$$

Next, if $W = (W_1, W_2)$ and $U_{\mathrm{mf}} = (u_{\mathrm{mf}}, v_{\mathrm{mf}})$, then (38) holds with (for instance)

$$N_2(U,V,W) = \begin{pmatrix} \frac{1}{2}(c_0 - 3U_1)V_1 W_1 - \frac{1}{2}V_1^2 W_1 \\ U_2(\gamma V_1 W_1 - 3V_2 W_2) + 2\gamma U_1 V_1 W_2 - V_2^2 W_2 + \gamma V_1 V_2 W_1 \end{pmatrix} \,.$$

The operator $\Lambda_\beta = L_{\beta,c_0} + DN(U_h)$ has the expression

$$\Lambda_\beta = \begin{pmatrix} (\partial_\xi - \beta)^2 + c_0(\partial_\xi - \beta) & 0 \\ 0 & -(1+(\partial_\xi - \beta)^2)^2 + c_0(\partial_\xi - \beta) \end{pmatrix}$$
$$+ \begin{pmatrix} \frac{1}{2}(1 + 2c_0 h - 3h^2) & 1 \\ 0 & \alpha - \gamma(1 - h^2) \end{pmatrix} \,.$$

Since $0 < \beta < \beta_0$, we know that $\lambda = 0$ is a simple isolated eigenvalue of $\Lambda_\beta$ with eigenfunction $e^{\beta\xi}\partial_\xi U_h$. The corresponding spectral projection is given by

$$\Pi_c W = \begin{pmatrix} h' e^{\beta\xi} \\ 0 \end{pmatrix} \int_{\mathbb{R}} (\psi_1^c W_1 + \psi_2^c W_2) \,\mathrm{d}\xi \,,$$

where
$$\psi_1^c = Kh'e^{(c_0-\beta)\xi}, \quad K = \left(\int_{\mathbb{R}} (h')^2 e^{c_0\xi}\,\mathrm{d}\xi\right)^{-1},$$

and $\psi_2^c$ is the unique solution in $H^2$ of the differential equation

$$-(1+(\partial_\xi+\beta)^2)^2\psi_2^c - c_0(\partial_\xi+\beta)\psi_2^c + (\alpha-\gamma(1-h^2))\psi_2^c + \psi_1^c = 0.$$

Finally, if $G = (G_1, G_2)$ and $U_{\mathrm{mf}} - U_{\mathrm{per}}^- = U^{\mathrm{d}} = (U_1^{\mathrm{d}}, U_2^{\mathrm{d}})$, then equation (40) holds with

$$\begin{aligned}
G_1(U^{\mathrm{d}}, V) &= c_0 U_1^{\mathrm{d}} V_1 - \tfrac{3}{2}(u_{\mathrm{mf}}+u_{\mathrm{per}}^-)U_1^{\mathrm{d}}V_1 - \tfrac{3}{2}U_1^{\mathrm{d}}V_1^2, \\
G_2(U^{\mathrm{d}}, V) &= 2\gamma(v_{\mathrm{mf}}U_1^{\mathrm{d}} + u_{\mathrm{per}}^- U_2^{\mathrm{d}})V_1 + \gamma(u_{\mathrm{mf}}+u_{\mathrm{per}}^-)U_1^{\mathrm{d}}V_2 - 3(v_{\mathrm{mf}}+v_{\mathrm{per}}^-)U_2^{\mathrm{d}}V_2 \\
&\quad -3U_2^{\mathrm{d}}V_2^2 + \gamma U_2^{\mathrm{d}}V_1^2 + 2\gamma U_1^{\mathrm{d}}V_1 V_2.
\end{aligned}$$

We now start the analysis of the evolution system (39), (40) for $V, W, q$. Our goal is to show that, if the initial data $V(0), W(0), q(0)$ are sufficiently small, then $W(t)$ and $\dot{q}(t)$ will go to zero exponentially, while $V(t)$ tends to zero diffusively as described in Theorem 2.4. As was explained in the previous section, we shall consider (40) as a perturbation of (43), which can be treated using the techniques developed in [Sch96], see [ES02]. In particular, we shall use the Bloch wave formalism introduced in section 3.2.

If $\hat{V} = \mathcal{T}V$ is the Bloch wave transform of $V$, we set $\hat{v}_{\mathrm{c}} = \hat{E}_{\mathrm{c}}\hat{V}$ and $\hat{v}_{\mathrm{s}} = \hat{E}_{\mathrm{s}}\hat{V}$, where $\hat{E}_{\mathrm{c}}, \hat{E}_{\mathrm{s}}$ are the mode filters defined by (21). In particular, $V = \mathcal{T}^{-1}(\hat{v}_{\mathrm{c}} + \hat{v}_{\mathrm{s}})$. For notational convenience, we also set $w \equiv W$. Then system (39), (40) is equivalent to

$$\begin{aligned}
\partial_t \hat{v}_{\mathrm{c}} &= \widehat{\mathcal{M}}\hat{v}_{\mathrm{c}} + \hat{E}_{\mathrm{c}}\widehat{\mathcal{N}}_1(\hat{v}_{\mathrm{c}},\hat{v}_{\mathrm{s}}) + \hat{E}_{\mathrm{c}}\widehat{\mathcal{H}}(\hat{v}_{\mathrm{c}},\hat{v}_{\mathrm{s}},w,q,t), \\
\partial_t \hat{v}_{\mathrm{s}} &= \widehat{\mathcal{M}}\hat{v}_{\mathrm{s}} + \hat{E}_{\mathrm{s}}\widehat{\mathcal{N}}_1(\hat{v}_{\mathrm{c}},\hat{v}_{\mathrm{s}}) + \hat{E}_{\mathrm{s}}\widehat{\mathcal{H}}(\hat{v}_{\mathrm{c}},\hat{v}_{\mathrm{s}},w,q,t), \\
\partial_t w &= \Lambda_\beta w + \Pi_s[(c-c_0)\partial_\xi w + \mathcal{N}_2(\hat{v}_{\mathrm{c}},\hat{v}_{\mathrm{s}},w,q,t)], \\
\dot{q} &= \mathcal{N}_3(\hat{v}_{\mathrm{c}},\hat{v}_{\mathrm{s}},w,q,t),
\end{aligned} \quad (44)$$

where $\mathcal{M} = L(\partial_x) + DN(U_{\mathrm{per}}^-)$, $\Lambda_\beta = L_{\beta,c_0} + DN(U_h)$, and $\widehat{\mathcal{M}} = \mathcal{T}\mathcal{M}\mathcal{T}^{-1}$. The remaining terms in (44) are given by

$$\begin{aligned}
\widehat{\mathcal{N}}_1(\hat{v}_{\mathrm{c}},\hat{v}_{\mathrm{s}}) &= \mathcal{T}N_1(U_{\mathrm{per}}^-, V), \\
\widehat{\mathcal{H}}(\hat{v}_{\mathrm{c}},\hat{v}_{\mathrm{s}},w,q,t) &= \mathcal{T}G(U_{\mathrm{mf}}-U_{\mathrm{per}}^-, V) + \dot{q}\mathcal{T}\partial_1 U_{\mathrm{mf}}, \\
\mathcal{N}_2(\hat{v}_{\mathrm{c}},\hat{v}_{\mathrm{s}},w,q,t) &= L_\Delta(U_{\mathrm{mf}})w + N_2(U_{\mathrm{mf}}, T_{-ct}V, w) + \dot{q}(\partial_1 U_{\mathrm{mf}})e^{\beta\xi}, \\
\mathcal{N}_3(\hat{v}_{\mathrm{c}},\hat{v}_{\mathrm{s}},w,q,t) &= -(\Pi_c(\partial_1 U_{\mathrm{mf}}e^{\beta\xi}))^{-1}\times \\
&\quad \Pi_c\Big[(c-c_0)\partial_\xi w + L_\Delta(U_{\mathrm{mf}})w + N_2(U_{\mathrm{mf}}, T_{-ct}V, w)\Big],
\end{aligned} \quad (45)$$

where $V = \mathcal{T}^{-1}(\hat{v}_{\mathrm{c}} + \hat{v}_{\mathrm{s}})$, $(T_{-ct}V)(\xi, t) = V(\xi+ct, t)$, and $\dot{q} = \mathcal{N}_3(\hat{v}_{\mathrm{c}},\hat{v}_{\mathrm{s}},w,q,t)$. For later use, we observe that $\mathcal{N}_2$ is a *linear* function of its third argument $w$, and so is $\mathcal{N}_3$.

As explained in [Sch96], it is useful to modify this system by eliminating the quadratic terms with respect to $\hat{v}_{\mathrm{c}}$ in the equation for $\hat{v}_{\mathrm{s}}$. We thus introduce the new variables $\hat{u}_{\mathrm{c}}, \hat{u}_{\mathrm{s}}$ defined by

$$\hat{u}_{\mathrm{c}} = \hat{v}_{\mathrm{c}}, \quad \hat{u}_{\mathrm{s}} = \hat{v}_{\mathrm{s}} - \tfrac{1}{2}\widehat{\mathcal{M}}^{-1}\hat{E}_{\mathrm{s}}\widehat{\mathcal{N}}_1''(0)[\hat{v}_{\mathrm{c}},\hat{v}_{\mathrm{c}}], \quad (46)$$

where $\widehat{\mathcal{N}}_1''(0) = D^2_{\hat{v}_c,\hat{v}_c}\widehat{\mathcal{N}}_1(\hat{v}_c,0)|_{\hat{v}_c=0}$. As we shall see in the next section, this change of variables eliminates the leading terms in the asymptotic behavior of $\hat{v}_s$ as $t \to \infty$ and simplifies therefore the analysis.

Applying this transformation to (44), we obtain our final system

$$\begin{aligned}
\partial_t \hat{u}_c &= \widehat{\mathcal{M}}\hat{u}_c + \widehat{\mathcal{N}}_c(\hat{u}_c, \hat{u}_s) + \widehat{\mathcal{H}}_c(\hat{u}_c, \hat{u}_s, w, q) \,, \\
\partial_t \hat{u}_s &= \widehat{\mathcal{M}}\hat{u}_s + \widehat{\mathcal{N}}_s(\hat{u}_c, \hat{u}_s) + \widehat{\mathcal{H}}_s(\hat{u}_c, \hat{u}_s, w, q) \,, \\
\partial_t w &= \Lambda_\beta w + \Pi_s((c-c_0)\partial_\xi w + \mathcal{N}_w(\hat{u}_c, \hat{u}_s, w, q)) \,, \\
\dot{q} &= \mathcal{N}_q(\hat{u}_c, \hat{u}_s, w, q) \,,
\end{aligned} \quad (47)$$

where

$$\begin{aligned}
\widehat{\mathcal{N}}_c(\hat{u}_c, \hat{u}_s) &= \hat{E}_c\widehat{\mathcal{N}}_1(\hat{v}_c, \hat{v}_s) \,, \\
\widehat{\mathcal{N}}_s(\hat{u}_c, \hat{u}_s) &= \hat{E}_s\widehat{\mathcal{N}}_1(\hat{v}_c, \hat{v}_s) - \tfrac{1}{2}\partial_t(\widehat{\mathcal{M}}^{-1}\hat{E}_s\widehat{\mathcal{N}}_1''(0)[\hat{v}_c, \hat{v}_c]) \,, \\
\widehat{\mathcal{H}}_c(\hat{u}_c, \hat{u}_s, w, q) &= \hat{E}_c\widehat{\mathcal{H}}(\hat{v}_c, \hat{v}_s, w, q, t) \,, \\
\widehat{\mathcal{H}}_s(\hat{u}_c, \hat{u}_s, w, q) &= \hat{E}_s\widehat{\mathcal{H}}(\hat{v}_c, \hat{v}_s, w, q, t) \,, \\
\mathcal{N}_w(\hat{u}_c, \hat{u}_s, w, q) &= \mathcal{N}_2(\hat{v}_c, \hat{v}_s, w, q, t) \,, \\
\mathcal{N}_q(\hat{u}_c, \hat{u}_s, w, q) &= \mathcal{N}_3(\hat{v}_c, \hat{v}_s, w, q, t) \,.
\end{aligned} \quad (48)$$

In the right-hand side of (48), the variables $\hat{v}_c, \hat{v}_s$ have to be replaced everywhere by their expression (46) in terms of $\hat{u}_c, \hat{u}_s$. For simplicity, we also dropped the explicit dependence on $t$ in $\widehat{\mathcal{H}}_c, \widehat{\mathcal{H}}_s, \mathcal{N}_w$, and $\mathcal{N}_q$.

**Remark 5.1** *If we set $w = q = 0$ and $U_{\mathrm{mf}} - U_{\mathrm{per}}^- = 0$, then (47) reduces to the system*

$$\partial_t \hat{u}_c = \widehat{\mathcal{M}}\hat{u}_c + \widehat{\mathcal{N}}_c(\hat{u}_c, \hat{u}_s) \,, \quad \partial_t \hat{u}_s = \widehat{\mathcal{M}}\hat{u}_s + \widehat{\mathcal{N}}_s(\hat{u}_c, \hat{u}_s) \,, \quad (49)$$

*which governs the evolution of the perturbations of the spatially periodic pattern $U_{\mathrm{per}}^-$. All the arguments in sections 5.3 and 5.4 apply a fortiori to (49), and show that these perturbations satisfy (9). Thus Theorem 2.4 follows as a particular case of Theorem 2.8.*

## 5.3 The renormalization procedure

In section 3.3, we described the asymptotic behavior of the solutions of the linearized equation (3) around the periodic steady state $U_{\mathrm{per}}^-$. We showed that the rescaled central part $\hat{E}_c\hat{V}(\ell t^{-1/2}, x, t)$ converges to the Gaussian function $Ae^{-d\ell^2}$ as $t \to \infty$, while the stable part is exponentially damped. This continuous rescaling of the Bloch variable $\ell$ can be used to treat the nonlinear problem also, see [EWW97]. However, it seems somewhat easier to use a *discrete* version of the above rescaling, which is known as the "renormalization group" method, see [BK92, Sch96]. The idea is to fix some (sufficiently small) $\sigma \in (0, 1)$ and to define, for all $n \in \mathbb{N}$, the rescaled quantities

$$\begin{aligned}
\hat{v}_{c,n}(\varkappa, x, \tau) &= \hat{u}_c(\sigma^n\varkappa, x, \sigma^{-2n}\tau) = (\widehat{\mathcal{L}}^n\hat{u}_c)(\varkappa, x, \sigma^{-2n}\tau) \,, \\
\hat{v}_{s,n}(\varkappa, x, \tau) &= \sigma^{-3n/2}\hat{u}_s(\sigma^n\varkappa, x, \sigma^{-2n}\tau) = \sigma^{-3n/2}(\widehat{\mathcal{L}}^n\hat{u}_s)(\varkappa, x, \sigma^{-2n}\tau) \,, \\
w_n(\xi, \tau) &= e^{\nu\sigma^{-2n}\tau}w(\xi, \sigma^{-2n}\tau) \,, \quad q_n(\tau) = q(\sigma^{-2n}\tau) \,,
\end{aligned} \quad (50)$$

where $\nu > 0$ is as in Proposition 2.2 and $\widehat{\mathcal{L}} \equiv \widehat{\mathcal{L}}_\sigma$ is the rescaling operator defined by

$$(\widehat{\mathcal{L}}\hat{f})(\ell, x) = \hat{f}(\sigma\ell, x).$$

Inserting (50) into (47), we obtain an evolution system for $\{\hat{v}_{c,n}, \hat{v}_{s,n}, w_n, q_n\}$ which we denote by $S_n$.

For notational convenience, we take our initial data for (47) at time $t = 1$ instead of $t = 0$. Starting from these data, we first solve the evolution system $S_1$ for $\{\hat{v}_{c,1}, \hat{v}_{s,1}, w_1, q_1\}$ on the time interval $\tau \in [\sigma^2, 1]$. Evaluating the result at time $\tau = 1$ provides the initial data for the rescaled system $S_2$ satisfied by $\{\hat{v}_{c,2}, \hat{v}_{s,2}, w_2, q_2\}$, which we again solve for $\tau \in [\sigma^2, 1]$. Repeating this procedure, we see that solving (47) for $t \in [1, \infty)$ is equivalent to solving a sequence of rescaled systems $S_n$ on the fixed time interval $[\sigma^2, 1]$. Since (47) is not autonomous, the rescaled systems $S_n$ will also depend explicitly on time.

The advantage of this iterative procedure is that the rescaled system $S_n$ becomes simpler as $n \to \infty$ because the asymptotically irrelevant terms in $S_n$ are multiplied by exponentially small factors such as $\sigma^n$. This is due to the prefactors in (50) which anticipate the decay of the quantities $\hat{u}_c, \hat{u}_s, w, q$. For instance, since the quadratic terms have been eliminated by (46), we expect that $\hat{u}_s(\ell, x, t)$ will decay like $1/t$ as $t \to \infty$. Now, since $t = \sigma^{-2n}\tau$ and $\tau$ varies in a bounded interval, we can think of $\sigma^n$ as being $1/\sqrt{t}$. Thus, the rescaled variable $\hat{v}_{s,n}(\tau) = \sigma^{-3n/2}(\widehat{\mathcal{L}}^n \hat{u}_s)(\sigma^{-2n}\tau)$ still converges to zero as $n \to \infty$. Moreover, replacing $\widehat{\mathcal{L}}^n \hat{u}_s$ with $\sigma^{3n/2}\hat{v}_{s,n}$ in the evolution system $S_n$ produces small factors $\sigma^{3n/2}$. Similarly, since $w(t) = \mathcal{O}(e^{-2\nu t})$ as $t \to \infty$, we expect that the rescaled quantity $w_n(\xi, \tau)$ defined in (50) will vanish rapidly as $n \to \infty$. In contrast, $\hat{u}_c$ converges to a Gaussian and $q$ to a finite limit $q_*$, so the corresponding rescaled quantities have no prefactors in (50).

Our starting point is the integral equation satisfied by the rescaled quantities (50) on the time interval $\tau \in [\sigma^2, 1]$. Using (47), we find

$$\begin{aligned}
\hat{v}_{c,n}(\varkappa, x, \tau) &= e^{\sigma^{-2n}\widehat{\mathcal{M}}_{c,n}(\tau-\sigma^2)}\hat{v}_{c,n-1}(\sigma\varkappa, x, 1) \\
&\quad + \sigma^{-2n}\int_{\sigma^2}^{\tau} e^{\sigma^{-2n}\widehat{\mathcal{M}}_{c,n}(\tau-\tau')}\widehat{\mathcal{N}}_{c,n}(\hat{v}_{c,n}, \hat{v}_{s,n})(\varkappa, x, \tau')\,d\tau' \\
&\quad + \sigma^{-2n}\int_{\sigma^2}^{\tau} e^{\sigma^{-2n}\widehat{\mathcal{M}}_{c,n}(\tau-\tau')}\widehat{\mathcal{H}}_{c,n}(\hat{v}_{c,n}, \hat{v}_{s,n}, w_n, q_n)(\varkappa, x, \tau')\,d\tau', \\
\hat{v}_{s,n}(\varkappa, x, \tau) &= e^{\sigma^{-2n}\widehat{\mathcal{M}}_{s,n}(\tau-\sigma^2)}\sigma^{-3/2}\hat{v}_{s,n-1}(\sigma\varkappa, x, 1) \\
&\quad + \sigma^{-7n/2}\int_{\sigma^2}^{\tau} e^{\sigma^{-2n}\widehat{\mathcal{M}}_{s,n}(\tau-\tau')}\widehat{\mathcal{N}}_{s,n}(\hat{v}_{c,n}, \hat{v}_{s,n})(\varkappa, x, \tau')\,d\tau' \quad (51) \\
&\quad + \sigma^{-7n/2}\int_{\sigma^2}^{\tau} e^{\sigma^{-2n}\widehat{\mathcal{M}}_{s,n}(\tau-\tau')}\widehat{\mathcal{H}}_{s,n}(\hat{v}_{c,n}, \hat{v}_{s,n}, w_n, q_n)(\varkappa, x, \tau')\,d\tau', \\
w_n(\xi, \tau) &= e^{\sigma^{-2n}(\Lambda_\beta+\nu)(\tau-\sigma^2)}w_{n-1}(\xi, 1) + \sigma^{-2n}\int_{\sigma^2}^{\tau} e^{\sigma^{-2n}(\Lambda_\beta+\nu)(\tau-\tau')}\times \\
&\quad \Pi_s\Big[(c-c_0)\partial_\xi w_n + \mathcal{N}_{w,n}(\hat{v}_{c,n}, \hat{v}_{c,n}, w_n, q_n)\Big](\xi, \tau')\,d\tau', \\
q_n(\tau) &= q_{n-1}(1) + \sigma^{-2n}\int_{\sigma^2}^{\tau} \mathcal{N}_{q,n}(\hat{v}_{c,n}, \hat{v}_{s,n}, w_n, q_n)\,d\tau',
\end{aligned}$$

where $\widehat{\mathcal{M}}_{c,n} = \widehat{\mathcal{L}}^n \hat{E}_c^h \widehat{\mathcal{M}} \widehat{\mathcal{L}}^{-n}$, $\widehat{\mathcal{M}}_{s,n} = \widehat{\mathcal{L}}^n \hat{E}_s^h \widehat{\mathcal{M}} \widehat{\mathcal{L}}^{-n}$, and

$$\begin{aligned}
\widehat{\mathcal{N}}_{c,n}(\hat{v}_{c,n}, \hat{v}_{s,n}) &= \widehat{\mathcal{L}}^n \widehat{\mathcal{N}}_c(\widehat{\mathcal{L}}^{-n}\hat{v}_{c,n}, \sigma^{3n/2}\widehat{\mathcal{L}}^{-n}\hat{v}_{s,n}) \,, \\
\widehat{\mathcal{N}}_{s,n}(\hat{v}_{c,n}, \hat{v}_{s,n}) &= \widehat{\mathcal{L}}^n \widehat{\mathcal{N}}_s(\widehat{\mathcal{L}}^{-n}\hat{v}_{c,n}, \sigma^{3n/2}\widehat{\mathcal{L}}^{-n}\hat{v}_{s,n}) \,, \\
\widehat{\mathcal{H}}_{c,n}(\hat{v}_{c,n}, \hat{v}_{s,n}, w_n, q_n) &= \widehat{\mathcal{L}}^n \widehat{\mathcal{H}}_c(\widehat{\mathcal{L}}^{-n}\hat{v}_{c,n}, \sigma^{3n/2}\widehat{\mathcal{L}}^{-n}\hat{v}_{s,n}, e^{-\nu\sigma^{-2n}\tau}w_n, q_n) \,, \\
\widehat{\mathcal{H}}_{s,n}(\hat{v}_{c,n}, \hat{v}_{s,n}, w_n, q_n) &= \widehat{\mathcal{L}}^n \widehat{\mathcal{H}}_s(\widehat{\mathcal{L}}^{-n}\hat{v}_{c,n}, \sigma^{3n/2}\widehat{\mathcal{L}}^{-n}\hat{v}_{s,n}, e^{-\nu\sigma^{-2n}\tau}w_n, q_n) \,, \\
\mathcal{N}_{w,n}(\hat{v}_{c,n}, \hat{v}_{s,n}, w_n, q_n) &= \mathcal{N}_w(\widehat{\mathcal{L}}^{-n}\hat{v}_{c,n}, \sigma^{3n/2}\widehat{\mathcal{L}}^{-n}\hat{v}_{s,n}, w_n, q_n) \,, \\
\mathcal{N}_{q,n}(\hat{v}_{c,n}, \hat{v}_{s,n}, w_n, q_n) &= e^{-\nu\sigma^{-2n}\tau}\mathcal{N}_q(\widehat{\mathcal{L}}^{-n}\hat{v}_{c,n}, \sigma^{3n/2}\widehat{\mathcal{L}}^{-n}\hat{v}_{s,n}, w_n, q_n) \,.
\end{aligned} \tag{52}$$

We recall that $\mathcal{N}_w$ and $\mathcal{N}_q$ are linear functions of the third argument $w$.

We shall control the evolution of $\hat{v}_{c,n}$, $\hat{v}_{s,n}$ in the $n$-dependent function spaces $\mathcal{K}^c_{\sigma^n}$, $\mathcal{K}^s_{\sigma^n}$, where

$$\mathcal{K}^c_\sigma = \mathcal{K}_{\sigma, 3/2} \,, \quad \mathcal{K}^s_\sigma = \mathcal{K}_{\sigma, 1} \,, \tag{53}$$

and $\mathcal{K}_{\sigma, \rho}$ is defined in (24). The reason for these particular choices of the parameter $\rho$ (which measures the decay in the Bloch variable $\ell$) will be explained in the proof of Lemma 5.4. On the other hand, we shall use the fixed spaces $H^2(\mathbb{R})$ and $\mathbb{R}$ to control the remaining variables $w_n$ and $q_n$, respectively.

To estimate the various terms in (51), we now list a number of lemmas which are very similar to the corresponding statements in [ES02, Section 7]. The proofs use exactly the same techniques, so we shall be rather brief and we refer the reader to [Sch96, ES02] for more details.

We first bound the linear semigroups generated by $\widehat{\mathcal{M}}_{c,n}$, $\widehat{\mathcal{M}}_{s,n}$, and $\Lambda_\beta + \nu$. The fact that applying $e^{\sigma^{-2n}\widehat{\mathcal{M}}_{c,n}\tau}$ improves the decay in the Bloch variable will be used to compensate the "loss of $\rho$" in estimate (55) below. Similarly, the smoothing property of $e^{\sigma^{-2n}(\Lambda_\beta+\nu)\tau}$ will allow to control the term $(c - c_0)\partial_\xi w_n$ in the equation for $w_n$.

**Lemma 5.2** *Fix $\rho_1 \geq \rho_2 \geq 0$. There exist positive constants $C_1$, $\nu_1$, and $\nu_2$ such that, for all $\sigma \in (0, 1]$, all $n \in \mathbb{N}$, and all $\tau, \tau' \in [\sigma^2, 1]$ with $\tau > \tau'$, one has*

$$\begin{aligned}
\|e^{\sigma^{-2n}\widehat{\mathcal{M}}_{c,n}(\tau-\tau')}\widehat{\mathcal{L}}^n \hat{E}_c^h \widehat{\mathcal{L}}^{-n}\hat{g}\|_{\mathcal{K}_{\sigma^n, \rho_1}} &\leq C_1(\tau - \tau')^{\rho_2 - \rho_1} \|\hat{g}\|_{\mathcal{K}_{\sigma^n, \rho_2}} \,, \\
\|e^{\sigma^{-2n}\widehat{\mathcal{M}}_{s,n}(\tau-\tau')}\widehat{\mathcal{L}}^n \hat{E}_s^h \widehat{\mathcal{L}}^{-n}\hat{g}\|_{\mathcal{K}_{\sigma^n, \rho_1}} &\leq C_1 e^{-\nu_1 \sigma^{-2n}(\tau-\tau')}(\tau - \tau')^{\rho_2 - \rho_1} \|\hat{g}\|_{\mathcal{K}_{\sigma^n, \rho_2}} \,, \\
\|e^{\sigma^{-2n}(\Lambda_\beta+\nu)(\tau-\tau')}\Pi_s w\|_{H^{\rho_1}} &\leq C_1 e^{-\nu_2 \sigma^{-2n}(\tau-\tau')}(\tau - \tau')^{(\rho_2 - \rho_1)/2} \|w\|_{H^{\rho_2}} \,.
\end{aligned}$$

**Proof.** These estimates follow directly from Lemma 3.1 and Proposition 2.2. We recall that $\nu_1 = \mathcal{O}(\varepsilon^2)$ and $\nu_2 = \mathcal{O}(\beta)$. □

Next, we estimate the nonlinear terms (52). Most of these terms are bounded in a straightforward way by "counting the powers of $\sigma$" and using estimate (26) as well as the identity

$$\widehat{\mathcal{L}}^n(\widehat{\mathcal{L}}^{-n}\hat{u} \star \widehat{\mathcal{L}}^{-n}\hat{v}) = \sigma^n(\hat{u} \star \hat{v}). \tag{54}$$

However, to bound the critical term $\widehat{\mathcal{N}}_{c,n}$, one has to use the structure of the system in a deeper way and to exploit some non-trivial cancellations. We need the following lemma which generalizes the fact that derivatives produce powers of $\sigma$ under scaling. As we shall see, although the nonlinearity $\mathcal{N}_{c,n}$ does not contain any derivative, it has a "derivative-like" structure, see also [Sch96].

**Lemma 5.3** *Let $\rho_1 \geq \rho_2 \geq 0$ and assume that $f \in C^2_{\text{per}}([-1/2, 1/2], C^2([0, 2\pi], \mathbb{C}))$ satisfies*

$$\|\partial^j_\ell f(\ell, \cdot)\|_{C^2([0,2\pi],\mathbb{C})} \leq C|\ell|^{2(\rho_1-\rho_2)-j}, \quad j = 0, 1, 2.$$

*Then there exists $C > 0$ such that, for all $\sigma \in (0, 1]$,*

$$\|(\widehat{\mathcal{L}}_\sigma f)\hat{u}\|_{\mathcal{K}_{\sigma,\rho_2}} \leq C\sigma^{2(\rho_1-\rho_2)}\|\hat{u}\|_{\mathcal{K}_{\sigma,\rho_1}}.$$

**Proof.** This follows from $\sup_{\ell \in [-1/2, 1/2]} |(\sigma\ell)^{2(\rho_1-\rho_2)}/(1+\ell^2)^{\rho_1-\rho_2}| \leq \sigma^{2(\rho_1-\rho_2)}$. □

**Lemma 5.4** *There exist positive constants $C_2, \bar{q}$, and $\nu_3$ such that for all $n \geq 2$, all $\sigma \in (0, 1]$ and all $\tau \in [\sigma^2, 1]$ the following estimates hold. If $\max\{\|\hat{v}_{c,n}\|_{\mathcal{K}^c_{\sigma^n}}, \|\hat{v}_{s,n}\|_{\mathcal{K}^s_{\sigma^n}}, \|w_n\|_{H^2}\} \leq 1$ and $|q_n| \leq \bar{q}$, then*

$$\|\widehat{\mathcal{N}}_{c,n}\|_{\mathcal{K}_{\sigma^n, 3/4}} \leq C_2 \sigma^{5n/2}(\|\hat{v}_{c,n}\|_{\mathcal{K}^c_{\sigma^n}} + \|\hat{v}_{s,n}\|_{\mathcal{K}^s_{\sigma^n}})^2, \tag{55}$$

$$\|\widehat{\mathcal{N}}_{s,n}\|_{\mathcal{K}^s_{\sigma^n}} \leq C_2 \sigma^{2n}(\|\hat{v}_{c,n}\|_{\mathcal{K}^c_{\sigma^n}} + \|\hat{v}_{s,n}\|_{\mathcal{K}^s_{\sigma^n}})^2, \tag{56}$$

$$\|\widehat{\mathcal{H}}_{c,n}\|_{\mathcal{K}^c_{\sigma^n}} \leq C_2 e^{-\nu_3 \sigma^{-2n}\tau}(\|\hat{v}_{c,n}\|_{\mathcal{K}^c_{\sigma^n}} + \|\hat{v}_{s,n}\|_{\mathcal{K}^s_{\sigma^n}} + \|w_n\|_{H^2}), \tag{57}$$

$$\|\widehat{\mathcal{H}}_{s,n}\|_{\mathcal{K}^s_{\sigma^n}} \leq C_2 e^{-\nu_3 \sigma^{-2n}\tau}(\|\hat{v}_{c,n}\|_{\mathcal{K}^c_{\sigma^n}} + \|\hat{v}_{s,n}\|_{\mathcal{K}^s_{\sigma^n}} + \|w_n\|_{H^2}), \tag{58}$$

$$\|\mathcal{N}_{w,n}\|_{H^2} \leq C_2(\varepsilon + |q_n| + \sigma^{n/2}(\|\hat{v}_{c,n}\|_{\mathcal{K}^c_{\sigma^n}} + \|\hat{v}_{s,n}\|_{\mathcal{K}^s_{\sigma^n}}))\|w_n\|_{H^2}, \tag{59}$$

$$|\mathcal{N}_{q,n}| \leq C_2 e^{-\nu\sigma^{-2n}\tau}\|w_n\|_{H^2}. \tag{60}$$

**Proof.** We start with (60). From (52), (50), (48) and (46) we have

$$\mathcal{N}_{q,n}(\hat{v}_{c,n}, \hat{v}_{s,n}, w_n, q_n) = e^{-\nu\sigma^{-2n}\tau}\mathcal{N}_3(\hat{v}_c, \hat{v}_s, w_n, q_n, t),$$

where $t = \sigma^{-2n}\tau$, $\hat{v}_c = \widehat{\mathcal{L}}^{-n}\hat{v}_{c,n}$ and $\hat{v}_s = \sigma^{3n/2}\widehat{\mathcal{L}}^{-n}\hat{v}_{s,n} + \frac{1}{2}\widehat{\mathcal{M}}^{-1}\hat{E}_s\widehat{\mathcal{N}}''_1(0)[\widehat{\mathcal{L}}^{-n}\hat{v}_{c,n}, \widehat{\mathcal{L}}^{-n}\hat{v}_{c,n}]$. In particular, if $V = \mathcal{T}^{-1}(\hat{v}_c + \hat{v}_s)$, it follows from (15), (25) that

$$\|V\|_{H^2} \leq C\sigma^{n/2}(\|\hat{v}_{c,n}\|_{\mathcal{K}^c_{\sigma^n}} + \sigma^{3n/2}\|\hat{v}_{s,n}\|_{\mathcal{K}^s_{\sigma^n}}).$$

Therefore, using (45), we find

$$|\mathcal{N}_{q,n}(\hat{v}_{c,n}, \hat{v}_{s,n}, w_n, q_n)| = e^{-\nu\sigma^{-2n}\tau}\left|\Pi_c(\partial_1 U_{\text{mf}}(\xi-q_n, \xi+ct)e^{\beta\xi})\right|^{-1}$$
$$\times \left|\Pi_c\left[(c-c_0)\partial_\xi w_n + L_\Delta(U_{\text{mf}})w_n + N_2(U_{\text{mf}}, V, w_n)\right]\right| \tag{61}$$
$$\leq Ce^{-\nu\sigma^{-2n}\tau}\left(\varepsilon + |q_n| + \sigma^{n/2}(\|\hat{v}_{c,n}\|_{\mathcal{K}^c_{\sigma^n}} + \|\hat{v}_{s,n}\|_{\mathcal{K}^s_{\sigma^n}})\right)\|w_n\|_{H^2},$$

since $\Pi_c(\partial_1 U_{\text{mf}}e^{\beta\xi}) = 1 + \mathcal{O}(\varepsilon + |q_n|)$, $|c - c_0| = \mathcal{O}(\varepsilon^2)$, $|L_\Delta(U_{\text{mf}})| = \mathcal{O}(\varepsilon + |q_n|)$, and since $N_2$ is a linear function of its third argument. This proves (60). Here and in the sequel, we also implicitly use the assumption $\max\{\|\hat{v}_{c,n}\|_{\mathcal{K}^c_{\sigma^n}}, \|\hat{v}_{s,n}\|_{\mathcal{K}^s_{\sigma^n}}, \|w_n\|_{H^2}\} \leq 1$. Following the same lines and using in addition estimate (61) for $\dot{q} = \sigma^{2n}\dot{q}_n = \mathcal{N}_{q,n}$, we obtain (59). Finally, from (53), (52), (48) and (26), we have

$$\|\widehat{\mathcal{H}}_{c,n}(\hat{v}_{c,n}, \hat{v}_{s,n}, w_n, q_n)\|_{\mathcal{K}^c_{\sigma^n}} \leq C\sigma^{-7n/2}\|\widehat{\mathcal{H}}(\hat{v}_c, \hat{v}_s, e^{-\nu\sigma^{-2n}\tau}w_n, q_n, t)\|_{\hat{H}^2_2},$$

where $t, \hat{v}_\text{c}, \hat{v}_\text{s}$ are as above, and $\widehat{\mathcal{H}}$ is given by (45). To bound the new term $G(U_\text{mf}-U_\text{per}^-, V)$, where $V = \mathcal{T}^{-1}(\hat{v}_\text{c} + \hat{v}_\text{s})$, we proceed as is suggested after (41) and obtain

$$\|G(U_\text{mf}-U_\text{per}^-, V)\|_{\mathbb{H}_2^2} \leq Ce^{-\nu_3 \sigma^{-2n}\tau}(\|\hat{v}_{\text{c},n}\|_{\mathcal{K}_{\sigma^n}^c} + \|\hat{v}_{\text{s},n}\|_{\mathcal{K}_{\sigma^n}^s} + \|w_n\|_{H^2}),$$

for some $\nu_3 \leq \min(\varepsilon\beta_2, \nu/2)$. Together with (61), this yields (57), and (58) is proved in a similar way.

To estimate $\widehat{\mathcal{N}}_{\text{c},n}$, we first note that $\widehat{\mathcal{N}}_1(\hat{v}_\text{c}, \hat{v}_\text{s}) = \mathcal{T}N_1(U_\text{per}^-, \mathcal{T}^{-1}(\hat{v}_\text{c} + \hat{v}_\text{s}))$ where

$$N_1(U_\text{per}^-, V) = \begin{pmatrix} \frac{1}{2}(c_0 - 3U_{\text{per},1}^-)V_1^2 \\ U_{\text{per},2}^-(\gamma V_1^2 - 3V_2^2) + 2\gamma U_{\text{per},1}^- V_1 V_2 \end{pmatrix} + \mathcal{O}(V^3).$$

Therefore,

$$\begin{aligned}
\widehat{\mathcal{N}}_{\text{c},n} &= \widehat{\mathcal{L}}^n \hat{E}_\text{c} \widehat{\mathcal{N}}_1(\widehat{\mathcal{L}}^{-n}\hat{v}_{\text{c},n}, \sigma^{3n/2}\widehat{\mathcal{L}}^{-n}\hat{v}_{\text{s},n} + \tfrac{1}{2}\widehat{\mathcal{M}}^{-1}\hat{E}_\text{s}\widehat{\mathcal{N}}_1''(0)[\widehat{\mathcal{L}}^{-n}\hat{v}_{\text{c},n}, \widehat{\mathcal{L}}^{-n}\hat{v}_{\text{c},n}]) \\
&= \widehat{\mathcal{L}}^n \hat{E}_\text{c} \begin{pmatrix} \frac{1}{2}(c_0 - 3U_{\text{per},1}^-)a_1(\hat{v}_{\text{c},n}, \hat{v}_{\text{s},n}) \\ U_{\text{per},2}^-(\gamma a_1(\hat{v}_{\text{c},n}, \hat{v}_{\text{s},n}) - 3a_2(\hat{v}_{\text{c},n}, \hat{v}_{\text{s},n})) + 2\gamma U_{\text{per},1}^- a_3(\hat{v}_{\text{c},n}, \hat{v}_{\text{s},n}) \end{pmatrix} \\
&\quad + \widehat{\mathcal{N}}_{\text{c},n}^{(3)}(\hat{v}_{\text{c},n}, \hat{v}_{\text{s},n}),
\end{aligned}$$

where $a_1, a_2, a_3$ are the quadratic terms in $\hat{v}_{\text{c},n}, \hat{v}_{\text{s},n}$ and $\widehat{\mathcal{N}}_{\text{c},n}^{(3)}$ contains the higher order convolutions. Setting $\hat{v}_{\text{c},n} = (\hat{v}_{\text{c},n}^{(1)}, \hat{v}_{\text{c},n}^{(2)})$, we find

$$a_1(\hat{v}_{\text{c},n}, \hat{v}_{\text{s},n}) = (\widehat{\mathcal{L}}^{-n}\hat{v}_{\text{c},n}^{(1)})^{\star 2} + 2\sigma^{3n/2}(\widehat{\mathcal{L}}^{-n}\hat{v}_{\text{c},n}^{(1)} \star \widehat{\mathcal{L}}^{-n}\hat{v}_{\text{s},n}^{(1)}) + \sigma^{3n}(\widehat{\mathcal{L}}^{-n}\hat{v}_{\text{s},n}^{(1)})^{\star 2} =: a_{11} + a_{12} + a_{13},$$

and similar expressions hold for $a_2$ and $a_3$. Most of these terms can be bounded in a straightforward way, using (54) together with the fact that $\|\widehat{\mathcal{L}}^n \hat{E}_\text{c} \widehat{\mathcal{L}}^{-n}f\|_{\mathcal{K}_{\sigma^n}^c} \leq C\|f\|_{\mathcal{K}_{\sigma^n}^c}$ for some $C > 0$ independent of $\sigma$ and $n$. For instance, $\|\widehat{\mathcal{L}}^n a_{12}\|_{\mathcal{K}_{\sigma^n}^c} \leq C\sigma^{5n/2}\|\hat{v}_{\text{c},n}\|_{\mathcal{K}_{\sigma^n}^c}\|\hat{v}_{\text{s},n}\|_{\mathcal{K}_{\sigma^n}^s}$. However, (54) is not sufficient for the quadratic (and cubic) convolutions involving $\hat{v}_{\text{c},n}$ only: proceeding as above, one would obtain $\|\widehat{\mathcal{L}}^n a_{11}\|_{\mathcal{K}_{\sigma^n}^c} \leq C\sigma^n \|\hat{v}_{\text{c},n}\|_{\mathcal{K}_{\sigma^n}^c}^2$, while we need $\sigma^{5n/2}$ in the right-hand side, see (55). Thus, more careful bounds are necessary.

Following [Sch96] we write $\widehat{\mathcal{N}}_{\text{c},n} = \hat{s}_1 + \hat{s}_2 + \widehat{\mathcal{N}}_{\text{c},n,r}$, where $\hat{s}_1$ and $\hat{s}_2$ contain the quadratic and cubic convolutions of $\hat{v}_{\text{c},n}$, respectively, and $\widehat{\mathcal{N}}_{\text{c},n,r}$ contains the higher convolutions of $\hat{v}_{\text{c},n}$ and all the convolutions with at least one factor $\sigma^{3n/2}\hat{v}_{\text{s},n}$. Thus

$$\|\widehat{\mathcal{N}}_{\text{c},n,r}\|_{\mathcal{K}_{\sigma^n}^c} \leq C\sigma^{5n/2}(\|\hat{v}_{\text{c},n}\|_{\mathcal{K}_{\sigma^n}^c} + \|\hat{v}_{\text{s},n}\|_{\mathcal{K}_{\sigma^n}^s})^2.$$

On the other hand, $\hat{s}_1$ is given by $\hat{s}_1 = \widehat{\mathcal{L}}^n \hat{E}_\text{c} \mathcal{T}B(\mathcal{T}^{-1}\widehat{\mathcal{L}}^{-n}\hat{v}_{\text{c},n}, \mathcal{T}^{-1}\widehat{\mathcal{L}}^{-n}\hat{v}_{\text{c},n})$, where $B$ is the symmetric bilinear form defined by $B(U_1, U_2) = \frac{1}{2}(D_{VV}^2 N_1(U_\text{per}^-, V)|_{V=0})[U_1, U_2]$. Thus

$$\mathcal{T}B(U, V)(\ell, x) = \int_{-1/2}^{1/2} B_2(x)[\mathcal{T}U(\ell - m, x), \mathcal{T}V(m, x)]\,\mathrm{d}m,$$

where $B_2$ is again a symmetric bilinear form, depending explicitly on $x$ due to the multiplication by $U_\text{per}^-$. Writing $\hat{v}_{\text{c},n}(\varkappa, x) = \widehat{\mathcal{L}}^n \hat{E}_\text{c}^\text{h} \widehat{\mathcal{L}}^{-n}\hat{v}_{\text{c},n}(\varkappa, x) = a_n(\varkappa)\varphi(\sigma^n\varkappa, x)$, where $\varphi(\ell, \cdot)$ is the critical eigenfunction of $\mathcal{M}(\ell)$ from Lemma 3.1, and using (20), (21), we obtain $\hat{s}_1(\varkappa, x) = \varphi(\sigma^n\varkappa, x)\chi(\sigma^n\varkappa)\bar{s}_1(\varkappa)$ with

$$\begin{aligned}
\bar{s}_1(\varkappa) &= \sigma^n \int_0^{2\pi} \left\langle \psi(\sigma^n\varkappa, x), \int_{-1/2\sigma^n}^{1/2\sigma^n} B_2(x)[\hat{v}_{\text{c},n}(\varkappa-m, x)\hat{v}_{\text{c},n}(m, x)]\,\mathrm{d}m \right\rangle_{\mathbb{C}^2} \mathrm{d}x \\
&= \sigma^n \int_{-1/2\sigma^n}^{1/2\sigma^n} (\widehat{\mathcal{L}}^n K_1(\varkappa, \varkappa - m, m))a_n(\varkappa - m)a_n(m)\,\mathrm{d}m,
\end{aligned}$$

and where $\psi(\ell, \cdot)$ is the critical eigenfunction of $\mathcal{M}(\ell)^*$, and

$$K_1(\ell, \ell - m, m) = \int_0^{2\pi} \left\langle \psi(\ell, x), B_2(x)[\varphi(\ell - m, x), \varphi(m, x)] \right\rangle_{\mathbb{C}^2} \, \mathrm{d}x \,. \tag{62}$$

Similarly, $\hat{s}_2(\varkappa, x) = \varphi(\sigma^n \varkappa, x) \chi(\sigma^n \varkappa) \bar{s}_2(\varkappa)$, where

$$\bar{s}_2(\varkappa) = \sigma^{2n} \iint (\widehat{\mathcal{L}}^n K_2(\varkappa, \varkappa - m, m - k, k)) a_n(\varkappa - m) a_n(m - k) a_n(k) \, \mathrm{d}k \, \mathrm{d}m \,.$$

The kernel $K_2(\ell, \ell - m, m - k, k)$ encodes the cubic convolutions of $\hat{v}_{c,n}$.

Next we show that $|K_1(\ell, \ell - m, m)| = \mathcal{O}(\ell^2 + m^2)$ and $|K_2(\ell, \ell - m, m - k, k)| = \mathcal{O}(\ell + m + k)$ as $\ell, m, k \to 0$. Then, using Lemma 5.3 with $\rho_1 = 3/2$ and $\rho_2 = 3/4$, we gain a factor $\sigma^{3n/2}$ when estimating $\hat{s}_1$. Similarly, choosing for instance $\rho_1 = 1$ and $\rho_2 = 3/4$, we gain $\sigma^{n/2}$ when estimating $\hat{s}_2$. This gives the correct power of $\sigma$ in (55), and explains why we estimate $\widehat{\mathcal{N}}_{c,n}$ in $\mathcal{K}_{\sigma^n,3/4}$ and not in $\mathcal{K}_{\sigma^n}^c$. As was already mentioned, this loss in $\rho$ will be compensated by the smoothing properties of $e^{\widehat{\mathcal{M}}_{c,n}\tau}$, see Lemma 5.2.

First, from the translation invariance of (3), we have $K_1(0, 0, 0) = 0$ and $K_2(0, 0, 0, 0) = 0$. Indeed, the case $\ell = m = k = 0$ corresponds to the spatially periodic case, in which there exists a center manifold $\Gamma = \{U_{\mathrm{per}}^-(x - a) \mid a \in \mathbb{R}\}$ of $2\pi$-periodic equilibria, cf. section 4. It is not difficult to verify (see [ES02]) that the flow induced by (3) on $\Gamma$ is given by

$$\dot{a} = K_1(0, 0, 0) a^2 + K_2(0, 0, 0, 0) a^3 + \mathcal{O}(a^4) \,,$$

which immediately yields $K_1(0, 0, 0) = 0$ and $K_2(0, 0, 0, 0) = 0$.

It remains to show $|K_1(\ell, \ell - m, m)| = \mathcal{O}(\ell^2 + m^2)$. From $B(U, V) = B(V, U)$ we have $K_1(\ell, \ell - m, m) = K_1(\ell, m, \ell - m)$ and hence $\frac{\mathrm{d}}{\mathrm{d}m} K_1(0, 0, 0) = 0$. Finally, $\frac{\mathrm{d}}{\mathrm{d}\ell} K_1(0, 0, 0) = 0$ is a consequence of $\mu_1(\ell) = \mathcal{O}(\ell^2)$, which can be seen as follows. As in [Sch98b] we write

$$\varphi(\ell) = \varphi_0 + \ell \varphi_1 + \mathcal{O}(\ell^2), \quad \psi(\ell) = \psi_0 + \ell \psi_1 + \mathcal{O}(\ell^2),$$
$$\mathcal{M}(\ell) = \mathcal{M}_0 + \ell \mathcal{M}_1 + \mathcal{O}(\ell^2), \quad \mathcal{M}(\ell)^* = \mathcal{M}_0^* + \ell \mathcal{M}_1^* + \mathcal{O}(\ell^2).$$

From $\mathcal{M}(\ell) \varphi(\ell) = \mathcal{M}_0 \varphi_0 + \ell(\mathcal{M}_1 \varphi_0 + \mathcal{M}_0 \varphi_1) + \mathcal{O}(\ell^2) = \mathcal{O}(\ell^2)$ and the similar relation for $\mathcal{M}(\ell)^* \psi(\ell)$, we obtain

$$\mathcal{M}_0 \varphi_0 = 0, \quad \mathcal{M}_0^* \psi_0 = 0, \quad \mathcal{M}_0 \varphi_1 + \mathcal{M}_1 \varphi_0 = 0, \quad \mathcal{M}_0^* \psi_1 + \mathcal{M}_1^* \psi_0 = 0. \tag{63}$$

Moreover, for any fixed $\ell$ and any $2\pi$-periodic function $\hat{V}$, we have

$$\begin{aligned}
\partial_x(\mathcal{M}(\ell) \hat{V}) &= \partial_x \left( e^{-\mathrm{i}\ell x} (L + DN(U_{\mathrm{per}}^-)) \hat{V} e^{\mathrm{i}\ell x} \right) \\
&= e^{-\mathrm{i}\ell x} \left( (L + DN(U_{\mathrm{per}}^-)) \partial_x \hat{V} + D^2 N(U_{\mathrm{per}}^-)(\hat{V}, \partial_x U_{\mathrm{per}}^-) \right) e^{\mathrm{i}\ell x} \\
&= \mathcal{M}(\ell) \partial_x \hat{V} + 2 B_2(x)[\hat{V}, \partial_x U_{\mathrm{per}}^-] \,.
\end{aligned} \tag{64}$$

If we now choose $\hat{V} = \varphi_0 + \ell \varphi_1$, we have $\mathcal{M}(\ell) \hat{V} = \mathcal{O}(\ell^2)$. Hence it follows from (64) that

$$\begin{aligned}
\mathcal{M}_0 \partial_x \varphi_0 + 2 B_2[\varphi_0, \partial_x U_{\mathrm{per}}^-] &= 0 \,, \\
\mathcal{M}_1 \partial_x \varphi_0 + \mathcal{M}_0 \partial_x \varphi_1 + 2 B_2[\varphi_1, \partial_x U_{\mathrm{per}}^-] &= 0 \,.
\end{aligned} \tag{65}$$

Using (63), (65), and $\varphi_0 = c_N \partial_x U^-_{\text{per}}$, where $c_N$ is the normalization factor from Lemma 3.1, we thus obtain

$$\frac{\mathrm{d}}{\mathrm{d}\ell} K_1(0,0,0) = \langle \psi_1, B_2[\varphi_0, \varphi_0] \rangle + \langle \psi_0, B_2[\varphi_0, \varphi_1] \rangle$$
$$= -\tfrac{1}{2c_N} \langle \mathcal{M}_0^* \psi_1 + \mathcal{M}_1^* \psi_0, \partial_x \varphi_0 \rangle - \tfrac{1}{2c_N} \langle \mathcal{M}_0^* \psi_0, \partial_x \varphi_1 \rangle = 0.$$

This concludes the proof of (55). Finally, estimate (56) on $\widehat{\mathcal{N}}_{s,n}$ follows from expressing $\partial_t \hat{u}_c$ in $-\tfrac{1}{2} \partial_t (\widehat{\mathcal{M}}^{-1} \hat{E}_s \widehat{\mathcal{N}}_1''(0)[\hat{v}_c, \hat{v}_c])$ by (47) and using (55), (54) and Lemma 5.3. This gives the desired power of $\sigma$, because $\widehat{\mathcal{N}}_{s,n}$ contains no quadratic terms in $\hat{v}_{c,n}$ since these terms were precisely eliminated by the change of variables (46). The proof of Lemma 5.4 is now complete. □

To estimate the integrals in (51), we introduce the quantities

$$R_{v,n} = \sup\nolimits_{\tau \in [\sigma^2, 1]} \|\hat{v}_{c,n}(\tau)\|_{\mathcal{K}^c_{\sigma n}} + \sup\nolimits_{\tau \in [\sigma^2, 1]} \|\hat{v}_{s,n}(\tau)\|_{\mathcal{K}^s_{\sigma n}},$$
$$R_{w,n} = \sup\nolimits_{\tau \in [\sigma^2, 1]} \|w_n(\tau)\|_{H^2}, \qquad R_{q,n} = \sup\nolimits_{\tau \in [\sigma^2, 1]} |q_n(\tau)|.$$

**Lemma 5.5** *There exist positive constants $C_3, C_4$ such that for all $\sigma \in (0, 1]$, all $\tau \in [\sigma^2, 1]$ and all $n \geq 2$ the following estimates hold. If $\max(R_{v,n}, R_{w,n}) \leq 1$ and $R_{q,n} \leq \bar{q}$, then*

$$\sigma^{-2n} \| \int_{\sigma^2}^\tau e^{\sigma^{-2n} \widehat{\mathcal{M}}_{c,n}(\tau - \tau')} \widehat{\mathcal{N}}_{c,n}(\tau') \, \mathrm{d}\tau' \|_{\mathcal{K}^c_{\sigma n}} \leq C_3 \sigma^{n/2} R_{v,n}^2,$$
$$\sigma^{-7n/2} \| \int_{\sigma^2}^\tau e^{\sigma^{-2n} \widehat{\mathcal{M}}_{s,n}(\tau - \tau')} \widehat{\mathcal{N}}_{s,n}(\tau') \, \mathrm{d}\tau' \|_{\mathcal{K}^s_{\sigma n}} \leq C_3 \sigma^{n/2} R_{v,n}^2,$$
$$\sigma^{-2n} \| \int_{\sigma^2}^\tau e^{\sigma^{-2n} \widehat{\mathcal{M}}_{c,n}(\tau - \tau')} \widehat{\mathcal{H}}_{c,n}(\tau') \, \mathrm{d}\tau' \|_{\mathcal{K}^c_{\sigma n}} \leq C_3 \sigma^{n/2} (R_{v,n} + R_{w,n}),$$
$$\sigma^{-7n/2} \| \int_{\sigma^2}^\tau e^{\sigma^{-2n} \widehat{\mathcal{M}}_{s,n}(\tau - \tau')} \widehat{\mathcal{H}}_{s,n}(\tau') \, \mathrm{d}\tau' \|_{\mathcal{K}^s_{\sigma n}} \leq C_3 \sigma^{n/2} (R_{v,n} + R_{w,n}), \quad (66)$$
$$\sigma^{-2n} \| \int_{\sigma^2}^\tau e^{\sigma^{-2n}(\Lambda_\beta + \nu)(\tau - \tau')} \Pi_s \big[ (c - c_0) \partial_\xi w_n + \mathcal{N}_{w,n} \big](\tau') \, \mathrm{d}\tau' \|_{H^2}$$
$$\leq (C_4(\varepsilon + \bar{q}) + C_3 \sigma^{n/2} R_{v,n}) R_{w,n},$$
$$\sigma^{-2n} | \int_{\sigma^2}^\tau \mathcal{N}_{q,n}(\tau') \, \mathrm{d}\tau' | \leq C_3 e^{-\nu \sigma^{-2n+2}} R_{w,n}.$$

**Remark 5.6** *The proof shows that $C_3 \to \infty$ as $\varepsilon \to 0$, while $C_4$ can be chosen independent of $\varepsilon$.*

**Proof.** These estimates are easily obtained by combining Lemmas 5.2 and 5.4. The assumption $n \geq 2$ is used to simplify the bounds on $\widehat{\mathcal{H}}_{c,n}$ and $\widehat{\mathcal{H}}_{s,n}$. For instance, in the estimate involving $\widehat{\mathcal{H}}_{c,n}$, we use

$$\sigma^{-2n} \int_{\sigma^2}^\tau C_1 C_2 e^{-\nu_3 \sigma^{-2n} \tau'} \, \mathrm{d}\tau' \leq \frac{C_1 C_2}{\nu_3} e^{-\nu_3 \sigma^{-2n+2}} \leq C_3 \sigma^{n/2}.$$

The last inequality is very crude if $n$ is large, but it is sufficient for our purposes because it matches what we have for $\widehat{\mathcal{N}}_{c,n}, \widehat{\mathcal{N}}_{s,n}$. Note however that we keep the exponential factor in the estimate for $\mathcal{N}_{q,n}$. □

Finally, the following bounds hold for the first terms on the right-hand side of (51).

**Lemma 5.7** *There exists $C_5, C_6 > 0$ and $m \in \mathbb{N}$ such that, for all $\sigma \in (0,1]$ and all $\tau \in [\sigma^2, 1]$, one has*

$$\|e^{\sigma^{-2n}\widehat{\mathcal{M}}_{c,n}(\tau-\sigma^2)}\widehat{\mathcal{L}}^n \hat{E}_c^h \widehat{\mathcal{L}}^{-n}\widehat{\mathcal{L}}\hat{g}\|_{\mathcal{K}_{\sigma^n}^c} \leq C_5 \sigma^{-m}\|\hat{g}\|_{\mathcal{K}_{\sigma^{n-1}}^c},$$
$$\|e^{\sigma^{-2n}\widehat{\mathcal{M}}_{s,n}(\tau-\sigma^2)}\widehat{\mathcal{L}}^n \hat{E}_s^h \widehat{\mathcal{L}}^{-n}\sigma^{-3/2}\widehat{\mathcal{L}}\hat{g}\|_{\mathcal{K}_{\sigma^n}^s} \leq C_5 \sigma^{-m}e^{-\nu_1\sigma^{-2n}(\tau-\sigma^2)}\|\hat{g}\|_{\mathcal{K}_{\sigma^{n-1}}^s},$$
$$\|e^{\sigma^{-2n}(\Lambda_\beta+\nu)(\tau-\sigma^2)}\Pi_s w\|_{H^2} \leq C_6 e^{-\nu_2\sigma^{-2n}(\tau-\sigma^2)}\|w\|_{H^2}.$$

**Proof.** This follows immediately from Lemma 5.2 and estimate (26). The constant $C_6$ can be chosen independent of $\varepsilon$. $\square$

In the sequel, we set $\mu = C_4(\varepsilon + \bar{q})$ and assume that $\mu$ is sufficiently small so that

$$(2C_6 + 1)\mu \leq 1/2. \tag{67}$$

Combining the above Lemmas, we are now ready to give a priori bounds on the solution of (51) in terms of the initial data. Let

$$\rho_{v,n} = \|\hat{v}_{c,n}(\cdot,\cdot,1)\|_{\mathcal{K}_{\sigma^n}^c} + \|\hat{v}_{s,n}(\cdot,\cdot,1)\|_{\mathcal{K}_{\sigma^n}^s}, \quad \rho_{w,n} = \|w_n(\cdot,1)\|_{H^2}, \quad \rho_{q,n} = |q_n(1)|.$$

**Proposition 5.8** *There exist $\sigma_0 > 0$ and $C_7 > 0$ such that, for all $\sigma \in (0, \sigma_0)$ and all $n \geq 2$, the following holds. If $\rho_{v,n-1} \leq C_7 \sigma^m$ with $m$ as in Lemma 5.7, and $\rho_{w,n-1} + \rho_{q,n-1} \leq C_7$, then (51) has a unique solution $(\hat{v}_{c,n}, \hat{v}_{s,n}, w_n, q_n) \in C([\sigma^2, 1], \mathcal{K}_{\sigma^n}^c \times \mathcal{K}_{\sigma^n}^s \times H^2 \times \mathbb{R})$. In addition, we have $\max(R_{v,n}, R_{w,n}) \leq 1$, $R_{q,n} \leq \bar{q}$, and*

$$\begin{aligned}
R_{v,n} &\leq C_5 \sigma^{-m}\rho_{v,n-1} + 2C_3 \sigma^{n/2}(R_{v,n}^2 + R_{v,n} + R_{w,n}), \\
R_{w,n} &\leq C_6 \rho_{w,n-1} + (\mu + C_3 \sigma^{n/2} R_{v,n})R_{w,n}, \\
R_{q,n} &\leq \rho_{q,n-1} + C_3\, e^{-\nu\sigma^{-2n+2}} R_{w,n}.
\end{aligned} \tag{68}$$

**Proof.** Let $U_n = (\hat{v}_{c,n}, \hat{v}_{s,n}, w_n, q_n)$ and let $X_n = C([\sigma^2, 1], \mathcal{K}_{\sigma^n}^c \times \mathcal{K}_{\sigma^n}^s \times H^2 \times \mathbb{R})$ equipped with the norm $\|U_n\|_{X_n} = R_{v,n} + R_{w,n} + R_{q,n}/\bar{q}$. Let also $B_n$ be the unit ball in $X_n$. Given initial data $\hat{v}_{c,n-1}(\cdot,\cdot,1), \hat{v}_{s,n-1}(\cdot,\cdot,1), w_{n-1}(\cdot,1), q_{n-1}(1)$, the right-hand side of (51) defines a map $F_n : X_n \to X_n$. From Lemmas 5.5 and 5.7, we know that, if $U_n \in B_n$, then

$$\|F(U_n)\|_{X_n} \leq C_5(\sigma^{-m}\rho_{v,n-1} + \rho_{q,n-1}/\bar{q}) + C_6 \rho_{w,n-1} + (\mu + C\sigma^{n/2})\|U_n\|_{X_n}.$$

Similarly, if $U_n, \tilde{U}_n \in B_n$, we find

$$\|F(U_n) - F(\tilde{U}_n)\|_{X_n} \leq (\mu + C\sigma^{n/2})\|U_n - \tilde{U}_n\|_{X_n}.$$

By (67), we have $\mu \leq 1/2$. If we now assume $C\sigma < \frac{1}{2}$ and $\sigma^{-m}\rho_{v,n-1} + \rho_{w,n-1} + \rho_{q,n-1}/\bar{q} \leq C_7$ for some sufficiently small $C_7 > 0$, we see that $F_n$ maps $B_n$ into itself and is a strict contraction there. The unique fixed point $U_n$ is the desired solution and satisfies (68). $\square$

## 5.4 Iteration and Conclusion

To show that the recursion relation (68) can be iterated and to conclude the proof of Theorem 2.8 we need a better control on the critical term $\hat{v}_{c,n}(\varkappa, x, \tau)$. For each $n \in \mathbb{N}$, we decompose the solution $\hat{v}_{c,n}(\cdot, \cdot, 1)$ into a Gaussian part and a remainder, i.e.

$$\hat{v}_{c,n}(\varkappa, x, 1) = A_n \Psi(\varkappa) \varphi(\sigma^n \varkappa, x) + \hat{r}_n(\varkappa, x) ,$$

where $A_n \in \mathbb{C}$ and $\Psi(\varkappa) = e^{-d\varkappa^2}$ with $d$ and $\varphi(\ell, x)$ from Lemma 3.1. The amplitude $A_n$ is determined by the condition $\hat{r}_n(0, x) = 0$, $x \in [0, 2\pi]$. Equivalently, $A_n = \widehat{\Pi} \hat{v}_{c,n}(\cdot, \cdot, 1)$, where $\widehat{\Pi} : \mathcal{K}_{\sigma,\rho} \to \mathbb{C}$ is defined by

$$\widehat{\Pi} f = \langle \psi(0, \cdot), f(0, \cdot) \rangle ,$$

see (20). Then (51) can be decomposed accordingly and takes the form

$$A_n = A_{n-1} + \sigma^{-2n} \widehat{\Pi} \left( \int_{\sigma^2}^1 e^{\sigma^{-2n} \widehat{\mathcal{M}}_{c,n}(1-\tau')} (\widehat{\mathcal{N}}_{c,n} + \widehat{\mathcal{H}}_{c,n})(\tau') \, \mathrm{d}\tau' \right) , \tag{69}$$

$$\begin{aligned}\hat{r}_n(\varkappa, x) &= e^{\sigma^{-2n} \widehat{\mathcal{M}}_{c,n}(1-\sigma^2)} \hat{r}_{n-1}(\sigma \varkappa, x) \\ &\quad + \sigma^{-2n} \int_{\sigma^2}^1 e^{\sigma^{-2n} \widehat{\mathcal{M}}_{c,n}(1-\tau')} (\widehat{\mathcal{N}}_{c,n} + \widehat{\mathcal{H}}_{c,n})(\varkappa, x, \tau') \, \mathrm{d}\tau' \\ &\quad + e^{\sigma^{-2n} \widehat{\mathcal{M}}_{c,n}(1-\sigma^2)} A_{n-1} \Psi(\sigma \varkappa) \varphi(\sigma^n \varkappa, x) - A_n \Psi(\varkappa) \varphi(\sigma^n \varkappa, x) . \end{aligned} \tag{70}$$

We also define $\rho_{r,n} = \|\hat{r}_n\|_{\mathcal{K}_{\sigma^n}^c} + \|\hat{v}_{s,n}(\cdot, \cdot, 1)\|_{\mathcal{K}_{\sigma^n}^s}$. By construction we have $\rho_{v,n} \leq C(|A_n| + \rho_{r,n})$. Our main estimate is now:

**Proposition 5.9** *Under the assumptions of Proposition 5.8, the solution of (51) satisfies*

$$\begin{aligned} |A_n - A_{n-1}| &\leq C_8 \sigma^{n/2} \sigma^{-m} (|A_{n-1}| + \rho_{r,n-1} + \rho_{w,n-1}) , \\ \rho_{r,n} &\leq C_8 \sigma \rho_{r,n-1} + C_8 \sigma^{n/2} \sigma^{-m} (|A_{n-1}| + \rho_{r,n-1} + \rho_{w,n-1}) , \\ \rho_{w,n} &\leq (C_8 \sigma + \tfrac{1}{2}) \rho_{w,n-1} , \\ |q_n - q_{n-1}| &\leq C_8 \, e^{-\nu \sigma^{-2n+2}} \rho_{w,n-1} , \end{aligned} \tag{71}$$

*for some $C_8 > 0$.*

**Proof.** Since $|\widehat{\Pi} f| \leq C \|f\|_{\mathcal{K}_{\sigma^n}^c}$ for some $C > 0$ independent of $n$, it follows immediately from (69) and (66) that

$$|A_n - A_{n-1}| \leq C \sigma^{n/2} (R_{v,n} + R_{w,n}) . \tag{72}$$

Next, since $\hat{r}_{n-1}(0, x) = 0$, there exists $C > 0$ such that

$$\|e^{\sigma^{-2n} \widehat{\mathcal{M}}_{c,n}(1-\sigma^2)} \widehat{\mathcal{L}} \hat{r}_{n-1}\|_{\mathcal{K}_{\sigma^n}^c} \leq C\sigma \|\hat{r}_{n-1}\|_{\mathcal{K}_{\sigma^{n-1}}^c} ,$$

see [Sch96, ES02]. This crucial estimate shows that, if $\sigma$ is sufficiently small, the linear semigroup *contracts* the remainder term $\hat{r}_n$, which is the reason for subtracting the Gaussian part from $\hat{v}_{c,n}$. Moreover, the last line in the right-hand side of (70) can be estimated by $C|A_n - A_{n-1}| + C\sigma^n R_{v,n}$, see [ES02]. Summarizing, we obtain

$$\rho_{r,n} \leq C\sigma \rho_{r,n-1} + C\sigma^{n/2}(R_{v,n} + R_{w,n}) \,. \tag{73}$$

On the other hand, it follows from Lemma 5.7 that

$$\|e^{\sigma^{-2n}(\Lambda_\beta + \nu)(1-\sigma^2)} w_{n-1}\|_{H^2} \leq C e^{-\nu_2 \sigma^{-2n}(1-\sigma^2)} \|w_{n-1}\|_{H^2} \leq C\sigma \|w_{n-1}\|_{H^2} \,,$$

hence
$$\rho_{w,n} \leq C\sigma \rho_{w,n-1} + (\mu + C_3 \sigma^{n/2}) R_{w,n} \,. \tag{74}$$

Finally,
$$|q_n - q_{n-1}| \leq C_3 \, e^{-\nu \sigma^{-2n+2}} R_{w,n} \,. \tag{75}$$

Now, if we assume that $C_3 \sigma_0 \leq \min(\mu, 1/8)$, it follows from (68) that

$$R_{w,n} \leq \frac{C_6}{1-2\mu} \rho_{w,n-1} \,, \quad R_{v,n} \leq 2C_5 \sigma^{-m} \rho_{v,n-1} + R_{w,n} \,. \tag{76}$$

Inserting these bounds into (72), (73), (74), (75), and using the fact that $2\mu C_6/(1-2\mu) \leq 1/2$ by (67), we obtain (71). $\square$

It is now straightforward to conclude the proof of Theorem 2.8. First, we choose $\sigma \in (0, \sigma_0)$ sufficiently small so that $C_8 \sigma \leq \sigma^{3/4}$. Then, we set $n_0 = 2m + 2$ and take initial data $\hat{v}_{c,n_0}, \hat{v}_{s,n_0}, w_{n_0}, q_{n_0}$ such that

$$|A_{n_0}| + \rho_{r,n_0} + \rho_{w,n_0} + \rho_{q,n_0} \leq \delta \,,$$

for some $\delta > 0$. Remark that, since the original problem is autonomous, we may without loss of generality take our initial data at time $t_0 = \sigma^{-2n_0}$.

As is easy to verify, if $\delta > 0$ is sufficiently small, the recursion relation (71) implies that $|A_n| + \rho_{r,n} + \rho_{w,n} + \rho_{q,n} \leq C\delta$ for all $n \geq n_0$. In particular, the assumptions of Proposition 5.8 are satisfied for all $n \geq n_0$. Moreover, there exists $C_9 > 0$, $A_* \in \mathbb{C}$, and $q_* \in \mathbb{R}$ such that, for all $n \geq n_0$,

$$|A_n - A_*| + \rho_{r,n} \leq C_9 \delta \sigma^{n/2} \,, \quad |q_n - q_*| \leq C_9 \delta e^{-\nu \sigma^{-2n}} \,, \quad \rho_{w,n} \leq \delta \,. \tag{77}$$

It remains to show that these estimates imply (9) and (10). First, combining (77) with (76) and using the last bound in (66), we obtain

$$R_{w,n} + e^{\nu \sigma^{-2n}} \sup_{\tau \in [\sigma^2, 1]} |q_n(\tau) - q_*| \leq C\delta \,, \quad n \geq n_0 \,.$$

If we undo the change of variables (50), we thus find

$$e^{\nu t}(\|W(\cdot, t)\|_{H^2} + |q(t) - q_*|) \leq C\delta \,, \quad t \geq t_0 \,,$$

which implies (10), see (37). Similarly, remarking that

$$\begin{aligned}
\hat{v}_{c,n}(\varkappa, x, \tau) - A_* \Psi(\varkappa\sqrt{\tau})\varphi(0, x) = \\
\hat{v}_{c,n}(\varkappa, x, \tau) - e^{\sigma^{-2n}\widehat{\mathcal{M}}_{c,n}(\tau - \sigma^2)}\Big(\hat{v}_{c,n-1}(\sigma\varkappa, x, 1) - \hat{r}_{n-1}(\sigma\varkappa, x)\Big) \\
+ A_{n-1}\Big(e^{\sigma^{-2n}\widehat{\mathcal{M}}_{c,n}(\tau - \sigma^2)}\Psi(\sigma\varkappa) - \Psi(\varkappa\sqrt{\tau})\Big)\varphi(\sigma^n\varkappa, x) \\
+ (A_{n-1} - A_*)\Psi(\varkappa\sqrt{\tau})\varphi(\sigma^n\varkappa, x) + A_*\Psi(\varkappa\sqrt{\tau})(\varphi(\sigma^n\varkappa, x) - \varphi(0, x)),
\end{aligned}$$

and using (77), (76), (66), (51) and Lemma 3.1, we arrive at

$$\sup_{\tau \in [\sigma^2, 1]} \Big(\|\hat{v}_{c,n}(\cdot, \cdot, \tau) - A_*\Psi(\cdot\sqrt{\tau})\varphi(0, \cdot)\|_{\mathcal{K}_{\sigma^n}^c} + \|\hat{v}_{s,n}(\cdot, \cdot, \tau)\|_{\mathcal{K}_{\sigma^n}^s}\Big) \leq C\delta\sigma^{n/2}, \quad n \geq n_0.$$

In particular, using (24), we have for all $n \geq n_0$,

$$\sup_{\tau \in [\sigma^2, 1]} \sup_{x \in [0, 2\pi]} \int_{-1/2\sigma^n}^{1/2\sigma^n} \Big(|\hat{v}_{c,n}(\varkappa, x, \tau) - A_*\Psi(\varkappa\sqrt{\tau})\varphi(0, x)| + |\hat{v}_{s,n}(\varkappa, x, \tau)|\Big)\, d\varkappa \leq C\delta\sigma^{n/2}.$$

If we now undo the changes of variables (50) and (46), we obtain

$$\sup_{x \in [0, 2\pi]} \int_{-1/2}^{1/2} |\hat{V}(\ell, x, t) - A_*\Psi(\ell\sqrt{t})\varphi(0, x)|\, d\ell \leq \frac{C\delta}{t^{3/4}}, \quad t \geq t_0. \tag{78}$$

Since $\int_{-1/2}^{1/2} e^{i\ell x}\hat{V}(\ell, x, t)\, d\ell = V(x, t)$ and

$$\int_{-1/2}^{1/2} e^{i\ell x} A_*\Psi(\ell\sqrt{t})\varphi(0, x)\, d\ell = \frac{A_*}{\sqrt{4\pi dt}} e^{-x^2/(4\pi dt)}\varphi(0, x) + \mathcal{O}(t^{-1}),$$

we see that (78) implies (9). The proof of Theorem 2.8 is now complete. □

# 6  Results in case II

In this section, we give the results about existence and stability of modulated fronts in case II, namely with $F(u) = 1 - u$ in (1). When $\alpha$ crosses zero, the homogeneous steady sate $U_+$ ahead of the front destabilizes and undergoes a Turing bifurcation, but the equilibrium $U_-$ behind the front remains asymptotically stable with some exponential rate. As explained in Remark 4.4, we have existence of a modulated front connecting the Turing pattern $U_{\text{per}}^+$ ahead of the front with the trivial solution $U_-$ behind.

**Theorem 6.1** *For $\varepsilon > 0$ sufficiently small, there exist a modulated front solution of (3) of the form*
$$U(x, t) = U_{\text{mf}}(x - ct, x), \quad x \in \mathbb{R}, \quad t \in \mathbb{R},$$
*where $U_{\text{mf}}(\xi, x)$ is $2\pi$-periodic in its second argument and $c = c_0 + \mathcal{O}(\varepsilon^2)$. Moreover, there exist positive constants $C, \beta_1, \beta_2$ (independent of $\varepsilon$) such that*

$$\sup_{\xi, x \in \mathbb{R}} |U_{\text{mf}}(\xi, x) - U_h(\xi)| \leq C\varepsilon,$$

*and*

$$\|U_{\mathrm{mf}}(\xi,\cdot) - U_{\mathrm{per}}^+(\cdot)\|_{(H^2(0,2\pi))^2} \le Ce^{-\beta_1\xi}, \quad \xi \ge 0,$$
$$\|U_{\mathrm{mf}}(\xi,\cdot) - U_-\|_{(H^2(0,2\pi))^2} \le Ce^{\beta_2\xi}, \quad \xi \le 0.$$

Remark that, in contrast with Theorem 2.6, the convergence rate of $U_{\mathrm{mf}}(\xi,x)$ towards $U_-$ as $\xi \to -\infty$ is independent of $\varepsilon$.

Proving the nonlinear stability of the modulated front is much easier here than in case I:

**Theorem 6.2** *If $\beta \in (0,\beta_0)$ and $\varepsilon > 0$ is sufficiently small, there exist positive constants $C$, $\nu$, $\delta$ such that the following holds. For all $V_0 : \mathbb{R} \to \mathbb{R}^2$ with $\|V_0(x)(1+e^{\beta x})\|_{H^2} \le \delta$, there exists a unique global solution $U(x,t)$ of (3) with initial data $U(x,0) = U_{\mathrm{mf}}(x,x) + V_0(x)$. Moreover, there exists a shift function $q : \mathbb{R}_+ \to \mathbb{R}$ and a real constant $q_*$ such that $U(x,t)$ can be represented as*

$$U(x,t) = U_{\mathrm{mf}}(x-ct-q(t),x) + V(x,t), \quad x \in \mathbb{R},\ t \ge 0,$$

*where*

$$\sup_{x\in\mathbb{R}} |V(x,t)| + \sup_{\xi\in\mathbb{R}} |V(\xi+ct,t)e^{\beta\xi}| + |q(t)-q_*| \le Ce^{-\nu t}, \quad t \ge 0. \tag{79}$$

**Proof.** We proceed exactly as in section 5.1. Setting

$$U(x,t) = U_{\mathrm{mf}}(x-ct-q(t),x) + V(x,t), \quad \text{and} \quad W(\xi,t) = V(\xi+ct,t)e^{\beta\xi},$$

we obtain equation (36) for $V$ and (39) for $W$ and $q$. However, we replace (40) with

$$\partial_t V = L(\partial_x)V + DN(U_-)V + N_1(U_-,V) + \dot q\, \partial_1 U_{\mathrm{mf}} + G(U_{\mathrm{mf}} - U_-, V).$$

In contrast with the previous case, the spectrum of the linear operator $\mathcal{L} = L(\partial_x) + DN(U_-)$ is strictly contained in the left-half plane. In addition, if we assume that $0 < \beta < \beta_2$ (where $\beta_2$ is defined in Theorem 6.1), we have

$$\sup_{\xi\in\mathbb{R}} \|U_{\mathrm{mf}}(\xi-q(t),\cdot) - U_-\|_{H^2(0,2\pi)}\, e^{-\beta\xi} \le C.$$

Thus, proceeding as in (41), we find

$$\|G(U_{\mathrm{mf}} - U_-, V)\|_{H^2} \le C\|W\|_{H^2}(1 + \|V\|_{H^2}),$$

for some $C > 0$. Summarizing, the evolution system for $V,W,q$ has the form

$$\begin{aligned}
\partial_t V &= \mathcal{L}V + \mathcal{O}(V^2) + \mathcal{O}(\dot q) + \mathcal{O}(W), \\
\partial_t W &= \Lambda_\beta W + \mathcal{O}((\varepsilon+|q|)W) + \mathcal{O}(VW), \\
\dot q &= \mathcal{O}((\varepsilon+|q|)W) + \mathcal{O}(VW).
\end{aligned} \tag{80}$$

By construction, there exists $\nu > 0$ such that $\|e^{\mathcal{L}t}V\|_{H^2} \le Ce^{-\nu t}\|V\|_{H^2}$ and $\|e^{\Lambda_\beta t}W\|_{H^2} \le Ce^{-\nu t}\|W\|_{H^2}$ for all $t \ge 0$. Therefore, if $\varepsilon > 0$ is small enough and if the initial data satisfy $\|V(\cdot,0)\|_{H^2} + \|W(\cdot,0)\|_{H^2} + |q(0)| \le \delta$ for some sufficiently small $\delta > 0$, it is clear that (80) has a global solution in $H^2 \times H^2 \times \mathbb{R}$. Moreover, there exists $q_* \in \mathbb{R}$ such that

$$\|V(\cdot,t)\|_{H^2} + \|W(\cdot,t)\|_{H^2} + |q(t)-q_*| \le Ce^{-\nu t},$$

as $t \to +\infty$. This concludes the proof. $\square$

# 7 Numerical simulations

Using numerical simulations, we illustrate the results from section 2 concerning modulated fronts, and, for a different model, the existence and stability of modulated pulses. These computer experiments give the impression that the assertions of Theorem 2.6 and 2.8 are also true for $\varepsilon$ and $\delta$ not necessarily small, i.e. the existence and stability of modulated structures also holds for non-small values of the bifurcation parameter and non-small perturbations.

## 7.1 Modulated fronts

To solve (5) numerically we subtract the (unmodulated) front $h(y) = \tanh(y/2)$ from $u$, i.e, we set $u = h_c + w$ and integrate the resulting system for $(w, v)$ using finite differences, periodic boundary conditions on the large domain $y \in (-60\pi, 60\pi)$ (see Remark 7.1 below), and implicit time stepping. The parameters $\gamma = 0.5$ and $c_0 = 0.5$ are kept fixed. We start with rather generic initial conditions and calculate the (discretized) modulated fronts dynamically which of course is only possible if they are stable (for the discretized system).

Figure 4 shows the evolution towards a modulated front for $\alpha = 0.1$ and initial conditions $(w, v)|_{t=0} = (1/\cosh(y), 0.01 \sin(y))$. The hump in $w|_{t=0}$ is transformed into a shift in $u = h + w$ rather quick. The transient time in which essentially $v(t)$ reaches its proper amplitude $\mathcal{O}(\varepsilon)$ and then couples back into the $u$ equation to produce the modulating pulse is about 50 units. Convergence of the solution to the modulated front with similar transient behavior was observed for more rather generic initial conditions. Starting with initial conditions $(w, v)|_{t=0} = (1/\cosh(y), 0.3 \sin(y))$ we get a much shorter transient. We can also, for instance, add $\mathcal{O}(1)$ humps away from $y = 0$ to $(w, v)|_{t=0}$. These get damped out very quickly.

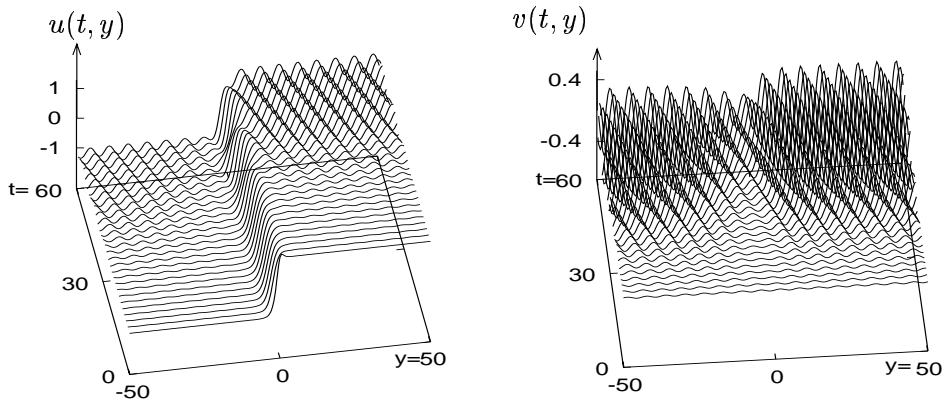

Figure 4: Evolution towards a modulated front, $\alpha = 0.1$.

Figure 5 shows snapshots of the solutions at some fixed times $t$. In a) ($t = 80$ with the solution from fig.4) with $\alpha = 0.1$ the different amplitudes of the periodic patterns at $y = \pm 50$ are clearly visible. In order to display the different decay rates to the periodic patterns ahead of and behind the front we take a smaller $\alpha = 0.01$; b) shows the modulated front, and c),d) the functions $w = u_{\text{mf}} - \tanh(y/2)$ and $v$. The effect of $-\gamma v F(h_c + w)$ is that the Turing pattern in $v$ gets damped while passing through the modulated front. It then converges with rate $\mathcal{O}(e^{-\mathcal{O}(\varepsilon)|y|})$ to the Turing pattern in the recovery zone behind the front.

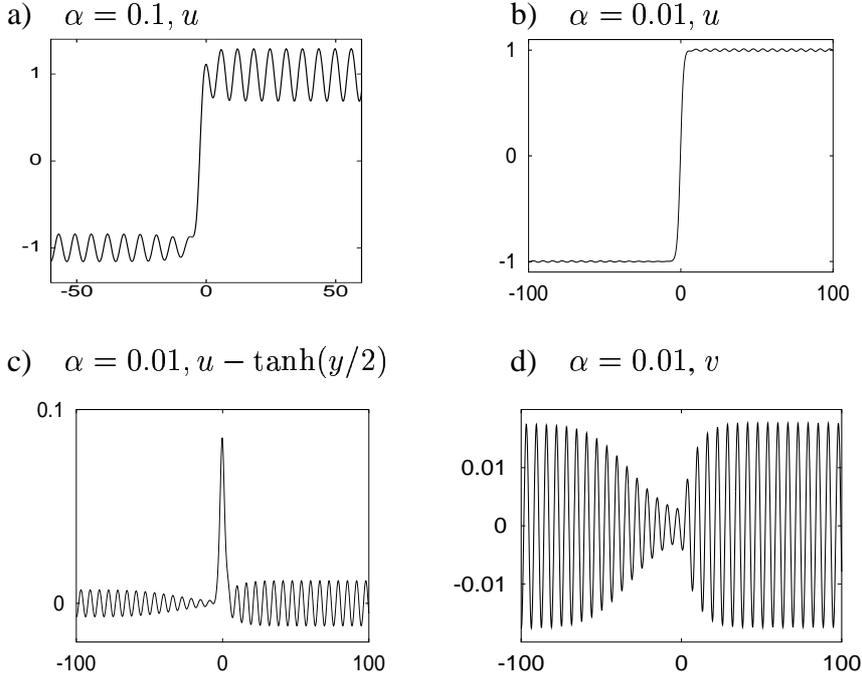

a) $\alpha = 0.1, u$
b) $\alpha = 0.01, u$
c) $\alpha = 0.01, u - \tanh(y/2)$
d) $\alpha = 0.01, v$

Figure 5: Snapshots of modulated fronts.

**Remark 7.1** *There is a conceptual problem with periodic boundary conditions, in particular for very small $\varepsilon > 0$. Strictly speaking, on any finite domain with periodic boundary conditions there is no "ahead and behind the front". To minimize this effect we chose a large domain, and see that at the center of mass, say $y \in (-50, 50)$ the analytically predicted dynamics of the modulated front are nicely recovered. Using even larger domains it can be checked that the influence of the boundary conditions near the center is indeed very small.*

Finally, we illustrate what happens in case III with $F(u) = 1 + u$. As an example of the typical evolution of the unstable front, fig.6 shows $u(t, y) - \tanh(y/2)$, i.e., the $u$ component of the perturbation of the front in the comoving frame, with $\alpha = 0.01$.

## 7.2 Modulated pulses

In order to obtain an example with modulated pulses we couple a non symmetric complex Ginzburg-Landau equation (nsGLe) with a Swift-Hohenberg like equation. The nsGLe for the complex field $A(x,t) \in \mathbb{C}$ reads

$$\partial_t A = c_1 \partial_x^2 A - (\alpha_0 + i\nu_0)A + \alpha_1 \overline{A} + 4i|A|^2 A \tag{81}$$

where $c_1 = c_{1r} + ic_{1i} \in \mathbb{C}$ and $\alpha_0, \alpha_1, \nu_0 \in \mathbb{R}$ are parameters with $\nu_0 > \alpha_1 > 0$ and in particular $0 < c_{1r}, \alpha_0, \alpha_1 \ll 1$. Thus, the nsGLe is a dissipative perturbation of the nonlinear

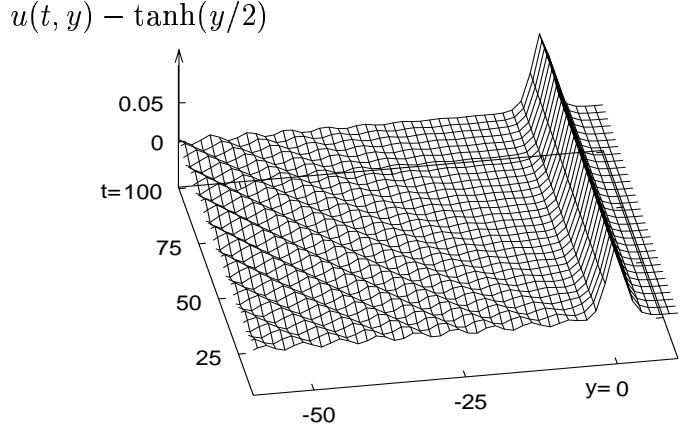

Figure 6: Typical evolution in case III; the homogenous rest state behind the front is unstable, and a Turing pattern develops behind the front. However, the front travels away from the pattern, and the resulting solution is not time periodic in any comoving frame.

Schrödinger equation, where due to the term $\alpha_1 \overline{A}$ the usual $\mathcal{S}^1$-symmetry $A \mapsto e^{i\vartheta} A$ is broken. It arises for instance as a modulation equation for optical fibers with phase-sensitive amplifiers, see [KK96] and the references therein, or for dissipative systems with a resonant spatially periodic forcing [Uec01]. For suitable parameters the nsGLe has an exponentially stable one parameter family $\{A_{\mathrm{pu}}(\cdot - x_0) : x_0 \in \mathbb{R}\}$ of pulse solutions. See [KS98] for this result and a comprehensive discussion of the nsGLe. For $c_{1r} = 0$ the pulse is explicitly given by

$$A_{\mathrm{pu}}(x) = \sqrt{b_1} \mathrm{sech}(\sqrt{b_2} x) e^{i\vartheta}, \\ \cos(2\vartheta) = \alpha_0/\alpha_1, \quad b_2 = (\nu_0 + \alpha_1 \sin(2\vartheta))/c_{1i}, \quad b_1 = c_{1i} b_2/2. \tag{82}$$

For small $c_{1r} > 0$ this pulse persists, and the spectrum of the linearization of (81) about $A_{\mathrm{pu}}$ is as follows. The continuous spectrum is given by the two curves

$$\lambda_{1,2}(k) = -\alpha_0 - c_{1r} k^2 \pm \mathrm{i} \sqrt{(c_{1i} k^2 + \nu_0)^2 - \alpha_1^2}.$$

Moreover, we obtain one simple eigenvalue 0 from the translational invariance and the rest of the discrete spectrum, consisting of 5 more simple eigenvalues, is in the left complex half plane, see [KS98].

We now couple the nsGLe for $A = u_1 + \mathrm{i} u_2$ with the SHe for $v \in \mathbb{R}$, i.e. we consider

$$\partial_t U = LU + N(U) \tag{83}$$

where $U = (u_1(x,t), u_2(x,t), v(x,t)) \in \mathbb{R}^3$,

$$L = \begin{pmatrix} \alpha_1 - \alpha_0 + c_{1r} \partial_x^2 - c_0 \partial_x & -c_{1i} \partial_x^2 + \nu_0 & \mu_1 \\ c_{1i} \partial_x^2 - \nu_0 & -(\alpha_1 + \alpha_0) + c_{1r} \partial_x^2 - c_0 \partial_x & 0 \\ 0 & 0 & -(1 + \partial_x^2)^2 + \alpha_2 \end{pmatrix},$$

$$N(U) = \left( -4|u|^2 u_2, \, 4|u|^2 u_1, \, -v^3 + \mu_2 u_1 v \right)^T, \quad |u|^2 = u_1^2 + u_2^2,$$

and where $\alpha_0, \alpha_1, \nu_0$ and $c_1$ are fixed in such a way that the nsGLe has a stable pulse solution. We use $\mu_{1,2} \in \mathbb{R}$ as coupling parameters and $\alpha_2$ as the bifurcation parameter. Moreover we set $h(\xi) = (u_1(\xi), u_2(\xi)) \in \mathbb{R}^2$, $\xi = x - c_0 t$, where $u_1(x) + iu_2(x) = A_{\text{pu}}(x)$ is the pulse solution of the nsGLe.

Using the spectral properties of $A_{\text{pu}}$ and setting $\mu_2 = 0$ it is clear that the analogue of Theorem 2.1 holds for (83), i.e., that for $\alpha_2 < 0$ the pulse $(h, 0)$ is exponentially stable. Moreover, numerical simulations of the nsGLe reveal that also for small $c_{1r} > 0$ the pulse $A_{\text{pu}}$ fulfills $\text{Re} A_{\text{pu}}(x) = u_1(x) \geq 0$ for all $x \in \mathbb{R}$. Using this, we can conclude as in the proof of Theorem 2.1 that for $\alpha_2 < 0$ the linearization of (83) around $(h, 0)$ has spectrum in the left complex half-plane for all $\mu_2 \leq 0$, and in fact even for $0 < \mu_2 < \mu_*$ for a sufficiently small $\mu_* = \mu_*(\alpha_2) > 0$. Hence we have roughly the same starting point for a bifurcation analysis for (83) as for (1), and may expect modulated pulses to bifurcate for $\alpha_2 > 0$.

This is now illustrated by numerical simulations of (83), where we fix $(\alpha_0, \alpha_1, \nu_0, c_1) = (0.6, 0.8, 1, 0.1 + 10i)$, $c_0 = 1$ and $\mu_1 = 1$, and integrate (83) in the moving frame $y = x - c_0 t$, again on the large domain $\xi \in (-60\pi, 60\pi)$ with periodic boundary conditions. In order to obtain nice graphs, the value $c_{1i} = 10$ has been chosen relatively large so that the pulse has a width larger than the period $2\pi$ of periodic pattern.

a) $\mu_2 = -1$ \hspace{4cm} b) $\mu_2 = 1$

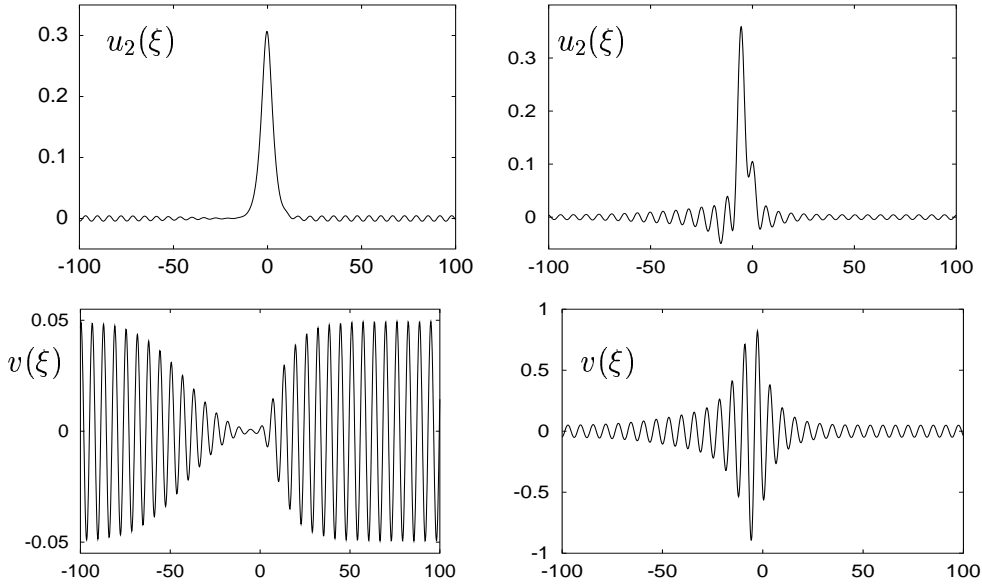

Figure 7: $u_2$ and $v$ for a modulated pulse, $\varepsilon = 0.05$. a) $\mu_2 = -1$, b) $\mu_2 = 1$

In the first simulation we let $\alpha_2 = \varepsilon^2 = 0.0025$ and $\mu_2 = -1$. Again we choose a small $\varepsilon = 0.05$ in order to resolve the different convergence rates ahead and behind the pulse. We start with an approximation of $U_{\text{mp}}$ in the form $u_1(\xi, 0) + iu_2(\xi, 0) = A_{\text{pu}}(\xi)$ with $A_{\text{pu}}$ from (82) and $v(\xi) = \varepsilon \cos(\xi)$. The solution converges quickly to a modulated pulse $U_{\text{mp}}$, which illustrates the stability of the modulated pulses. In figure 7a) we show $u_2(t = 50)$ and $v(t = 50)$ for $\xi \in (-100, 100)$, i.e., roughly half the computational domain; cf. Remark 7.1. The first component $u_1$ is similar to $u_2$ but $\max u_1(\xi, t) \approx 1$ so that $u_2$ is more eligible for graphical purposes. The effect of $\mu_2 = -1$ is that the Turing pattern gets damped while passing through the pulse, with an $\mathcal{O}(e^{-\mathcal{O}(\varepsilon)|y|})$ convergence to the Turing pattern behind the pulse.

For $\mu_2 = 1$ we get the converse effect, see figure 7b). The pattern gets amplified while passing through the pulse and decays to amplitude $\mathcal{O}(\varepsilon)$ in the recovery zone behind the pulse. Note however that this is now also a nonlinear effect of the damping $-v^3$ in the $v$ equation.

Finally, in figure 8 we present snapshots of $u_1$ and $v$ from a numerical simulation with $\varepsilon=0.9$, $\mu_2=-1$ and the remaining parameters as above. The convergence to the periodic pattern is equally fast ahead and behind the pulse. We remark that numerically we could produce stable modulated pulses $U_{\text{mp}}^\varepsilon$ up to $\varepsilon \approx 2$. This works most easily by continuation in $\varepsilon$, i.e., by slowly increasing $\varepsilon$, integrating and waiting until the solution settles to $U_{\text{mp}}^\varepsilon$, then increasing $\varepsilon$ again. For $\varepsilon > 2$ this breaks down and the solution disintegrates into wave–packets.

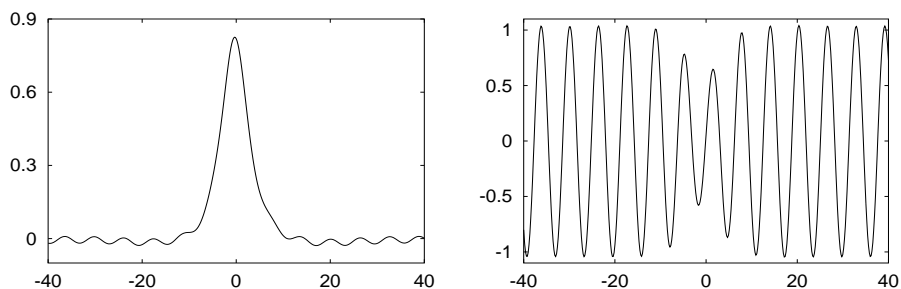

Figure 8: $u_1$ and $v$ for a modulated pulse, $\mu_2 = -1$, $\varepsilon = 0.9$